\documentclass[aps,preprint,amssymb,amsmath,eqsecnum,nofootinbib]{revtex4}
\usepackage{graphics}
\newtheorem{question}{Question}[section]
\newtheorem{definition}{DEFINITION}[section]
\newtheorem{theorem}{Theorem}[section]
\newtheorem{lemma}{Lemma}[section]
\newtheorem{answer}{Answer}[section]
\newtheorem{corollary}{Corollary}[section]

\newtheorem{remark}{Remark}[section]
\newtheorem{axiom}{AXIOM}[section]

\newenvironment{hypothesis}{HP: \begin{center}} {\end{center}}
\newenvironment{thesis}{TH: \begin{center}} {\end{center}}
\newenvironment{proof}{\begin{center}PROOF: \end{center}} {$ \blacksquare $}
\newtheorem{example}{Example}[section]
\begin{document}
\title{On the mathematical structure of Tonal Harmony}
\author{Gavriel Segre}
\homepage{http://www.gavrielsegre.com}
\email{info@gavrielsegre.com}
\date{15-2-2004  16:07}
\maketitle
\newpage
\tableofcontents
\newpage
\section{Acknowledgments}
I would like to thank Enrico Fubini for stimulating discussions on
the Jewish nature of Sch\"{o}nberg's way of relationing with the
tradition, i.e. in recognizing that any true evolution or even
overcoming of a tradition requires its deep knowledge.

I would then like to thank Gilberto Bosco for many stimulating
discussions about the meaning of mathematically formalizing (not
necessary classical) Harmony

I would then like to thank Francesco Spagnolo for many useful
remarks concerning Jewish Music and for allowing me to access to
the Grove's thesaurus

I would then like to thank my father for uncountably many fruitful
discussions concerning the themes discussed in this paper

A final remark, with a retrospective value, concerns my gratitude
to Massimo De Benedetti whose skills on using Bibtex were
unvaluable

This work is dedicated to Euler \cite{Euler-26a},
\cite{Euler-26b}, \cite{Euler-26c} and Lagrange
\cite{Lagrange-67a}, \cite{Lagrange-67b}, the fathers of
Mathematical Physics of Music

\section{Notation}
\begin{tabular}{|c|c|}
  $ p_{Nature} (t) $ & natural pound at time t \\
  $ p_{Culture} (t) $ & cultural pound at time t  \\
  $ {\mathcal{R}}_{ac} [i(t)] $ & acoustical response to the input i(t) \\
  $ harmonic[\nu , n ] $  & $n^{th}$ harmonic of the note $ \nu $ \\
  $ {\mathbf{mi}} ( \nu_{1} \, , \,  \nu_{2} ) $ & musical
  interval among $ \nu_{1} $ and  $ \nu_{2} $ \\
   $ \nu_{ref} $ & reference note \\
   $ {\mathbf{p}}( \nu ) $ & pitch of the note $ \nu $ \\
  $ S(\nu) $ & Fourier transform of the signal $ s(t) $ \\
  $ {\mathcal{D}} $ & tranfer function of the human ear \\
  $ scale-range ( \nu ) $  & scale-range of the note $ \nu $ \\
  $ R_{\nu} $ & rescaling function to the scale-range of $\nu$ \\
  $ a \, = \, b \, mod c^{{\mathbb{Z}}} $  & a is congruent to b modulo 
powers of c  \\
  $ [ a ]_{c^{{\mathbb{Z}}}} $ & residue class of a modulo powers of  c \\
  $ {\mathbb{Z}}_{c^{{\mathbb{Z}}} } $  &  set of the residue classes of a 
modulo powers of  c   \\
  $ c_{n} $ & n level Euler coordination \\
  $ {\mathcal{N}}_{Euler} $ & Euler's notes \\
  $ {\mathcal{P}}_{Euler} ( \nu)  $ & Euler's point of $ \nu $  \\
  $ {\mathcal{S}}_{Euler} $ & Euler's space \\
  $ H_{prime} $ & prime vector of extended Euler's space \\
$ {\mathcal{N}}_{Euler}^{just-tuned} $ & just-tuned Euler's notes \\
  $ {\mathcal{S}}_{Euler}^{just-tuned} $  & just-tuned subset of Euler's 
space \\
  $ {\mathcal{N}}_{Euler}^{Pyt-tuned} $ & pytagorically-tuned Euler's notes  
\\
  $ {\mathcal{S}}_{Euler}^{Pyt-tuned} $  &  pytagorically-tuned subset of 
Euler's space \\
  $ {\mathcal{N}}_{Euler}^{n-temp-tuned} $ & n-tempered-tuned Euler's notes 
\\
$ {\mathcal{S}}_{Euler}^{n-temp-tuned} $ & n-tempered-tuned subset of 
Euler's space \\
  $ {\mathcal{N}}_{Euler}^{n_{1},n_{2},n_{3}-temp-tuned} $ & $ 
n_{1},n_{2},n_{3} $ -tempered-tuned Euler's notes \\
$ {\mathcal{S}}_{Euler}^{n_{1},n_{2},n_{3}-temp-tuned} $ &  $ 
n_{1},n_{2},n_{3} $- tempered-tuned subset of Euler's space \\
$ {\mathbb{E}}_{Euler}^{notes} $ & canonical notes' basis of
Euler's space \\
$ {\mathbb{E}}_{Euler}^{int} $ & canonical intervals' basis of
Euler's space \\
\hline
\end{tabular}
\newpage
\begin{tabular}{|c|c|}
  $ \hat{Kf} $ & fifth comma (pytagorical comma) \\
$ \hat{Kt} $ & third comma (syntonic comma)  \\
$ {\mathbf{Kf}} $ & pitch of the fifth comma \\
$ {\mathbf{Kt}} $ & pitch of the third comma \\
$ G_{Euler} ( \nu_{1} \, , \, \nu_{2} ) $ & Euler's gradus
suavitatis of  $ \nu_{1} \, , \, \nu_{2} $ \\
$ I[ s_{1} \, \cdots \, s_{n} ] $ & index of physical consonance
of $ s_{1} \, , \cdots \, , \, s_{n} $ \\
$ \omega_{1} \sim_{{\mathbb{Q}}} \,  \cdots \, \sim_{{\mathbb{Q}}}
\omega_{n} $ &
commensurability of $ \omega_{1} \, , \cdots \, , \, \omega_{n} $ \\
$ I[ \omega_{1} \, \cdots \, \omega_{n} \, | \bar{a} ] $ & index
of physical  consonance among  $ \omega_{1} \, ,
\cdots \, , \, \omega_{n} $ w.r.t. $ \bar{a} $ \\
$ scale( \omega \, , \, R ) $ & scale of $ \omega $ at fixed
interval R \\
$ ord( \{ a_{n} \}) $ & ordered permutation of $  \{ a_{n} \} $  \\
$ \downarrow t $ & downward closure of the tree t \\
$ TC(s) $ & transitive closure of the set s \\
support(s) & support of the set s \\
$ V_{afa} $ & class of the pure sets \\
$ G $ & nodes' set of the graph  $ {\mathbf{G}} $ \\
$ \rightarrow_{G} $ & edges' binary relation of the graph $
{\mathbf{G}} $ \\
$ GRAPHS $ &  proper class all  graphs \\
$ D ({\mathbf{G}}) $ & decoration of the graph $ {\mathbf{G}} $ \\
$ \Sigma^{n} $ & words of length n on the alphabet $ \Sigma $ \\
$ \Sigma^{\star} $ & words on the alphabet $ \Sigma $  \\
$ \Sigma^{\infty} $ &  sequences on the alphabet $ \Sigma $ \\
$ \Sigma_{NR}^{n} $ & nonrepetitive  words of length n on the alphabet $ 
\Sigma $ \\
$ \Sigma_{NR}^{\star} $ & nonrepetitive words on the alphabet $ \Sigma $  \\
  $ \vec{x} $ & word \\
  $ \lambda $ & empty word \\
  $ | \vec{x} | $ &  length of the word $ \vec{x} $ \\
  $ string(n) $ & $ n^{th} $ word in quasilexicographic order \\
  $ |n | $ & length of the $ n^{th} $ word in quasilexicographic order \\
  $ \bar{x} $ & sequence \\
  $ <_{p} $ & prefix order relation \\
  $ \cdot $ & concatenation operator \\  \hline
\end{tabular}
\newpage
\begin{tabular}{|c|c|}
$ x_{n} $ & $ n^{th} $ digit of the word $\vec{x} $ or of the sequence $ 
\bar{x} $ \\
  $ \vec{x}(n) $ & prefix of length n of the word $ \vec{x} $ or of the 
sequence $ \bar{x} $  \\
  $ \vec{x}(n,m) $ & subword of the sequence $ \bar{x}
  $ obtained taking the digits from the $n^{th}$ to the $m^{th}$ \\
  $ \vec{x}^{n} $  &  word made of n repetitions
of the word  $ \vec{x} $ \\
$ [i]_{12} $ & $ i^{th} $ letter of the musical alphabet $
{\mathbb{Z}}_{12} $ \\
$ C \, , \, \cdots \, , \, B $ & musical notation for, respectively, $
[0]_{12} \, , \,  \cdots \, , \,  [11]_{12}$ \\
$ T_{y} $ & translation operator by y \\
$ Inv $ & inversion operator \\
$ mode( \vec{x} \, , \, i ) $ & mode of the word $ \vec{x} $ of $
i^{th} $ degree \\
$ chord ( \vec{x} \, , \, i \, , \, n ) $ & chord of $ \vec{x} $
of $ i^{th} $ degree at level n \\
$ maxlevel ( \vec{x} ) $ & maximum possible level of $ \vec{x}
$'s chords \\
$ \hat{g} $ & map on words induced by the map g on letters \\
$ \sim_{T} $ & translational equivalence relation on words \\
$ \sim_{Inv} $ & inversion's equivalence relation on words \\
$ d( x,y ) $ & distance among x and y \\
$ I( \vec{x} ) $ &  interval vector of $ \vec{x} $ \\
$ {\mathcal{S}}_{greg} $ & gregorian words \\
$ {\mathcal{S}}_{cl} $ & classical words \\
$ {\mathcal{T}}_{Maz} $  & Mazzola tonalities \\
$ tonality( \vec{x} \, , \, n) $ & tonality of the word $ \vec{x}
$ at level n \\
$ {\mathcal{T}}_{n} $ & tonalities at level n \\
$ {\mathcal{T}} $ & tonalities \\
$ {\mathcal{HW}} (t) $ &  harmonic words of t \\
${\mathcal{P}} ( t_{1} , t_{2} ) $ & pivotal degrees of $ t_{1} , t_{2} $ \\
$ { \mathcal{C}}(t \, , \, {\mathcal{T}}_{context} ) $ & cadences
of t w.r.t. to the context $ {\mathcal{T}}_{context}  $ \\
$  {\mathcal{T}}_{n.c.} ( \vec{x} , n ) $ & natural context for $
tonality(\vec{x} , n) $ \\
$ { \mathcal{MC}}(t \, , \, {\mathcal{T}}_{context} ) $ & minimal
cadences of t w.r.t. to the context $ {\mathcal{T}}_{context}  $ \\
\hline
\end{tabular}
\newpage
\begin{tabular}{|c|c|}
$ SYM( \vec{x} ) $ & symmetry group of the word $ \vec{x} $ \\
$ {\mathcal{M}}_{Maz} ( t_{1} \, , \, t_{2} ) $ & Mazzola's
modulations from $ t_{1} $ to $ t_{2} $ \\
$ {\mathcal{M}} ( t_{1} \, , \, t_{2} ) $ &
modulations from $ t_{1} $ to $ t_{2} $ \\
$ {\mathcal{MP}}_{Maz} $ & Mazzola's tonal musical pieces \\
$ {\mathcal{MP}} $ &  tonal musical pieces \\
$ \Xi $ & pytagoric musical alphabet \\
$ [i]_{12}^{(n)} $ & $ i^{th} $ letter at the $ n^{th} $  cycle of
$ \Xi $ \\
$ C^{(n)} \, , \, \cdots \, , \,  B^{(n)} $ & musical notation
for, respectively, $
[0]_{12}^{(n)} \, , \, \cdots \, , \, [11]_{12}^{(n)} $ \\
$ [i]_{12}^{(n) \, just} $ & $ i^{th} $ letter at the $ n^{th} $  cycle of
of the just-intonation alphabet  \\
$ C^{(n)}_{just} \, , \,  \cdots \, , \,  B^{(n)}_{just} $ &
musical notation for, respectively, $
[0]_{12}^{(n) \, just} \,  , \, \cdots \, , \,  [11]_{12}^{(n) \, just} $ \\
$ T_{z}^{Pyt} $ & pytagoric translation by z \\
$ Inv^{Pyt} $ & pytagoric inversion \\
$ C_{+} $ & fifth cycle's raising \\
$ C_{-} $ & fifth cycle's lowering \\
$ \sim_{T^{Pyt}} $ & pytagoric translation's equivalence relation
\\
$ \tilde{x} $ & pytagoric word \\
$ \tilde{pw} ( \vec{x} ,n ) $ & pytagoric word of $ \vec{x} $ at
cycle n \\
$ {\mathcal{T}}^{Pyt} $ & pytagoric tonalities \\
$ {\mathcal{T}}_{n} $ & pytagoric tonalities at level n \\
$ {\mathcal{HW}}^{Pyt} (t) $ &  pytagoric harmonic words of t \\
$ {\mathcal{P}}^{Pyt} ( t_{1} , t_{2} ) $ & pytagoric pivotal
degrees of $ t_{1} , t_{2} $ \\
$ { \mathcal{C}}^{Pyt}(t \, , \, {\mathcal{T}}^{Pyt}_{context} ) $
& pytagoric cadences
of t w.r.t. to the  context $ {\mathcal{T}}^{Pyt}_{context}  $ \\
$ { \mathcal{MC}}^{Pyt}(t \, , \, {\mathcal{T}}^{Pyt}_{context} )
$ & pytagoric minimal cadences of t w.r.t. to the context $ 
{\mathcal{T}}^{Pyt}_{context}  $ \\
$ {\mathcal{M}}^{Pyt}  ( t_{1} \, , \, t_{2} ) $ &
pytagoric modulations from $ t_{1} $ to $ t_{2} $ \\
$ {\mathcal{M}}_{K}^{Pyt}  ( t_{1} \, , \, t_{2} ) $ & pytagoric
comma-modulations from $ t_{1} $ to $ t_{2} $ \\
$ {\mathcal{MP}}^{Pyt} $ & pytagoric tonal musical pieces \\
$ t_{1} \, \sim_{K} \,  t_{2} $ & $ t_{1} $ is a
comma-displacement of $ t_{2} $ \\
$ s_{1} \, \sim_{P} \, s_{2} $ & $  s_{1} $ is P-metamere to $
s_{2} $ \\
$ \sim_{P} \, s $ & P-valence of s \\ \hline
\end{tabular}
\newpage
\section{Introduction} \label{sec:Introduction}
No example could be given of the radical dichotomy existing
nowadays among Science and Humanities as the intellectual analysis
on the structure of Contemporary Music.

Let us start analyzing the overwhelming confusion existing in most
of the discussions concerning the concept of \emph{musical
consonance}.

The correct conceptual approach would consist in distinguishing
among the following two notions:
\begin{itemize}
    \item \emph{physical consonance}: a known phenomenon in
    Acoustics concerning a particular mathematical structure of
    some patterns of sounds played simultaneously
    \item \emph{esthetic consonance}:  the set of esthetic rules
    codified by a given Harmomy, i.e. by a given human
    formalization of the admissible patterns of simultaneous
    sounds
\end{itemize}
As a matter of principle, no  reason exist in assuming that
\emph{esthetic consonance} has to be constrained by \emph{physical
consonance}:

is the natural phenomenon of animals tearing to pieces themselves
in the jungle an index of the fact that we should make the same ?

Indeed a great peculiarity of the human specie consists in that it
has developed \emph{Culture}: a certain (continously evolving in
time) amount of information.

Our behavior is, conseguentially, determined nowadays both by
Nature and by Culture,  with  relative proportions $ p_{Nature} $
and $ p_{Culture} $ evolving with time as:
\begin{eqnarray} \label{eq:Nature versus Culture}
  p_{Nature} (t) \, + \, p_{Culture} (t) \;  & =& \; 1 \\
  \frac{d p_{Nature} }{d t } \; &  < &  \; 0 \\
  \frac{d p_{Culture} }{d t } \; &  > &  \;  0
\end{eqnarray}
The estimation of the contributions of $ p_{Nature} $ and $
p_{Culture} $ is, typically, a matter of great discussion among
researchers, each one tending to over-estimate the role of its own
research area.

It should be said, with this regard, that the colossal amount of
sloveliness shown by "human scientists" whenever they try to adopt
precise, rigorously defined scientific notions \cite{AAVV-00},
\cite{Sokal-Bricmont-99} streghtens scientists' common opinion
that nothing serious is going on in their exponential productions
of words without an underlying rigorous mathematical language.

This state of affairs occurs, first of all, in all the analyses by
Human Scientists concerning the old question: " is music
\footnote{The idle issue concerning the definition of the term
music has led to curious theoretizations (we could call them
\emph{mafious esthetics})  such as the one expressed in
\cite{Nattiez-02a} according to which the esthetical value of a
work of art depends on whether the artist has paid the required
percentage to the suitable corporation that conseguentially
attributed to him the social status of artist,this fact
determining, for example, that $ 4'33''$ by John Cage, whose score
is : " I TACET ; II TACET ; III TACET" is an avant-garde
masterpiece instead of being a fraud; according to these curious
esthetics Bach, Mozart and Beethoven weren't great composer owing
to the fact they produced great music: they composed great music,
by definition, owing to the fact that they were composers} a
language "?

Let us observe, with this regard, that a misconception of the
mathematical notion of language, i.e. the semantical side of
Mathematical-Logic based on the notion of model of a formal system
about which I demand to the Appendix G.5 "Formal Logic" of
\cite{Mazzola-02a} or to the $ 10^{th} $ chapter "Geometric Logic
and Classifying Topoi" of \cite{Mc-Lane-Moerdijk-92}, has led some
musicologists (cfr. e.g. the second chapter "L'occidente cristiano
e l'idea di musica" of the first part "I problemi estetici e
storici della musica" of \cite{Fubini-95}) to see a mutual
exlusivity among  conceptions based on the link
\emph{music-poetry} (with an emphasis to the emotional valence of
this latter) and conceptions based on the link
\emph{music-mathematics} (with an emphasis to formal aspects);
since the semantical side of Mathematics allows, as a matter of
principle, to encode both \emph{intrinsic musical meanings} and
\emph{extrinsic musical meanings}  \cite{Nattiez-02b}.
such a misconception lacks of any conceptual ground.

An ever more catastrophic situation occurs, anyway, as to the application
of \emph{Information Theory} to \emph{Esthetics}:

that the \emph{esthetic level} of an art of work is someway
related to its \emph{amount of information} is something very
intuitive; the lack of mathematical rigor in many attempt to apply
Information Theory to Esthetics \cite{Arnheim-74} \footnote{and
more generally the plethora of nonsense written by "social
scientists"  on speaking about information and codes (cfr. e.g.
the cap.1 "L'universo dei segnali" and the cap.2 "L'universo del
senso" of \cite{Eco-68} or the cap.5 "La famiglia dei codici" of
\cite{Eco-84} based on the ravings about codes and Eco's s-codes
of \cite{Eco-75})} kept the scientific community away from the
whole matter.

This is a pity since, whenever applied in a serious way \footnote{
A particular analysis should be made as to the section 3.10
"Peculiarit\'{a} e struttura statistica del messaggio musicale" of
\cite{Angeleri-99} where an inconsistent dichotomy among the
\emph{techno-semantic component} and the \emph{esthetic component}
of information is introduced}, Information Theory could allow to
get some insight, e.g. through the following footsteps:
\begin{enumerate}
    \item \emph{Shannon's probabilistic information} has to be adopted to
    quantify the information's amount of an art of work relative to
    some cultural data, its absolute, intrinsic information's amount being 
quantified by its  \emph{Kolmogorov-Solomonoff-Chaitin's
    algorithmic information} \cite{Calude-02}; for example the cultural data 
of
    Tonal Harmony may be codified as a probability distribution on
    the space of monophonic musical pieces w.r.t. which the
    pieces containing chromatic notes not belonging to the
    diatonic scale of the underlying tonality are strongly less
    probable than pieces not containing alterations and have,
    conseguentially an higher \emph{probabilistic information}
     relative to Tonal Harmony though their absolute, \emph{algorithmic
    information} is not higher \footnote{A good example of such a kind of 
information-theoretical esthetical estimation
    is given by Gerard Hoffnung's not too flattering cartoon about Webern's 
music\cite{Hoffnung-00}  }.
    \item the esthetical value of an art of work, seen as a
    function of its amount of information I, has the
    following qualitative behavior:

\includegraphics{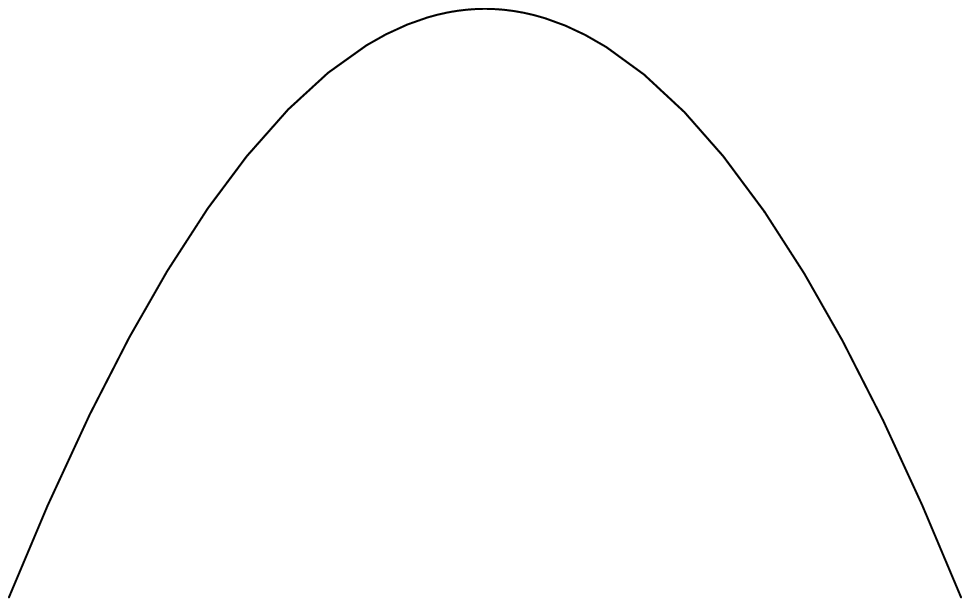}

\smallskip

i.e. it is very low as  $ I \approx 0 $ (e.g. no one would
consider the telephone "free-signal" as a musical
    masterpiece), it increases  monotonically as I increases
    reaching a maximum after which it becomes to decrease tending
    to zero for maximal information (no one would consider a
    random art of work an esthetical masterpiece).

    Since also the \emph{physical complexity} of an object,  quantified
    by Bennett's notion of \emph{logical depth} \cite{Bennett-88},
    \cite{Li-Vitanyi-97} have such a qualitative behavior as a
    function of the information, one could be tempted to
    conjecture to consider  \emph{physical complexity} as an
    approximate, empirical quantitative
    measure of \emph{esthetical value}, giving a firm foundation to the
    many speculations about the inter-relation among Esthetics and
    Complexity (as to Music consider, for example, the debate
    about the complexity of contemporary music and Xenakis' notion
    of \emph{mass} \cite{Battier-01} )
\end{enumerate}

Returning to the issue of estimating the relative pounds $ p_{Nature} \, , 
\, p_{Culture}  $ in eq.\ref{eq:Nature versus Culture} let us observe , 
furthermore, that the increasing rule of the cultural pound, anyway, cannot
delete the fact that the production, the propagation and the
reception of musical sounds are physical phenomena occurring in
Nature and described by the Laws of Acoustics.

The confusion among the concepts of \emph{esthetic consonance} and
of  \emph{physical consonance} has unfortunately led the
overwhelming majority of the existing literature to make confusion
between two completelly different issues:
\begin{enumerate}
    \item the issue whether the phenomenon of musical consonance
    has a physical ground
    \item the issue whether a positive answer to the previous
    issue would imply constraints on the Harmony's rules as to  
\emph{esthetic consonance}
\end{enumerate}
The confusion among these two, logically independent issues has:
\begin{itemize}
    \item led some scientists to claim that the fact that
    Atonal Music is based on a notion of \emph{esthetic
    consonance} differing from \emph{physical consonance} more
    than Tonal Music implies that Atonal Music is esthetically
    inferior than Tonal Music \cite{Frova-02}
    \item led some musicians to forget that the codification of
    new Harmonies with new \emph{esthetic consonances} cannot
    delete the fact that \emph{physical consonance} is a physical
    phenomenon occurring in Nature \footnote{a celebrated exposition of this 
argument was expressed,curiously, in "Human Sciences" by
    Claude Levi-Strauss, the father of all the discussions concerning the 
antinomy Nature versus Culture, in the "Overture" of
    \cite{Levi-Strauss-98} where he critizes the attitude by Pierre Boulez 
and many "serial thinkers" ,  of forgetting and  deleting the existence of 
what he calls \emph{"the first level of
    articulation"}; it should be observed, with this regard, that such a 
charge cannot be ascribed to Arnold Sch\"{o}nberg himself who clearly 
expressed, cfr. the cap.12 "Valutazione
    apollinea di un epoca dionisiaca" of \cite{Schoenberg-03} and the here 
cited analysis contained in the $7^{th}$ section "Sospensione ed 
eliminazione della tonalit\'{a}" and the $ 8^{th}$  section "La scala 
cromatica come fondamento della tonalit\'{a}" of the $19^{th}$ chapter 
"Alcune aggiunte e schemi che integrano
    il sistema" of \cite{Schoenberg-02}, that his theory of dissonance's 
emancipation considered dissonances as farer consonances in the sequences of 
harmonics}
\end{itemize}
It is interesting, at this point, to observe \cite{Fubini-01} that
many of the supporters of both the naturality and the esthetical
superiority of Tonal Orthodoxy  identified in Andreas
Werckmeister's \cite{Werckmeister-91} introduction \footnote{We
adhere here, for simplicity, to the common attribution; from an
historical perspective, anyway, we recall that the final
Werckmeister's codification was significantly preceded by a
plethora of other authors (cfr. e.g. the third section
"Necessit\'{a} del temperamento" of the third chapter "I rapporti
fra i suoni" of \cite{Tammaro-03})} of the chromatic
equally-tempered scale the first semen of the atonal heresy, as it
is condensed in the following sentence by Paul Hindemith about the
discovery of the equable-tempering:
\begin{center}
  \emph{"Anyone who has ever tasted the delights of pure intonation by the 
continual displacement of the comma in string-quartet
  playing \footnote{As remarked by Pierce \cite{Pierce-87} both the thesis 
here supported by Hindemith according to which the players of arc-strings
  instruments adopt natural intervals in their solos as well as the thesis 
often opposed to it, i.e.  the thesis claiming that they adopt well-tempered 
intervals seems to have been falsificated
  by concrete measurements: measurements show that the players of 
arc-strings
  instruments adopt nor natural neither well-tempered intervals} , must come 
to the conclusion that there can be no such
  thing as atonal music, in which the existence of tone
  relationship is denied. The decline in the value placed upon
  tonality is based on the system of equal temperament, a
  compromise which is presented to us by the keyboard as an aid in
  mastering the tonal world, and then pretends to be the world
  itself. One needs only to have seen how the most fanatical lower
  of the piano will close his ears in horror at the falseness of
  the tempered chords of his instrument, once he has compared them
  a few times with those produced by a harmonium in pure
  intonation, to realize that with the blessing of equal
  temperament there entered into the world of music - lest the bliss of 
musical mortals be complete - a curse as well: the curse of
  too easy achievement of tone-connections. The tremendous growth
  of piano music in the last century is attributable to it, and in
  the "atonal style" I see its final fullfillment - the uncritical
  idolatry of tempered tuning"; (cited from chapter 4 "Harmony", section 10 
"Atonality and Polytonality" of \cite{Hindemith-42})}
\end{center}

It is precisely this point that makes the whole subject
interesting.

To show the inconsistence  of  argumentations, such as Molly
Gustin's one \cite{Gustin-69} waved in \cite{Frova-02}, claiming
to furnish a mathematical proof of the esthetical superiority of
Tonal Music on Atonal Music would be a trivial and non interesting
matter if it didn't lead to the further conceptual step of
appreciating that:
\begin{enumerate}
    \item Classical Tonal Harmony is itself  invalidated by a mathematical
    inconsistence
    \item such a mathematical inconsistence affects the same
    formalization of the concept of Tonal Music,

    \item the issue is deeply linked with  the causes and the effects of  
Werckmeister's revolution
    \item it corroborates  Sch\"{o}nberg original
    viewpoint according to which Atonal Music was not a
    revolution but simply the last step in an historical process
    begun centuries before
\end{enumerate}
\newpage
\section{On the Musical Relativity Theory} \label{sec:On the Musical 
Relativity Theory}
The passage from Special Relativity to General Relativity \cite{Wald-84} may 
be
formalized as the passage from the:
\begin{axiom} \label{ax:principle of special relativity}
\end{axiom}
PRINCIPLE OF SPECIAL RELATIVITY
\begin{center}
\emph{All the laws of Physics are the same in all the inertial frames}
\end{center}
to the:
\begin{axiom} \label{ax:principle of general relativity}
\end{axiom}
PRINCIPLE OF GENERAL RELATIVITY
\begin{center}
\emph{All the laws of Physics are the same in all the frames}
\end{center}
Each of the two principles leads, in the respective theory, to the
corollaries:
\begin{corollary} \label{cor:corollary of special non-inferability}
\end{corollary}
COROLLARY OF SPECIAL NON-INFERABILITY
\begin{center}
\emph{No experiment allows an observer enclosed in a box to infer which is 
the particular inertial frame
of the box}
\end{center}
\begin{corollary} \label{cor:corollary of general non-inferability}
\end{corollary}
COROLLARY OF GENERAL NON-INFERABILITY
\begin{center}
\emph{No experiment allows an observer enclosed in a box to infer which is 
the particular frame
of the box}
\end{center}
Let us now suppose to pass from Physics to Music through the
following:

TABLE OF ANSATZS: $  PHYSICS \,  \mapsto  \, MUSIC $:

\bigskip

\begin{tabular}{|c|c|} \label{tab:table of ansatzs from Physics to Music}
  Physics & Harmony \\
  \hline
  frame & scale \\
  inertial frame & diatonic scale \\
  \hline
\end{tabular}

\bigskip

One obtains two musical theories consisting in a set of contraints
about the laws of harmony, the Special Theory of Musical
Relativity and the General Theory of Musical Relativity
\footnote{these theories shouldn't be confused with the Henry
Cowell's "Musical Relativity Theory" exposed in 1930 by the author
in \cite{Cowell-98} (curiously the stimulating locution "Musical
Relativity Theory" wasn't used by Cowell as a title for any of his
many writings \cite{Saylor-89}) that constitutes  a completely
different thing about which I will return later}, based,
respectively, on the following principles:
\begin{axiom} \label{ax:principle of special musical relativity}
\end{axiom}
PRINCIPLE OF SPECIAL MUSICAL RELATIVITY
\begin{center}
\emph{All the laws of Harmony are the same in all the diatonic scales}
\end{center}
and:
\begin{axiom} \label{ax:principle of general musical relativity}
\end{axiom}
PRINCIPLE OF GENERAL MUSICAL RELATIVITY
\begin{center}
\emph{All the laws of Harmony are the same in all the scales}
\end{center}

\bigskip

Driven by the analogy  one could be immediately  tempted to conjecture that 
axiom\ref{ax:principle of special musical
relativity} and axiom\ref{ax:principle of general musical
relativity} imply the following:
\begin{corollary} \label{cor:corollary of special musical non-inferability}
\end{corollary}
COROLLARY OF SPECIAL  MUSICAL NON-INFERABILITY
\begin{center}
\emph{A listener cannot infer the particular diatonic scale of a tonal 
musical piece}
\end{center}
\begin{corollary} \label{cor:corollary of general musical non-inferability}
\end{corollary}
COROLLARY OF GENERAL MUSICAL  NON-INFERABILITY
\begin{center}
  \emph{A listener cannot infer the scale of a musical piece}
\end{center}
A conceptual bug, anyway, inficiates the passage from axiom\ref{ax:principle 
of special musical
relativity} and axiom\ref{ax:principle of general musical
relativity} to, respectively, corollary\ref{cor:corollary of special musical
non-inferability} and corollary\ref{cor:corollary of general musical
non-inferability}: it implicitely assumes that the natural pound $ 
p_{Nature} $ of eq.\ref{eq:Nature versus Culture}, as to the net effect on 
the
listener, vanishes, a fact that we know to be false.

Let us leave aside, for a moment, such a conceptual bug precisely
in order of looking where it would take us.

The  adoption of the Special Relativity Theory and, hence, the imposition of 
corollary\ref{cor:corollary of general musical non-inferability}
would then  consist in the revolution realized by  Andreas
Werckmeister's equable-tempering who deleted the particular
"colours" of the different tonalities \footnote{E.g., as reported
by H. Keller in the section "Il carattere della tonalit\'{a} nel
Clavicembalo Ben Temperato" \cite{Keller-91}, D Major was thought
to be the tonality in which arcs are most shining, G minor was
associated in Baroque to the pompous music of official situations,
the preferred tonality for pastoral music was F major in the North
and G major in  the Sud etc. (for a more detailed historical
analysis cfr. the fourth section "Differenze espressive nel mondo
tonale" of the fifth chapter "Tra senso e metafisica" of
\cite{Tammaro-03})} and was, for this reason, strongly opposed; as
it is well known a determinant factor to the victory of
Werckmeister's ideas was the support they received \footnote{It
must be recalled, anyway, that from an historical point of view,
it is not clear whether Bach referred to the equable-tempering, to
the third, usually denoted as Werckmeister-III, of the six
temperings introduced by Werckmeister or to something else (cfr.
the section 7.3 "Scale temperate" of \cite{Frova-02} and  the
third section "Necessit\'{a} del temperamento" of the third
chapter "I rapporti fra i suoni" of \cite{Tammaro-03}) } from
Johann Sebastian Bach's \emph{Well-tempered Clavier} each of whose
two books (respectively BWV 846-869 AND BWV 870-893 \footnote{such
an indexing of the tonalities has led Piergiorgio Odifreddi, in
the  title of his divulgative article \cite{Odifreddi-99c}, to
talk of their "well-numbering"; such an indexing is, obviously,
recursive, its concrete implementation on computer being given by
the Miscellaneous - Music standard package of Mathematica
\cite{Wolfram-96} as well as by the Mathematica expressions
PlayScale and PlayChord described in \cite{Neuwirth-02}. One could
indeed think to go further comparing the notions of well-tempering
and well-ordering; since the tempering operation consists in
requiring the closure of the circle of fifths after 12 steps, i.e.
in altering the dynamics on the circle transforming it from a
quasi-periodic one to a periodic one by distributing, in a
suitable way, the pytagorical comma (the difference among 12
pytagorical fifths and 7 octaves, equal to the difference among
enharmonic notes in the pytagorical scale), such a procedure
results in the passage from a non-numerable infinity of different
tonalities to a finite number of different tonalities; while,
before tempering, a well-ordering of the different tonalities is
possible only assuming the Axiom of Choice, i.e. in an
intrinsically non-constructive way \cite{Schechter-98}, the
tempering's procedure transforms the well-ordering of different
tonalities into a constructive, and as we have seen even recursive
\cite{Bridges-98}, business} ) is made of 24 couple of preludes
and fugues in each of the 24 ionic and aeolian tonalities .

Exactly as the passage from Special Relativity to General
Relativity consists in  deleting the existence of a special,
privileged class of frames (the inertial ones), the passage from
Special Musical Relativity to General Musical Relativity would be
the  departure from Tonal Harmony  consisting in deleting the
existence of a special, privileged class of scales.

A first step in this direction follows Guerino Mazzola's
observation \cite{Mazzola-01} \cite{Mazzola-02a}
\cite{Mazzola-02b} that the mathematical structure underlying the
concept of tonality, whose intuitive content may be summarized as:
\begin{enumerate}
    \item the existence, at each instant of time, of a  tonal centre around 
which the
    musical discourse gravitates
    \item the changes of the tonal centres are performed by
    Modulation Theory codifying them as symmetry transformations
\end{enumerate}
is well larger than the classical notion of tonality
characterizing 24 particular ways of gravitating around the tonal
centre; a first enlargement of the notion of tonality consists in
passing from the notion of \emph{classical-tonality} to that of
\emph{gregorian-tonality} considering all the 84 modes of
gregorian music as distinct tonalities; the elimination of any
constraint in the characterization of the way of gravitating
around the tonal centre leads to a further radical enlargement of
the notion of tonality obtained introducing the notion of
\emph{Mazzola's tonality} and  taking into account all the 792
Mazzola's tonalities.

A partial passage from the Special Musical Relativity Theory to
the General Musical Relativity Theory would then consist in
replacing  the axiom\ref{ax:principle of special musical
relativity} with an analogous Principle of Special Musical
Relativity referred to such a generalized definition of the notion
of tonality.

It must be observed, with this regard, that nobody has composed
yet a sort of "Generalize well-tempered clavier", i.e. a
collection of two books of ordered musical composition, each book
consisting of an ordered collection of musical compositions in all
these possible "generalized tonalities".

As I will show later, anyway, the notion of Mazzola tonality may
itself be  generalized further.

\bigskip

Returning now to the mentioned conceptual bug, it allows
immediately to explain why  cor.\ref{cor:corollary of special
musical non-inferability} is trivially false:

our psycho-acoustic perception of sound is absolute and not
relative: the \emph{level of sonorous sensation (phon)}, the
\emph{level of subjective sonorous sensation (son)} , the
\emph{critical band}, the \emph{subjective height (mel)} and many
other phycho-acustic physical quantities depend on the frequency
(cfr. e.g. the cap.4 of \cite{Frova-02}, the cap.7 of
\cite{Pierce-87} and the cap. 3 of \cite{Everest-96}): so it is
not surprising that the phenomenon of \emph{perfect pitch}  (the
ability of some people to identify a  single note without a
reference to an external diapason) is an experimental evidence
that falsificates cor.\ref{cor:corollary of special musical
non-inferability}.

The contribution to $ p_{Nature} $ of all these psycho-acoustic
factors is, anyway, strongly lower than the contribution of an
external acoustic factor: the role of \emph{harmonics}, whose
discussion requires a brief report on the foundations of
Acoustics.

\newpage

\section{Elements of Acoustics}

Demanding to \cite{Everest-96},\cite{Rossing-Moore-Wheeler-02}, 
\cite{Pierce-87}, \cite{Frova-02} for any further information let us recall 
that from a physical point of view, a \emph{sonore signal} is codified by 
the associated-function $ s \; : \;  {\mathbb{R}}
\mapsto {\mathbb{R}} $ such that s(t) expresses the ratio among the 
perturbation
at time t produced by the \emph{sonore signal} in the pressure  of a fixed
point of the propagating medium and the  Pascal (the unity measure
of pression  in the International System I will assume from here
and beyond).

A \emph{sonore signal} s(t) may be seen, more precisely, as the
\emph{acoustical response}:
\begin{equation}
    s[t] \; = \; {\mathcal{R}}_{ac} [ i (t) ]
\end{equation}
(in pression) of the medium to an input $ i(t) $ produced by some
\emph{sonore source}.

The functional \emph{acoustical response} $ {\mathcal{R}}_{ac} $
is usually, nonlinear; beside in the phenomenon of
\emph{combinational notes}  I will discuss in section
\ref{sec:Horizontal and vertical rules in Tonal Harmony} the
deviation from linearity of the  \emph{acoustical response} is
neglibible, having no physical relevance \footnote{   The
introduced distinction among linear and nonlinear acoustic
response is, obviously, of a more general nature, as it  can be
appreciated introducing the mechanical and electric equivalents of
an acoustic system and looking at non-linear acoustic response
functions in terms of resistors with non-linear characteristic
(cfr. e.g. the section1.6  "Analogie elettriche, meccaniche ed
acustiche" of \cite{Everest-96} and the $ 10^{th} $ chapter
"Differential equations for electrical circuits" of
\cite{Hirsch-Smale-74})}.

It is then  natural, in these cases, to assume the validity of the
following:
\begin{axiom} \label{ax:principle of superposition}
\end{axiom}
PRINCIPLE OF SUPERPOSITION:
\begin{equation}
    {\mathcal{R}}_{ac} [ \lambda_{1} i_{1} (t) + \,  \cdots \,   + 
\lambda_{n} i_{n} (t)
    ] \; = \; \lambda_{1} {\mathcal{R}}_{ac} [ i_{1} (t) ] \,  +
    \cdots \, \lambda_{n} {\mathcal{R}}_{ac} [ i_{n} (t) ]
\end{equation}

\smallskip

A sonore signal will be said to be a \emph{sound} if its
\emph{associated function} is periodic.

the \emph{note} of a \emph{sound} of period T is then defined as
the \emph{frequency}  of its associated function  $ \nu :=
\frac{1}{T} $.

Given two notes $ \nu_{1} $ and $ \nu_{2} $ with  $ \nu_{1} \,
\leq \nu_{2} $
\begin{definition} \label{def:musical interval among two note}
\end{definition}
MUSICAL INTERVAL AMONG  $ \nu_{1} $ and $ \nu_{2} $:
\begin{equation}
    {\mathbf{mi}} ( \nu_{1} \, , \, \nu_{2} ) \; := \; 
\frac{\nu_{2}}{\nu_{1}}
\end{equation}

Fixed once and for all a reference note $ \nu_{ref} $ (e.g. by
the assumption  $ \nu_{ref} \; := \; C_{2} \; := \;  132 Hz $
adopted in the appendix 2.3 "Frequency and Glissando " of
\cite{Mazzola-02a}):
\begin{definition} \label{def:pitch of a note}
\end{definition}
PITCH OF THE NOTE $ \nu $ (IN CENTS):
\begin{equation}
    {\mathbf{p}}( \nu ) \; := \; \frac{1200}{ \log(2)} \log ( {\mathbf{mi}} 
( \nu \, , \, \nu_{ref} ) )
\end{equation}

A \emph{sound} of frequency  $ \nu $ is said to be \emph{pure} if
its associated function s(t) is a trigonometric function, i.e. it
is of the form $ s(t) = A \, \sin( 2 \pi \nu t ) $ or $ s(t) = A
\cos( 2 \pi \nu t) $, the positive factor $ A \in {\mathbb{R}}_{+}
$ being called the \emph{amplitude}.

By Fourier transform any \emph{sonore signal} may be seen as the
uncountable infinite superposition of \emph{pure sounds} of any
possible frequency:
\begin{theorem} \label{th:Fourier decomposition of a sonore signal}
\end{theorem}
FOURIER DECOMPOSITION OF A SONORE SIGNAL
\begin{equation}
    s(t) \; = \; \int_{0}^{\infty} d \nu ( a( \mu ) \cos( 2 \pi \mu t) \,
    + \, b( 2 \pi \mu t ) \sin ( 2 \pi \mu t ) ) \, d \mu
\end{equation}
where the two functions:
\begin{equation}
    a ( \mu ) \; = \;  \int_{- \infty} ^{+ \infty} d t
    s(t) \cos ( 2 \pi \mu t )
\end{equation}
and
\begin{equation}
    b ( \mu ) \; = \; \int_{- \infty} ^{+ \infty}
    d t
    s(t) \sin ( 2 \pi \mu t )
\end{equation}
are called the \emph{Fourier components} of the \emph{sonore
signal}

If the \emph{sonore signal} is itself a \emph{sound} of
\emph{note} $ \nu $ , furthermore, theorem\ref{th:Fourier
decomposition of a sonore signal} immediately implies that:
\begin{enumerate}
    \item only a countable infinity of \emph{pure sounds} contribute to its 
\emph{Fourier components}
    \item the notes of this infinity of Fourier-contributing  
\emph{pure-sounds} (called the \emph{harmonics} of the sound) , are the 
multiple
    integers of the note of the considered sound
\end{enumerate}
as is stated by the following:
\begin{corollary} \label{cor:harmonic decomposition of a sound}
\end{corollary}
HARMONIC DECOMPOSITION OF A SOUND:
\begin{equation}
  a ( \mu ) \; = \; \sum_{n=0}^{\infty} a_{n} \delta ( \nu - n \mu
  )
\end{equation}
\begin{equation}
  b ( \mu ) \; = \; \sum_{n=0}^{\infty} b_{n} \delta ( \nu - n \mu)
\end{equation}

It may be useful, a this point, to introduce the following:
\begin{definition}
\end{definition}
HARMONIC SEQUENCE OF THE NOTE $ \nu $:

the sequence of notes $ \{ harmonic( \nu , n ) \}_ { \{ \: n \in
{\mathbb{N}}_{+} \} } $:
\begin{equation}
    harmonic( \nu , n ) \; :=  \; (n+1) \, \nu
\end{equation}

I will adopt, from here and beyond, the usual complex
representation of harmonic motion, representing a sound $ s(t) :=
\frac{a_{0}}{2} + \sum_{n=1}^{\infty} ( a_{n} \cos ( n \omega t )
+ b(n) \sin ( n \omega t ) ) $ through the complex function:
\begin{equation}
    \tilde{s} (t) \; := \; \sum_{n = - \infty}^{+ \infty} c_{n} e^{i n 
\omega t}
\end{equation}
where:
\begin{equation}
    c_{n} \; := \; \left\{%
\begin{array}{ll}
    \frac{a_{0}}{2}, & \hbox{if $ n = 0 $;} \\
    \frac{a_{n} - i b_{n}}{2}, & \hbox{if $ n > 0 $;} \\
    \frac{a_{n} + i b_{n}}{2}, & \hbox{if $ n < 0 $} \\
\end{array}%
\right.
\end{equation}
dropping the tilde from here and beyond.

As to arbitrary sonore signals, furthermore, I will adopt, from
here and bejond, the analogous complex representation of Fourier
integrals:
\begin{equation}
    s(t) \; = \; \int_{- \infty}^{+ \infty} d \nu  S( \nu) e ^{ 2 \pi i \nu 
t}
\end{equation}
\begin{equation}
   S( \nu ) \; = \; \int_{- \infty}^{+ \infty} d t s(t) e ^{-  2 \pi i \nu 
t}
\end{equation}

The concept that a \emph{sound} is the composition of other sounds
is the source of much confusion, that may be completelly avoided
taking strongly into account the following remarks:
\begin{enumerate}
    \item the \emph{amplitude} of a \emph{pure sound}, ruling the
   \emph{intensity} of the involved pressure-wave, may be seen as
   a parameter measuring its \emph{volume}
    \item the only harmonic  contributing to the harmonic decomposition of  
a \emph{pure
    sound} with note $ \nu $ is $ \nu $ itself \footnote{Among all the 
physical nonsense of Sch\"{o}nberg's analysis
    about harmonics in the first section of the fourth chapter of 
\cite{Schoenberg-02} this is a point in which he makes
    trivially false statements ending up with the "harmonics of harmonics" 
and the  involved regressum ad
    infinitum (SIC!)}
    \item if the  amplitudes of two \emph{pure sounds}  are, roughly 
speaking, of the same
     order of size, they are perceived by our ears as
    two distinct \emph{sonore signals}.
    \item if the  amplitudes of two \emph{pure sounds}  are, roughly 
speaking, of
    different order of size, our ears perceive them as a unique \emph{sonore
     signal} whose associate function is the sum of the associated
     functions of the two \emph{pure sounds}
     \item the reason why the \emph{pure-sounds} contributing to
     the harmonic decomposition of a  \emph{non-pure sound} are
     perceived as a unique  \emph{sonore
     signal} is that the amplitude $ a_{n} $ and  $ b_{n} $
     decrease enough quickly as n grows, as is implied by the
     Bessel's inequality:
\begin{equation}\label{eq:Bessel's inequality}
    \frac{a_{0}^{2}}{2} \, + \, \sum_{k=1}^{n} ( a_{k}^{2} +
    b_{k}^{2} ) \; \leq \; \frac{2}{T} \int_{-
    \frac{T}{2}}^{\frac{T}{2}} \, s^{2} (t) dt \; \; \forall n \in
    {\mathbb{N}}
\end{equation}
     so that, for different values of n, they are always
     of different order of size
\end{enumerate}

\bigskip

The real situation is, anyway, more complex since:
\begin{enumerate}
    \item any concretely  occurring \emph{sonore signal} has
    finite support,\emph{sounds} are just a mathematical
    idealization not occurring in reality
    \item one usually adopts  the previously introduced
    metaphorical locution according to which our ears hear the
    harmonical components of a sonore signal; the human auditory
    anatomic behavior is not that of a computer
    getting from the ear device the  input consisting in the whole 
specification of the function
    s(t) and after that computing its whole Fourier transform
\end{enumerate}
These observations have led Guerino Mazzola (cfr. the $ 2^{th} $
chapter "Topography" and the $ 15^{th} $ part "Sound" of
\cite{Mazzola-02a}), referring  to Jean Molino's scheme concerning
the three-parts a musical communication's stream is made of:
\begin{description}
    \item[Poiesis]  the production of music by a creator
    \item[Neutral Niveau] the transmitted musical message
    \item[Esthesis] the reception of the musical message by a
    listener
\end{description}
to contest the traditional viewpoint of Acoustics giving for
granted that the notions of note and pitch belong to the neutral
level.

Mazzola doesn't anyway seems to appreciate the role of the Theorem
of Tonal Indetermination we are going to introduce.

Given a sonore signal s(t) let us introduce the following
probability distributions over:  $ ( {\mathbb{R}} \, , \, {
\mathcal{B}} ({\mathbb{R}}) )$ \footnote{$ { \mathcal{B}}
({\mathbb{R}}) $ denoting the Borel-$\sigma $-algebra of $
{\mathbb{R}} $}:
\begin{equation} \label{eq:probability distribution of a signal sonore}
    P_{s} (t) \; := \;  \frac{ | s(t) |^{2} } { E  }
\end{equation}
\begin{equation} \label{eq:probability distribution of the Fourier transform 
of a signal sonore}
    P_{S} ( \nu ) \; := \;  \frac{ | S( \nu ) |^{2} } { E }
\end{equation}
where, owing to Parseval's equation, one has that:
\begin{equation}
    E \; := \; \int_{ - \infty }^{  +  \infty  } \, \, dt \, | s(t) |^{2} \; 
=
    \; \int_{ - \infty }^{  +  \infty  } \, d\nu | S( \nu) |^{2}
\end{equation}
It may be proved that (cfr. the section 5.10.6 "Sul concetto di
trasformazione" of the $ 5^{th} $ chapter "La trasmissione
dell'informazione" of \cite{Angeleri-99}, the appendix N-3
"Analytic Signals and the Uncertainty Relation" of \cite{Reza-94},
the section 8.10 "The Uncertainty Principle" of \cite{Hamming-89},
the section 19.4 "Wavelet transformation" of
\cite{Harris-Stocker-98} and the part D "Wavelet Analysis"  of
\cite{Bremaud-02}) that:
\begin{theorem} \label{th:theorem of note's indetermination}
\end{theorem}
THEOREM OF NOTE'S INDETERMINATION:

\begin{hypothesis}
\end{hypothesis}
\begin{equation}
s( t ) \; = \; o_{t \rightarrow \pm \infty}(t)
\end{equation}
\begin{thesis}
\end{thesis}
\begin{equation}
    \Delta t \, \Delta \nu \; \geq \; \frac{1}{4 \pi}
\end{equation}
where $ \Delta t $ and $ \Delta \nu $  are the standard-deviations
of, respectively, the  probability distributions
\ref{eq:probability distribution of a signal sonore} and
\ref{eq:probability distribution of the Fourier transform of a
signal sonore}.

The conceptual deepness of theorem\ref{th:theorem of note's
indetermination} was fully analyzed by Dennis Gabor who:
\begin{enumerate}
    \item conjectured that is may be seen as a particular case $ out(t) \, = 
\, in (t) $  of an
    information-theoretic indetermination's theorem generalizing
    to arbitrary channels with band-width $ \Delta B $ such
    that:
\begin{equation}
    supp ( T ( \omega ) ) \; = \; \Delta \nu
\end{equation}
Nyquist's criterium:
\begin{equation}  \label{eq:Nyquist's criterium}
    \Delta B \, \Delta t \; = \; constant
\end{equation}

on the minimum temporal distance $ \Delta
    t $ among two Dirac-delta's impulses such that they can be distinguished 
at the output of a Nyquist channel, i.e. a channel with transfer function:
\begin{equation}
    T ( \omega ) \; := \; \frac{Out(\omega )}{In ( \omega) } \; =
    \; | A ( \omega ) | \, e^{ - \, i \, \beta ( \omega ) }
\end{equation}
with:
\begin{eqnarray}
   A ( \omega ) & \; = \; \left\{%
\begin{array}{ll}
    A_{0}, & \hbox{if $ | \omega | \, \leq \, \omega_{0} $;} \\
    0, & \hbox{otherwise.} \\
\end{array}%
\right.     \\\sqrt{a }
  \beta ( \omega ) \; = \;  \omega \, t_{0}
\end{eqnarray}
   \item understood its role in determining the invariance of the acquired 
information in any signal's analysis performed by a central computational 
unit (e.g. the human brain) with an inner time-clock $ \Delta t $ getting 
the sound s(t)  through the sampling performed by a device (e.g. the human 
ear) with  linear transfer function $ {\mathcal{D}}_{L}$ \footnote{Exactly 
as the acoustical response of the medium, the transfer function  of  the 
human ear may
be decomposed in its linear and non linear parts:
\begin{equation}
    {\mathcal{D}} \;  = \;   {\mathcal{D}}_{L} \, + \, {\mathcal{D}}_{NL}
\end{equation}
$  {\mathcal{D}}_{NL}  $  conspires together with
the nonlinear part of the medium acoustic response to generate the
non-linear effects such as \emph{combinational tones} discussed in
the section \ref{sec:Horizontal and vertical rules in Tonal
Harmony} and many other effects about which we demand to the $
3^{th} $ chapter "L'orecchio e la percezione del suono" of
\cite{Everest-96}, to the appendix B "Auditory Physiolosy and
Psychology" of \cite{Mazzola-02a} or the text
\cite{Zwicker-Fastl-99} enterely ddvoted to such issue.Exactly as
to the acoustcial response, we will ignore from here and beyond
the physical effect involved with  $  {\mathcal{D}}_{NL}  $ not
being relevant for the issue discussed in this paper; we will
furthermore refer to this approximation as the assumption of the
following generalization of axiom\ref{ax:principle of
superposition}, i.e.:
\begin{axiom} \label{ax: extented principle of superposition}
\end{axiom}
EXTENDED PRINCIPLE OF SUPERPOSITION:
\begin{equation}
    {\mathcal{D}} \, \circ \,  {\mathcal{R}}_{ac} [ \lambda_{1} i_{1} (t) + 
\,  \cdots \,   + \lambda_{n} i_{n} (t)
    ] \; = \; \lambda_{1}  {\mathcal{D}} \, \circ \,  {\mathcal{R}}_{ac} [ 
i_{1} (t) ] \,  +
    \cdots \, \lambda_{n} {\mathcal{D}} \, \circ \,  {\mathcal{R}}_{ac} [ 
i_{n} (t) ]
\end{equation}} whose band-width is less or equal to the one admitted by
Nyquist's
criterium\ref{eq:Nyquist's criterium} for sampling;

such an invariance follows, precisely, from the fact that the Gabor's 
transformation  of the signal
(that I am going to introduce) saturates eq.\ref{th:theorem of note's 
indetermination}, i.e. it respects
it with the equal sign, so that its  time-frequency representation is 
constitued of  rectangular  windows   of linear
dimensions $ ( 2 \, \sqrt{a } \, , \, \frac{1}{ \sqrt{a} }    )  $
and hence of constant area $ A \, = \, 2 $. The concept's of Gabor
transformation is based on the idea that, if one is interested in
getting the local-information about a signal localized around the
time $ t \, = \, b $ it is reasonable to generalize the Fourier's
tranform introducing into it a localized  gaussian weight-function
centered around $ t \, = \, b $. The generalization of this
concept to a more general class of possible small localized
weight-functions led Jean Morlet to introduce the following
notions:

\begin{definition}
\end{definition}
MOTHER WAVELET:

a function $ \psi \, \in \, L^{2} (  {\mathbb{R}} ) $ such that:
\begin{equation}
    0 \; < \; \int_{- \infty}^{+ \infty} \frac{ | \Psi ( \nu )   |^{2} }{ | 
\psi |
} \; < \; + \infty
\end{equation}
where, according to our general notation, $ \Psi ( \nu ) $ denotes
the Fourier transform of $   \psi (t) $.

Given a mother wavelet $  \psi  $:
  \begin{definition}
  \end{definition}
WAVELET GENERATED BY $  \psi  $:

the two parameter family of functions $ \{ \psi_{a,b} (t) \}_{b \in 
{\mathbb{R}} , a \in {\mathbb{R}} \, - \, \{ 0 \} } $ given by:
\begin{equation}
    \psi_{a,b} (t) \; := \; \frac{1}{ | a |^{ \frac{1}{2}}}  \,
\psi ( \frac{ t-b }{a} )
\end{equation}
\begin{definition}
\end{definition}
WAVELET TRANSFORM W.R.T.  $  \psi  $ OF $ f \, \in \, L^{2}  ( {\mathbb{R}} 
) $:
\begin{equation}
    T_{\psi} ( a \, , \, b )  \; := \; \int_{- \infty} ^{+ \infty}
dt \, f(t) \psi^{\star} (t)
\end{equation}

We can at last introduce the following:
\begin{definition}
\end{definition}
GABOR TRANSFORMATION:

the \emph{wavelet tranform} generated by a gaussian
\emph{\emph{mother wavelet}}

\end{enumerate}

It is important to stress that, though the Mathematics
underlying theorem\ref{th:theorem of note's indetermination} is
exactly that of the particular application of the Heisenberg's
Indetermination Theorem in Quantum Mechanics:
\begin{equation}
    < \alpha | ( \Delta \hat{A} )^{2} ( \Delta \hat{B} )^{2} |
    \alpha > \; \geq \; \frac{1}{4} | < \alpha | [ \hat{A} ,
    \hat{B} ]  | \alpha > |^{2}
\end{equation}
to the couple of observables \emph{position} and \emph{momentum} $
\hat{x} \, , \, \hat{p} $ \footnote{I would like to advise the
reader that in $ 5^{th} $ chapter "La trasmissione
dell'informazione" of \cite{Angeleri-99} it is erroneously stated
that theorem\ref{th:theorem of note's indetermination} implies the
indetermination's relation \emph{time} and \emph{energy} therein
referred as Heisenberg's Principle. As  stressed by any elementary
Quantum Mechanics' manual (cfr. e.g.  the section 5.6  "Teoria
perturbativa dipendente dal tempo" of \cite{Sakurai-90}), anyway,
the indetermination's relation \emph{time} and \emph{energy},
inferred in the framework of first-order perturbative theory, is
conceptually radically different from Heisenberg's Principle; in
particular, time is not an observable in Quantum Mechanics as it
is significantly testified by the conceptual bugs one has to face
on trying to define a \emph{quantum clock} about which I demand to
section 12.7 "The measurement of time" of \cite{Peres-95}, as well
as to \cite{Busch-02} and \cite{Sala-Mayato-Alonso-Egusquiza-02}},
the issue is quite different from a  physical point of view; it
should be mentioned, with this regard, that the adoption of the
quantities $ \Delta t $ and $ \Delta \nu $ as measures of the
uncertainty in, respectively, time and frequency of a sonore
signal s(t) is someway arbitary, since no physical ground exists
in the assumption of the involved probability, respectively,
distributions $ P_{s} (t) $ and  $ P_{S} (\nu) $ mathematically
assumed by mimicking the probability distribution of
quantum-mechanical wave-functions.

The standard deviation of $ P_{S} (\nu) $ is, anyway, certainly a
better measure of the uncertainty in frequency of s(t) than the
semi-distance among its first two zeroes of a function only
approximating $ S(\nu) $  that has been adopted in the section 3.3
"Pacchetti d'onda. Indeterminazione tonale" of \cite{Frova-02} on
the basis of the section 12.7 "Pacchetti d'onda. Velocit\'{a} di
gruppo " of \cite{Amaldi-Bizzarri-Pizzella-86}.

Let us now return to the information-theoretical point of
view to Music briefly introduced in section\ref{sec:Introduction}, taking 
into account, for simplicity,
the production of a monofonic musical piece made of notes of the
same duration:

such a musical piece may then be seen as a  \emph{string} or  as a
\emph{sequence} on the alphabet $ {\mathbb{R}}_{+} $ of all the
possible notes.

Since:
\begin{equation}
    card ( ( {\mathbb{R}}_{+} ) ^{\star} ) \; = \;  2^{ card 
({\mathbb{R}}_{+} )
    } \; = \; 2^{\aleph_{1}} \; = \; \aleph_{2}
\end{equation}
\begin{equation}
    card ( ( {\mathbb{R}}_{+} ) ^{\infty} ) \; = \;  2^{ card(  
({\mathbb{R}}_{+} )^{\star} )
    } \; = \; \aleph_{3}
\end{equation}
both considering terminating and nonterminating musical pieces one
has that the set of all the possible musical messages has a
cardinality higher than the continuum one, i.e. is made of "a more
than continuously infinity" of possible musical executions with
the consequential troubles as far as codification of musical
messages is concerned.

This simple consideration allows to understand why, from the
beginning of its structurally organized formalization of Music,
humanity tried to limit the space of the possible musical messages
introducing strongly smaller alphabets.

Different cultures pursued this goal in different ways: the
structural similarity of the obtained results, anyway, is a very
strong evidence of the role that the common underlying   Laws of
Nature had in such a process.

We have touched, at this point, the key argument of this paper
enclosed in the following:

\begin{question} \label{qu:question on the more natural musical messages}
\end{question}
QUESTION ON THE "MORE NATURAL MUSICAL MESSAGES" DETERMINED BY THE
HARMONIC SEQUENCES

\textbf{At which extent does the "naturality"  of the harmonic
sequences determines the "naturality" of some particular musical
messages ?}

\bigskip

The naive answer to question\ref{qu:question on the more natural
musical messages} one usually reads many times is the following:

\begin{answer} \label{an:the naife answer}
\end{answer}
THE NAIFE ANSWER:

\textbf{The natural role played, in the Laws of Nature, by the
first five harmonics determines the Nature-induced maximally
consonant triad: the major triad}

\bigskip

Since, as I will extensively discuss, the naife answer\ref{an:the
naife answer} is false,  let us proceed more carefully starting to
discuss, first of all, the "natural" role of the first harmonic.

A strong  empirical evidence of its pervasivity in different
cultures strongly leads to assume the following:

\begin{axiom} \label{ax:axiom of perception of repetition for the first 
harmonic}
\end{axiom}
AXIOM OF PERCEPTION OF REPETITION FOR THE FIRST HARMONIC

\textbf{the natural role played in the Laws of Nature by the first
harmonic is such that, given a sequences of notes $ \{ \nu_{n} \}
_{n \in {\mathbb{N}}}  $, the sequence of notes $ \{ harmonic(
\nu_{n} , 1 )  \} _{n \in {\mathbb{N}}}  $ is perceived as a
repetition of the sequence $ \{ \nu_{n} \} _{n \in {\mathbb{N}}} $
at an higher range}

\bigskip

Axiom\ref{ax:axiom of perception of repetition for the first
harmonic}  played a basic role in all the ways in which different
cultures constructed \emph{scale of notes} they adopted as musical
alphabet; conseguentially, in all cultures, a basic role was
played, in such a construction, by the following notion:
\begin{definition}
\end{definition}
SCALE RANGE OF A NOTE $ \nu $:
\begin{equation}
    scale-range(\nu) \; := \; [ \nu \, , \, harmonic(\nu , 1) ]
\end{equation}
Let us now  observe, that assuming the axiom\ref{ax:axiom of
perception of repetition for the first harmonic}, the construction
of any kind of \emph{scale of notes} starting from a note $ \nu $
may be limited to its \emph{scale range}: outside the \emph{scale
range} everything is simply repeated by multiplying each note of
the considered scale for a suitable power of two.

Contrary, any note lying  outside the scale range of $ \nu $ may
be rescaled to it by dividing it for a suitable power of two.

These operations may be easily performed through the following:
\begin{definition} \label{def:rescaling function}
\end{definition}
RESCALING FUNCTION TO THE SCALE-RANGE OF  $ \nu $:

the function $  R_{\nu} \; : \; {\mathbb{R}}_{+} - scale-range
(\nu ) \, \mapsto \; {\mathbb{R}}_{+} $ identified by the
following constraints:
\begin{itemize}
  \item for $ \mu \; < \; \nu $ :
\begin{equation}
  R_{\nu} ( \mu ) \; = \; 2^{ f_{\nu} ( \mu ) } \, \cdot \, \mu
\end{equation}
with:
\begin{enumerate}
  \item
\begin{equation}
  f_{\nu} ( \mu ) \; \in \; {\mathbb{N}}
\end{equation}
  \item
\begin{equation}
  2^{ f_{\nu} ( \mu ) } \, \cdot \, \mu \; \geq \; \nu
\end{equation}
  \item
  \begin{equation}
  2^{ f_{\nu} ( \mu ) } \, \cdot \, \mu \; <  \;  2 \nu
\end{equation}
\end{enumerate}
  \item  for $ \mu \; > \;  2 \nu $ :
\begin{equation}
  R_{\nu} ( \mu ) \; = \; 2^{- \, f_{\nu} ( \mu ) } \, \cdot \, \mu
\end{equation}
\begin{enumerate}
  \item
\begin{equation}
  f_{\nu} ( \mu ) \; \in \; {\mathbb{N}}
\end{equation}
  \item
\begin{equation}
  2^{ - \, f_{\nu} ( \mu ) } \, \cdot \, \mu \; \geq \; \nu
\end{equation}
  \item
  \begin{equation}
  2^{ - \, f_{\nu} ( \mu ) } \, \cdot \, \mu \; <  \;  2 \nu
\end{equation}
\end{enumerate}
\end{itemize}

\bigskip

Some simple consideration leads to the following:
\begin{lemma} \label{le:explicit formula for the rescaling function}
\end{lemma}
EXPLICIT FORMULA FOR THE RESCALING FUNCTION:
\begin{equation}
  R_{\nu} ( \mu ) \; = \;
  \begin{cases}
    2^{ Int ( \frac{ \log \nu \, - \, \log \mu }{ \log 2} \, + \, 1 )} \: 
\cdot \: \mu & \text{if $ \mu \, < \, \nu $  }, \\
    2^{ - \, Int ( \frac{ \log \mu \, - \, \log \nu    }{ \log 2}     )}  
\cdot \: \mu  & \text{if $ \mu \, > \,  2 \, \nu $}.
  \end{cases}
\end{equation}
\begin{proof}
\begin{itemize}
  \item for $ \mu \; < \; \nu $ :

  passing to the logarithms in the disequalities one obtains that:
\begin{equation}
  f_{\nu} ( \mu ) \; \in \; [ A_{\nu} ( \mu ) \, , \, A_{\nu} ( \mu
  ) + 1 )
\end{equation}
where:
\begin{equation}
  A_{\nu} ( \mu ) \; = \; \frac{ \log \mu \, - \, \log \nu    }{ \log 2}
\end{equation}
The thesis immediately follows imposing the constraint that $
f_{\nu} ( \mu )  $ has to be integer

  \item for $ \mu \; > \;  2 \nu $ :

  passing to the logarithms in the disequalities one obtains that:
\begin{equation}
  f_{\nu} ( \mu ) \; \in \; [ A_{\nu} ( \mu ) \, , \, A_{\nu} ( \mu
  ) + 1 )
\end{equation}
where:
\begin{equation}
  A_{\nu} ( \mu ) \; = \; \frac{ \log \mu \, - \, \log \nu }{ \log
  2} \; - \; 1
\end{equation}
The thesis immediately follows imposing the constraint that $
f_{\nu} ( \mu )  $ has to be integer
\end{itemize}
\end{proof}

Both the \emph{scale-range} of a note $ \nu $ and the rescaling
function to its scale-range may be computed through the following
Mathematica's expressions:
\begin{verbatim}
harmonic[nu_,n_]:=(n+1)*nu

scalerange[nu_]:=Interval[nu,harmonic[nu,1]]

rescalingtorange[nu_,mu_]:= If[mu < nu,
    mu*Power[2, IntegerPart[Times[Log[mu] - Log[nu], Power[2, -1]]] +
1], If[mu>2nu,
mu*Power[2,-IntegerPart[Times[Log[mu]-Log[nu],Power[2,-1]]]],
"undefined"]]
\end{verbatim}

taken from my Mathematica \cite{Wolfram-96} noteboook
\emph{Mathematical Music Theory} reported in the
section\ref{sec:Mathematica's notebook Mathematical Music Toolkit}
\bigskip

The rescaling function allows to introduce the following useful
generalization of Congruence Theory \cite{Manin-Panchishkin-95}:

given three real numbers $ a   ,  b , c  \in {\mathbb{R}} $:
\begin{definition} \label{def:congruence modulo powers}
\end{definition}
a IS CONGRUENT TO b MODULO POWERS OF c ( $ a \,  =  \, b \; mod \,
c^{{\mathbb{Z}}} $):
\begin{equation}
    \exists n \in {\mathbb{Z}} \; : \; a \, = \, b \, mod \, c^{n}
\end{equation}
\begin{definition} \label{def:residue class modulo powers}
\end{definition}
RESIDUE CLASS OF a MODULO POWERS OF c:
\begin{equation}
    [ a ]_{c^{{\mathbb{Z}}}} \; := \; \{ \, x \in {\mathbb{R}} \: :
    \: x \, = \,  a \,   mod \, c^{{\mathbb{Z}}} \, \}
\end{equation}
For $ c \in {\mathbb{N}} $ we can introduce, furthermore, the
following:
\begin{definition} \label{def:set of the residue classes modulo powers}
\end{definition}
SET OF THE RESIDUE CLASSES MODULO POWERS OF c:
\begin{equation}
  {\mathbb{Z}}_{c^{{\mathbb{Z}}}} \; := \;  { [a]_{c^{{\mathbb{Z}}}} \, : \, 
a \in [ 0 , c-1 ] }
\end{equation}
\newpage
\section{Horizontal and vertical rules in Tonal Harmony} 
\label{sec:Horizontal and vertical rules in Tonal Harmony}
Tonal Harmony is made of two ingredients:
\begin{description}

    \item[VERTICAL RULES IN SCORES]   The theory of physical consonance
    among sounds inside a fixed tonality, ruled according to the
    following:

\begin{axiom} \label{ax:axiom of the naturality of esthetics}
\end{axiom}
AXIOM OF THE NATURALITY OF ESTHETICS:
\begin{equation}
    \text{esthetic consonance} \; = \; \text{physical consonance}
\end{equation}
The existence of a physical acoustical reason underlying the
phenomonological evidence that certain couples of notes are
perceived by our ears as consonant, i.e. the net acoustical input
sounds good, while other couples of notes are perceived by our
ears as dissonant,i.e. the net acoustical input sounds bad, was
first suggested in 1638 by Galileo Galilei \cite{Galilei-59} in
the following famous passage:
\begin{center}
  \emph{"SALVIATI: $ \cdots $ Returning now to the original subject of 
discussion, I assert that the ratio of a musical interval
  is not immediately determined by the length, size, or tension of the 
strings, but rather by the ratio of their frequencies, that
  is, by the number of pulses of air waves which strike the tympanum of the 
ear, causing it also to vibrate vith the
  same frequency. This fact established, we may possibly explain why certain 
pairs of notes differing in pitch produce a pleasing
  sensation, other a less pleasant effect, and still others a
  disagreeable sensation. Such an explanation would be tantamount
  to an explanation of the more or less perfect consonances and of
  dissonances. The unpleasant sensation produced by the latter
  arises, I think, from the discordant vibrations of two different
  tones which strike the ear out of time. Especially harsch is the
  dissonance between notes whose frequencies are incommensurable;
  such a case occurs when one has two strings in unison and sounds
  one of them open, together with a part of the other which bears
  the same ratio to its whole length as a side of a square bears
  to the diagonal; this yelds a dissonance similar to the
  augmented fourth or diminished fifth.}
\end{center}

\begin{center}
  \emph{Agreeable consonance are pairs of tones which strike the ear with a 
certain regularity;
  this regularity consists in the fact that the pulses delivered by the two 
notes, in the same interval of time, shall be commensurable in
  number, so as not keep the ear drum in perpetual torment,
  bending in two different directions in in order to yeld to the
  ever-discordant impulses ". from part 146-147, pagg.103-104 of 
\cite{Galilei-59}}
\end{center}

We see that:
\begin{itemize}
    \item Galileo didn't know the phenomenon of harmonics and made the error 
of
    thinking that strings produce pure sounds
    \item he identified the source of the physical consonance between
these supposed pure sounds in the  condition of  commensurability of their 
frequencies (and
the consequent regularity  of the tympanum's sollecitation).
\end{itemize}

\bigskip
The first attempt of furnishing a  quantitative measure of
physical consonance was performed by L. Euler \cite{Euler-26a}
through the introduction of a function \emph{gradus suavitatis}
expressing the degree of consonance among two notes in terms of
the "simplicity" of their ratio. Even more remarkable, was
,anyway, the introduction  of the space on which the \emph{gradus
suavitatis} is  defined, that I will call, following Mazzola (cfr.
the $ 15^{th} $ part "Sound" of \cite{Mazzola-02a}) the
\emph{Euler's space}.

Gixen an integer positive number $ n \in {\mathbb{N}} $ and
denoted by $ \pi( i) $ the $ i^{th} $ prime number:
\begin{definition}
\end{definition}
n-LEVEL EULER'S COORDINATION

the  map $ c_{n} \, : \, {\mathbb{R}}^{n} \, \mapsto \,
{\mathbb{R}} $:
\begin{equation}
    c_{n}( x_{1} \,  , \,  \cdots  \, , \, x_{n} ) \; := \;
    \prod_{i=1}^{n} \pi( i) ^{ x_{i} }
\end{equation}
One has that:
\begin{lemma} \label{lem:on the injectivity of Euler's coordinations}
\end{lemma}
ON THE INJECTIVITY OF EULER'S COORDINATIONS:

\begin{enumerate}
    \item $ c_{n} $ is not injective
    \item $ c_{n} |_{ {\mathbb{Q}}^{n} } $ is injective
\end{enumerate}
\begin{proof}
\begin{enumerate}
    \item  the set $ c_{n}^{- 1} (y) $ of the counterimages of a
    positive real $ y \, \in \, [ 0 \, , \, + \infty ) $ is
    given by the solutions of the following equation:
\begin{equation}
    \prod_{i=1}^{n} \pi( i) ^{ x_{i} } \; = \; y
\end{equation}
i.e.:
\begin{equation}
   \sum_{i=1}^{n} \, x_{i} \, \log \pi( i) \; = \; \log y
\end{equation}
or:
\begin{equation}
  x_{2} \; = \; \frac{ \log y \, - \, \sum_{i=1}^{n} \, x_{i} \, \log \pi( 
i) }{ \log  \pi(1)  }
\end{equation}
so that:
\begin{equation}
    card ( c_{n}^{- 1} (y) ) \; = \; \aleph_{1} \; > \; 1
\end{equation}
    \item let us proceed by induction:
\begin{itemize}
    \item proving the injectivity of $ C_{2} $:

    since the integer solutions of the equations:
\begin{equation}
     2^{x_{1}} 3^{ x_{2} } \; = \; 2^{y_{1}} 3^{ y_{2} }
\end{equation}
are the same of the equation:
\begin{equation}
     \frac{ x_{1} \, -  \, y_{1}   }{ x_{2} \,  -  \, y_{2}   } \; = \; - \, 
\log_{2} 3
\end{equation}
the thesis immediately follows by the fact that $ \frac{ x_{1} \,
- \, y_{1}   }{ x_{2} \,  -  \, y_{2}   } $ is rational or
undetermined while $ - \, \log_{2} 3 $ is irrational
    \item proving that the injectivity of $ c_{n} $ is implied by
    the injectivity of $ c_{n-1} $:

    since by the inductive hypothesis:
\begin{equation}
    card ( c_{n-1}^{- 1} (y) ) \; \leq \; 1 \; \; \forall y \in [ 0
    , \infty )
\end{equation}
to prove that:
\begin{equation}
    card ( c_{n}^{- 1} (y) ) \; \leq \; 1 \; \; \forall y \in [ 0
    , \infty )
\end{equation}
it is sufficient to show that:
\begin{equation}
    ( \, \pi (n) ^{x_{n}} \, = \, \pi (n) ^{y_{n}} \, ) \;
    \Leftrightarrow \;  ( \, x_{n} \,  = \,  y_{n} )
\end{equation}
But this is a trivial consequence of the bijectivity of the
exponential function
\end{itemize}
\end{enumerate}
\end{proof}

\smallskip

\begin{definition} \label{def:Euler's notes}
\end{definition}
EULER'S NOTES:
\begin{equation}
    {\mathcal{N}}_{Euler} \; := \; \{ \nu_{ref} \,  c_{3} ( \vec{x} )  \, , 
\, \vec{x} \in {\mathbb{Q}}^{3} \}
\end{equation}
Given an Euler's note $ \nu $ the lemma\ref{lem:on the injectivity
of Euler's coordinations} allows to introduce the further notions:
\begin{definition}
\end{definition}
EULER'S POINT OF  $ \nu $ :
\begin{equation}
    {\mathcal{P}}_{Euler}( \nu ) \; := \; c_{3}^{- 1} (  \frac{ \nu }{ 
\nu_{ref} })
\end{equation}
\begin{definition}
\end{definition}
EULER'S SPACE:
\begin{equation}
    {\mathcal{S}}_{Euler} \; := \; {\mathcal{P}}_{Euler}( 
{\mathcal{N}}_{Euler} ) \; = \; {\mathbb{Q}}^{3}
\end{equation}
    One has clearly that:
\begin{theorem}
\end{theorem}
\begin{enumerate}
    \item $ {\mathcal{S}}_{Euler} $ is a linear space on the field
    $ {\mathbb{Q}}  $
    \item
\begin{equation}
   {\mathcal{P}}_{Euler}( \nu_{R} ) \; = \; ( 0 \, , \, 0 \, , \,
   0 )
\end{equation}
\begin{equation}
    {\mathcal{P}}_{Euler}( 2 \nu_{R} ) \; = \; ( 1 \, , \, 0 \, , \,
   0 ) \; =: \; \hat{o}
\end{equation}
\begin{equation}
    {\mathcal{P}}_{Euler}( \frac{3}{2} \nu_{R} ) \; = \; ( 0 \, , \, 1 \, , 
\,
   0 ) \; =: \; \hat{q}
\end{equation}
\begin{equation}
   {\mathcal{P}}_{Euler}( \frac{5}{4} \nu_{R} ) \; = \; ( 0 \, , \, 0 \, , 
\,
   1 ) \; = \;  \hat{t}
\end{equation}
\end{enumerate}
The passage from points in Euler's space  to the associated pitchs
can, by construction, be obtained through the following:
\begin{theorem} \label{th:from Euler points to pitchs}
\end{theorem}
FROM EULER POINTS TO PITCHS:
\begin{equation}
    {\mathbf{p}} ( \nu ) \; = \; \frac{1200}{\log(2)} H_{prime}
    \, \cdot \, ( {\mathcal{P}}_{Euler} ( \nu ) ) \; \;
\end{equation}
where:
\begin{equation}
    H_{prime} \; := \;  (  \log (2) , \log(3) , \log(5) )
\end{equation}
is called the \emph{prime vector}

\smallskip
\begin{remark}
\end{remark}
ON THE PLANE ORTHOGONAL TO THE PRIME VECTOR IN THE EXTENDED EULER
SPACE

Introduced the:
\begin{definition}
\end{definition}
EXTENDED EULER SPACE:
\begin{equation}
    {\mathcal{S}}_{Euler}^{EXT} \; := \; \{ x \hat{o} + y \hat{f}
    + z \hat{t} \, \, , x,y,z \in {\mathbb{R}} \} \; = \; {\mathbb{R}}^{3}
\end{equation}

the theorem\ref{th:from Euler points to pitchs} implies,
obviously, that:
\begin{corollary} \label{cor:on adding vectors orthogonal to the prime 
vector}
\end{corollary}
\begin{equation}
   {\mathbf{p}} ( \nu ) \; = \;  H_{prime} \, \cdot \,
   ( {\mathcal{P}}_{Euler}(\nu) + E ) \, \; \forall E \in
   H_{prime}^{\perp} \, , \, \forall \nu  \in {\mathcal{N}}_{Euler}
\end{equation}
Since as a function from $ {\mathbb{R}}_{+} $ to $ {\mathbb{R}} $
$ {\mathbf{p}} ( \nu )  $ is obviously an injective function, one
could get a moment of surprise as to  the corollary\ref{cor:on
adding vectors orthogonal to the prime vector} and erroneously be
led to think that an Euler note is represented by many Euler
points.

Corollary\ref{cor:on adding vectors orthogonal to the prime
vector} is, instead, simply nothing but a consequence  of the
lemma\ref{lem:on the injectivity of Euler's coordinations} stating
that Euler's coordination, defined on the extended Euler's space,
is not injective

\smallskip

\begin{definition}  \label{def:canonical notes' basis of Euler space}
\end{definition}
CANONICAL NOTES' BASIS OF EULER SPACE:
\begin{equation}
    {\mathbb{E}}_{Euler}^{notes} \; := \; \{\hat{o} \, , \,\hat{f} \, , \, 
\hat{t}  \}
\end{equation}

The chosen names for  versors in the definition\ref{def:canonical
notes' basis of Euler space} is owed to the fact that:
\begin{enumerate}
    \item $ \hat{o} $ is the Euler point of the note being distant
    an interval of octave from $ \nu_{R} $
    \item $ \hat{f} $ is the Euler point of the note being distant
    an interval of perfect fifth  from $ \nu_{R} $ (both in
    pytagoric and in just tunings, concepts that we will introduce later)
    \item $ \hat{t} $ is the Euler point of the note being distant
    an interval of major third  from $ \nu_{R} $ (in just tuning, a concept 
we will introduce later)
\end{enumerate}
Let us then consider some following important subsets of Euler
notes and the corresponding subsets of representative Euler points
in Euler's space:
\begin{definition} \label{def:just-tuned Euler notes}
\end{definition}
JUST-TUNED EULER NOTES:
\begin{equation}
    {\mathcal{N}}_{Euler}^{just-tuned} \; := \; \{ \nu_{ref} \,  c_{3} ( 
\vec{x} )  \, , \, \vec{x} \in {\mathbb{Z}}^{3} \}
\end{equation}
\begin{definition}
\end{definition}
JUST-TUNED SUBSET OF EULER'S SPACE:
\begin{equation}
  {\mathcal{S}}_{Euler}^{just-tuned} \; := \; {\mathcal{P}}_{Euler} ( 
{\mathcal{N}}_{Euler}^{just-tuned}
  ) \; = \; \{ x_{o} \hat{o} \, + \, x_{f} \hat{f} \, + \, x_{t}
  \hat{t} \, : \,  x_{o} , x_{f} , x_{t}  \in {\mathbb{Z}} \}
\end{equation}

\begin{definition} \label{def:pytagorically-tuned Euler notes}
\end{definition}
PYTAGORICALLY-TUNED EULER NOTES:
\begin{equation}
   {\mathcal{N}}_{Euler}^{Pyt-tuned} \; := \; \{ \nu_{ref} \,  c_{3} ( x_{o} 
\hat{o} \, + \, x_{f} \hat{f} )  \, , \, x_{o} , x_{f}  \in {\mathbb{Z}}  \}
\end{equation}
\begin{definition}
\end{definition}
PYTAGORICALLY-TUNED SUBSET OF EULER'S SPACE:
\begin{equation}
    {\mathcal{S}}_{Euler}^{Pyt-tuned} \; := \; {\mathcal{P}}_{Euler} ( 
{\mathcal{N}}_{Euler}^{Pyt-tuned}
    ) \; = \; \{ x_{o} \hat{o} \, + \, x_{f} \hat{f} \, : \,  x_{o} , x_{f}  
\in {\mathbb{Z}}   \}
\end{equation}

Given an integer number $ n \in {\mathbb{N}} $:
\begin{definition}
\end{definition}
n-TEMPERED TUNED EULER NOTES:
\begin{equation}
    {\mathcal{N}}_{Euler}^{n-temp-tuned} \; := \; \{ \nu_{ref} \,  c_{3} (  
\frac{x_{o}}{n} \hat{o}  )  \, , \, x_{o} \in {\mathbb{Z}}  \}
\end{equation}
\begin{definition}
\end{definition}
n-TEMPERED-TUNED SUBSET OF EULER'S SPACE:
\begin{equation}
    {\mathcal{S}}_{Euler}^{w-temp-tuned} \; := \;  {\mathcal{P}}_{Euler} (
    {\mathcal{N}}_{Euler}^{w-temp-tuned}) \; = \; {\mathbb{Z}}
    \frac{1}{n} \hat{o}
\end{equation}
Given three integer numbers $ n_{1} \, , \,n_{2} \, , \, n_{3} \in
{\mathbb{N}} $:
\begin{definition}
\end{definition}
$ n_{1},n_{2},n_{3} $-TEMPERED-JUST-TUNED EULER NOTES
\begin{equation}
    {\mathcal{N}}_{Euler}^{n_{1},n_{2},n_{3}-temp-just-tuned} \;
    := \; \{ \nu_{ref} \,  c_{3} (  \frac{x_{o}}{n_{1}} \hat{o} \, + \, 
\frac{x_{f}}{n_{2}} \hat{f} \, + \, \frac{x_{t}}{n_{3}} \hat{t}  )  \, , \, 
x_{o} , x_{f} , x_{t} \in {\mathbb{Z}}  \}
\end{equation}
\begin{definition}
\end{definition}
$ n_{1},n_{2},n_{3} $-TEMPERED-JUST-TUNED SUBSET OF EULER'S SPACE:
\begin{equation}
    {\mathcal{S}}_{Euler}^{n_{1},n_{2},n_{3}-temp-just-tuned} \;
    := \;  {\mathcal{P}}_{Euler}  (  
{\mathcal{N}}_{Euler}^{n_{1},n_{2},n_{3}-temp-just-tuned} )
\end{equation}

\medskip

It is important, at this point, to observe that  Euler's space $
{\mathcal{S}}_{Euler} $ may be used not only to denote
\emph{notes} but \emph{interval among notes}. With this respect is
is then useful to introduce the following:

\begin{definition}  \label{def:canonical intrevals' basis of Euler space}
\end{definition}
CANONICAL INTERVALS' BASIS OF EULER SPACE:
\begin{equation}
    {\mathbb{E}}_{Euler}^{int} \; := \; \{\hat{o} \, , \,\hat{f} - \hat{o} 
\, , \, \hat{t} - 2 \hat{o}  \}
\end{equation}
whose elements represents, respectively, the \emph{octave
interval}, the \emph{perfect fifth interval} and the \emph{major
third interval}.

Let us then introduce two  intervals whose theoretical relevance
will soon appear:
\begin{definition} \label{def:fifth commma}
\end{definition}
FIFTH COMMA (PYTAGORICAL COMMA):
\begin{equation}
    \hat{Kf} \; := \; - \, 7 \hat{o} \, + \,  12 ( \hat{f} - \hat{o} )
\end{equation}
\begin{definition} \label{def:third commma}
\end{definition}
THIRD COMMA (SYNTONIC COMMA):
\begin{equation}
    \hat{Kt} \; := \; 2 \hat{o} \, - \, 4  ( \hat{f} - \hat{o} ) \, + \,
( \hat{t} - 2 \hat{o} )
\end{equation}
Applying the theorem\ref{th:from Euler points to pitchs} one finds
that:
\begin{equation}
    {\mathbf{Kf}} \, := {\mathbf{p}}( {\mathcal{P}}_{Euler}^{- 1} ( \hat{Kf} 
)) \; \approx  \;
    23.46 \, Cents
\end{equation}
\begin{equation}
    {\mathbf{Kf}} \, := {\mathbf{p}}( {\mathcal{P}}_{Euler}^{- 1} ( \hat{Kf} 
)) \; \approx  \;
    - 21.61 \, Cents
\end{equation}

\smallskip

Let us now analyze how all the previously discussed modulo-octave
stuff appears in the framework of Euler space.

The computations about Euler's space may be performed using he
following Mathematica expression of the mentioned notebook of
sec.\ref{sec:Mathematica's notebook Mathematical Music Toolkit}:
\begin{verbatim}

eulercoordination[x_List]:=
  Product[Power[Prime[i],Part[x,i]], {i,1,Length[x]}]

FROMeulerpointTOnote[eulpoint_List] :=
  referencenote*2^Part[eulpoint,1]*3^Part[eulpoint,2]*5^Part[eulpoint,3]

FROMnoteTOpitch[nu_]:= Times[1200,Power[Log[2],-1]]
*Log[nu]-Log[Times[nu,Power[referencenote,-1]]]

FROMpitchTOnote[pitch_]:=
referencenote*Exp[Times[Log[2],Power[1200,-1]]*pitch]

Hprime=Table[Log[Prime[i]],{i,1,3}];

FROMeulerpointTOpitch[eulpoint_List]:=
Times[1200,Power[Log[2],-1]]*Dot[eulpoint,Hprime]

lastpartisequaltosomethingQ[x_,something_]:=Equal[Last[x],something]

rationalsupto[n_]:=Flatten[ Table[Rational[i,j],{i,0,n},{j,1,n}] ]

FROMindexTOrational[index_,n_]:=Part[rationalsupto[n],index]

FROMlistofindicesTOlistofrationals[list_,n_]:=
  Table[FROMindexTOrational[Part[list,i],n],{i,1,Length[list]}]

FROMnoteTOeulerpoint[nu_,n_]:=
  FROMlistofindicesTOlistofrationals[
    First[generalizedselect[
        Flatten[Table[ {{i,j,k},
              Dot[{Part[ rationalsupto[2],i ],Part[ rationalsupto[2],j],
                  Part[ rationalsupto[2],k ]},Hprime]},{i,1,
              Length[ rationalsupto[2] ]},{j,1,Length[ rationalsupto[2] 
]},{k,
              1,Length[ rationalsupto[2] ]}],2]  ,
        lastpartisequaltosomethingQ,FROMnoteTOpitch[nu]]],n]

FROMwordTOlistofeulerpoints[w_]:=
  Table[FROMnoteTOeulerpoint[Part[FROMwordTOscale[w] ,i]],{i,1,Length[w]}]
octaveepoint={1,0,0}; fifthpoint={0,1,0}; thirdpoint={0,0,1};

canonicalnotesbasis:={octaveepoint,fifthpoint,thirdpoint}

canonicalintervalsbasis:={octaveepoint,fifthpoint-octaveepoint,
    thirdpoint-2octaveepoint}

\end{verbatim}

Euler's formalization of the consonance's issue was based on
the introduction of the following function:

given a positive integer number $ x \,  \in \, {\mathbb{N}}_{+} $,
the Theorem of Prime Factorization assures the existence and the
unicity of a sequence $  \{ e_{n} \}_{n \in {\mathbb{N}}_{+}} $ of
positive integer numbers such that:
\begin{enumerate}
    \item $  \{ e_{n} \}_{n \in {\mathbb{N}}_{+}} $ ends with a
    countable infinity of zeros, i.e.:
\begin{equation}
    \exists \, N \, \in \, {\mathbb{N}}_{+} \; : ( \, \; e_{n} \, = \,
    0 \; \forall n > N \,)
\end{equation}
    \item
\begin{equation}
    x \; = \; \prod_{n=1}^{\infty} \pi(n)^{e(n)}
\end{equation}
\end{enumerate}
\begin{definition} \label{def:gradus suavitatis of a positive integer}
\end{definition}
GRADUS SUAVITATIS OF  THE  POSITIVE INTEGER x:
\begin{equation}
    G_{Eul}(x) \; := \; 1 \, + \, \sum_{k=1}^{\infty} e_{k} ( \pi(k) \, - \, 
1 )
\end{equation}
Given a positive rational number $ r \in {\mathbb{Q}}_{+} $ we
have by definition that:
\begin{equation}
    \exists ! ( p , q ) \in {\mathbb{N}}_{+}^{2} \; : \; ( \, x \, = \,
    \frac{p}{q} \, and \, gcd(p , q) = 1 \, )
\end{equation}
We can consequentially generalize the definition \ref{def:gradus
suavitatis of a positive integer} in the following way:
\begin{definition} \label{def:gradus suavitatis of a positive rational}
\end{definition}
GRADUS SUAVITATIS OF  THE  POSITIVE RATIONAL x:
\begin{equation}
   G_{Eul}(x) \; := \;  G_{Eul} ( p \, \cdot \, q )
\end{equation}
The definition \ref{def:gradus suavitatis of a positive rational}
may be then applied to the ratio of Just-tuned Euler notes
resulting in the following function
\begin{definition} \label{def:gradus suavitatis of a bichord of  just tuned 
Euler notes}
\end{definition}
GRADUS SUAVITATIS OF  BICHORDS OF JUST TUNED EULER NOTES:

the map $  G_{Eul} \, : \, {\mathbb{N}}_{Euler}^{just-tuned} \,
\times  \, {\mathbb{N}}_{Euler}^{just-tuned}  \, \mapsto \,
{\mathbb{N}}_{+} $ :
\begin{equation}
  G_{Eul} ( \mu \, , \, \nu     ) \; := \; G_{Eul} ( \frac{ \nu }{
  \mu} )
\end{equation}

As it has been observed by Mazzola, since the higher is the
consonance the lower is the value of the gradus function, $
G_{Eul} $ should be more properly called gradus dissuavitatis:

the total ordering $ <_{Eul} $ of dissonant intervals induced by
definition \ref{def:gradus suavitatis of a bichord of  just tuned
Euler notes}:
\begin{multline}
    1 \; <_{Eul} \; 5P \; <_{Eul} \;
    4P \; <_{Eul} \; 3M \; =_{Eul} \; 6M  \;  <_{Eul} \; 3m \\
     =_{Eul} \; 6m \;  =_{Eul} \; 2M \; <_{Eul} \; 7m  \; <_{Eul} \\ \; 7M  
\;
<_{Eul} \; 2m \; <_{Eul} \; 5^{\circ}
\end{multline}
where I have adopted the standard musical notation for fourth and
fifth perfect (P) as well as for major (M) , minor (m) and
diminished ($ \circ $)  intervals; this can been directly verified
through the following expression of the mentioned Mathematica
notebook of section\ref{sec:Mathematica's notebook Mathematical
Music Toolkit}:
\begin{verbatim}
gradussuavitatis[n_Integer]:=
  1+Sum[  ( FactorInteger[n][[i]][[1]]-1)* FactorInteger[n][[i]][[2]]  , 
{i,1,
        Length[FactorInteger[n]]}]

gradussuavitatis[r_Rational]:=
  gradussuavitatis[
    Times[Numerator[r]*Denominator[r],
      Power[GCD[Numerator[r],Power[Denominator[r],2]] , -1 ]]]
\end{verbatim}

\smallskip

The next step in the history of the comprehension of
\emph{physical consonance} was made in 1877 by Hermann L. F.
Helmhotz in his masterpiece \cite{Helmholtz-54} whose conclusions
are condensed in the following passage:

\begin{center}
  \emph{"When two musical tones are sounded at the same time, their united 
sound is generally disturbed by beats of the upper partials,
   so that a greater or less part of the whole mass of sound is broken up 
into pulses of tone, and the joint effect is rough.
   These exceptional cases are cold \textit{Consonances}. From the Part II: 
"On the interruptions of harmony", cap. 10: "Beats of the Upper Partial 
Tones", section "Order of Consonances with respect to Harmoniousness", pagg. 
194-197 of \cite{Helmholtz-54}}
\end{center}
\begin{center}
\emph{"When two or more compound tones are sounded at the same time beats 
may arise from the combinational tones as well
from the harmonic upper partials. In Chapter VII it was shown that the 
loudest combinational tones resulting from two
generating tones is that corresponding to the difference of their pitch 
numbers, or the differential tone of the first order. It is this 
combinational tone, therefore, which is chiefly effective in producing 
beats" From the Part II: "On the interruptions of harmony", cap.11: "Beats 
due to combinational tones" pagg.197-198 of \cite{Helmholtz-54} }
\end{center}

To explain Helmholtz's point of view it is essential to
distinguish among \emph{interference phenomena} prescribed by
Linear Acoustics and the \emph{nonlinear combinational effects}
owed to violations of the axiom \ref{ax:principle of
superposition}

let us suppose to have two distinct sonore sources $ i_{1}(t) $
and $ i_{2}(t) $ suct that each one, from its own, would produce,
if it was alone, respectively the sounds:
\begin{equation}
    s_{1} (t) \; :=  \; {\mathcal{R}}_{ac} [  i_{1}(t)  ] \; = \;  a_{1} 
e^{i \omega_{1} t}
\end{equation}
and
\begin{equation}
    s_{2} (t) \; :=  \; {\mathcal{R}}_{ac} [  i_{2}(t)  ] \; = \;  a_{2} 
e^{i \omega_{2} t}
\end{equation}
\begin{description}
    \item[LINEAR REGIME]
Assuming the linearity of medium's acoustic response, i.e. the
axiom \ref{ax:principle of superposition} it follows that:
\begin{multline}
    s(t) \; = \; {\mathcal{R}}_{ac} [ i_{1}(t) + i_{2}(t) ] \; =
    \; = {\mathcal{R}}_{ac} [ i_{1}(t) ] \, + \, {\mathcal{R}}_{ac} [
    i_{2}(t) ] \\
     = \; s_{1} (t) + s_{2} (t) \; = a_{1} e^{i \omega_{1} t} +  a_{2} e^{i 
\omega_{2} t}
\end{multline}
Posed:
\begin{equation}
    \Delta \, \omega \; := \; \omega_{2} \, - \, \omega_{1}
\end{equation}
and:
\begin{eqnarray}
  a_{1} &=& | a_{1} | e^{ i \phi_{1}} \\
   a_{2} &=& | a_{2} | e^{ i \phi_{2}}
\end{eqnarray}
it follows that:
\begin{equation} \label{eq:beats1}
    s(t) \; = \;   a(t) \,   e^{i \omega_{1} t}
\end{equation}
where:
\begin{equation} \label{eq:beats2}
    | a (t) | \; = \; \sqrt{ | a_{1} |^{2} \,  +  \, | a_{2} |^{2} \, + \, 2 
  | a_{1} |  \, | a_{2} | \, \cos (  \phi_{1} - \phi_{2} -  \Delta \, \omega 
)  )  }
\end{equation}
and:
\begin{equation} \label{eq:beats3}
    Arg( a (t) ) \; = \; \arctan [ \frac{ | a_{1} | \sin ( \phi_{1} ) \,  +  
\, | a_{2} | \sin ( \phi_{2} +  \Delta \, \omega \, t )  }{ | a_{1} | \cos ( 
\phi_{1} ) \,  +  \, | a_{2} | \cos ( \phi_{2} +  \Delta \, \omega \, t )  } 
      ]
\end{equation}
The qualitative behavior prescribed by equation
eq.\ref{eq:beats1}, eq.\ref{eq:beats2} and eq.\ref{eq:beats3} is
with some respect similar to that of an \emph{amplitude
modulation} in which an harmonic \emph{carrier signal} of
pulsation $ \omega_{1} $ is modulated by an harmonic
\emph{envelope signal} of pulsation $ \Delta \omega $.

Though similar, the signal described by eq.\ref{eq:beats1},
eq.\ref{eq:beats2} and eq.\ref{eq:beats3} is not an
\emph{amplitude modulating signal} that would have not only
spectral component $ \omega_{1} \, + \, \Delta \omega $ but also
spectral component $ \omega_{1} \, - \, \Delta \omega $.

    \item[NON LINEAR REGIME]
\end{description}
Since, usually the  acoustical response of the medium may be
expressed as the sum of a linear and nonlinear part, one has that:
\begin{multline}
   s(t) \; = \;  {\mathcal{R}}_{ac} [  i_{1}(t) + i_{2}(t)  ] \; =
    \; {\mathcal{R}}^{L}_{ac}
    [  i_{1}(t) + i_{2}(t) ] \, + \, {\mathcal{R}}^{NL}_{ac}
    [  i_{1}(t) + i_{2}(t) ]    \\
    = {\mathcal{R}}^{L}_{ac} [  i_{1}(t) ] \, + \, {\mathcal{R}}^{L}_{ac}
    [ i_{2}(t) ] \, + \, {\mathcal{R}}^{NL}_{ac} [  i_{1}(t) + i_{2}(t)
    ]
\end{multline}
It may be worth observing that, in the nonlinear regime, one
has to give up the complex representation of harmonic motion
since:
\begin{equation}
    \Re {\mathcal{R}}^{NL}_{ac} [  \tilde{i} (t) ] \; \neq \; 
{\mathcal{R}}^{NL}_{ac}
    [i(t)]
\end{equation}
\emph{Combinational tones} arise when the non linear part of the
acoustical response function is given by the product of the two
input signals.

One has that
\begin{multline}
      {\mathcal{R}}^{NL}_{ac} [  i_{1}(t) + i_{2}(t)  ] \; = \; 
{\mathcal{R}}^{NL}_{ac} [  i_{1}(t)
     ] \, \cdot \, {\mathcal{R}}^{NL}_{ac} [  i_{2}(t) ] \\
     = \;  s_{1} (t) \,  \cdot \,   s_{2} (t) \; = a_{1} \cos ( \omega_{1} t 
  + \phi_{1} ) \, \cdot \,  a_{2} \cos ( \omega_{2} t  + \phi_{2}
     ) \\
     = \; \frac{ a_{1} a_{2} }{2} \cos [ ( \omega_{1} + \omega_{2}
     ) + ( \phi_{1} + \phi_{2} ) ] \, + \, \frac{ a_{1} a_{2} }{2}
     \cos [ ( \omega_{1} - \omega_{2}
     ) + ( \phi_{1} -  \phi_{2} ) ]
\end{multline}
Consequentially the net acoustic response consists not only of the
two notes with pulsations $ \omega_{1} $ and $ \omega_{2} $ but
also of other 2 notes with pulsations $ \omega_{1} \, + \,
\omega_{2} $ and $ \omega_{1} \, - \, \omega_{2} $.
    \item[HORIZONTAL RULES IN SCORES] The theory of dissonance resolution 
inside a fixed tonality and the theory of modulation
    ruling how to connect different tonalities
\end{description}

\newpage
\section{Vertical rules require the commensurability of notes}
Given two sounds:
\begin{equation}
    s_{1} (t) \; = \;  \sum_{n = - \infty}^{ + \infty } a_{n}^{(1)} e^{i n 
\omega_{1} t}
\end{equation}
\begin{equation}
    s_{2} (t) \; = \; \sum_{n = - \infty}^{ + \infty } a_{n}^{(2)} e^{i n 
\omega_{2} t}
\end{equation}
let us formalize the original Galilei's viewpoint mixed with a bit
of Helmotz, and let us consider the etymology of the term
consonance: we are led to select the harmonics of the two sounds
that are commensurable.
\begin{definition} \label{def:physical index of consonance of two sounds}
\end{definition}
PHYSICAL INDEX OF CONSONANCE AMONG $ s_{1} $ and  $  s_{2}  $:
\begin{equation}
    I[ s_{1} \, , \, s_{2} ] \; := \; \sum_{n = - \infty}^{+
    \infty}  \sum_{m = - \infty}^{+
    \infty} \, \delta_{Kronecker} ( n  \omega_{1} + m \omega_{2} ) \,
    ( a_{n}^{(1)} + a_{m}^{(2)} )
\end{equation}
where:
\begin{equation}
    \delta_{Kronecker} ( \omega ) \; := \; \left\{%
\begin{array}{ll}
    1, & \hbox{if $ \omega \, = \, 0 $ ;} \\
    0, & \hbox{otherwise;} \\
\end{array}
\right.
\end{equation}
\begin{definition}
\end{definition}
$ \omega_{1} $ AND $ \omega_{2} $ ARE COMMENSURABLE ( $
\omega_{1}\, \sim_{ {\mathbb{Q}}} \, \omega_{2} $ ):
\begin{equation}
    \frac{ \omega_{1}  }{ \omega_{2} } \; \in \; {\mathbb{Q}}
\end{equation}

\begin{theorem}
\end{theorem}
CONSONANCE REQUIRES THE COMMENSURABILITY OF  2 NOTES:
\begin{equation}
  \omega_{1} \, \nsim_{ {\mathbb{Q}}} \, \omega_{2} \; \Rightarrow
  \; I [ s_{1} , s_{2} ] \, = \, 0
\end{equation}

More generally, given   N sounds $ s_{1} \,  \cdots \, s_{N} $:
\begin{equation}
    s_{i} (t) = \sum_{n = - \infty}^{ + \infty } a_{n}^{(i)} e^{i n 
\omega_{i}
    t} \; i = 1 , \cdots \, N
\end{equation}
\begin{definition}  \label{def:physical index of consonance of N sounds}
\end{definition}
PHYSICAL INDEX OF CONSONANCE AMONG $ s_{1} \,  \cdots \, s_{N} $:
\begin{equation}
    I[ s_{1} \, , \cdots \, s_{N} ] \; := \; \sum_{n_{1}, \cdots \, n_{N} = 
- \infty}^{+
    \infty}  \delta_{Kronecker} ( \sum_{i=1}^{N} \omega_{i} n_{i}  ) \,
    ( \sum_{i=1}^{N} a_{n}^{(i)})
\end{equation}
\begin{definition}
\end{definition}
$ \omega_{1} \, \cdots \,  \omega_{N} $ ARE COMMENSURABLE ( $
\omega_{1} \, \sim_{ {\mathbb{Q}}} \, \cdots \, \sim_{
{\mathbb{Q}}}  \omega_{N} $ ):
\begin{equation}
    \exists n_{1} \, \cdots \,  n_{N} \in {\mathbb{Z}} \; : \; 
\sum_{i=1}^{N} \omega_{i}
    n_{i} \,  = \,  0
\end{equation}
\begin{theorem}
\end{theorem}
CONSONANCE REQUIRES THE COMMENSURABILITY OF  N NOTES:
\begin{equation}
  \omega_{1} \, \nsim_{ {\mathbb{Q}}} \, \cdots \, \nsim_{
{\mathbb{Q}}}  \omega_{N}  \; \Rightarrow
  \; I [ s_{1} , \cdots , s_{N} ] \, = \, 0
\end{equation}

\begin{definition} \label{def:musical instrument}
\end{definition}
MUSICAL INSTRUMENT:
\begin{equation}
   \bar{a} \; = \; \{ a_{n} \}_{ n \in {\mathbb{Z}}} \; \in \;
   {\mathbb{C}}^{\infty}  \, : \, \sum_{n = - \infty} ^{\infty}
   a_{n} a_{-n} \, < \, + \infty
\end{equation}
Given N notes $ \omega_{1} \, \cdots \,  \omega_{N} $ and a
musical instrument $ \bar{a} $:
\begin{definition}  \label{def:physical index of consonance of N notes 
w.r.t. a musical instrument}
\end{definition}
PHYSICAL INDEX OF CONSONANCE OF $ \omega_{1} \, \cdots \,
\omega_{N} $ W.R.T.  $ \bar{a} $:
\begin{equation}
    I (  \omega_{1} \, \cdots \,  \omega_{N} \,  | \,  \bar{a} )
    \; := \; I (  \sum_{n = - \infty}^{ + \infty } a_{n} e^{i n \omega_{1}
    t} \, , \,  \cdots \, , \, \sum_{n = - \infty}^{ + \infty } a_{n} e^{i n
    \omega_{n} t} )
\end{equation}
\begin{example} \label{ex:physical consonance of notes w.r.t. a pure 
oscillator}
\end{example}
PHYSICAL CONSONANCE OF NOTES W.R.T. A OF PURE OSCILLATOR:

A pure oscillator is, by definition, a musical instrument of the
form:
\begin{equation}
    a_{n} \; := \;   a \, \delta_{n,k}
\end{equation}
Then:
\begin{equation}
  I ( \omega_{1} , \omega_{2} \, | \, a_{n} )  \; = \; \left\{%
\begin{array}{ll}
    2 \, a  , & \hbox{if $ \omega_{1} \, = \, \omega_{2}$} \\
    0, & \hbox{otherwise.} \\
\end{array}%
\right.
\end{equation}
\begin{example} \label{ex:physical consonance of notes w.r.t. an ideal  
instrument}
\end{example}
PHYSICAL CONSONANCE OF NOTES W.R.T. AN IDEAL MUSICAL INSTRUMENT

The example\ref{ex:physical consonance of notes w.r.t. a pure
oscillator} shows that, as to the possibility of generating
consonants sounds, a pure oscillator is the most inefficient
musical instrument: the only consonant notes it can product are
unisonous ones.

Looking at definition\ref{def:physical index of consonance of N
notes w.r.t. a musical instrument} it should be clear why,
contrary, an ideal instrument, i.e. the most efficient  musical
instrument as to production of consonance, is one in which the
amplitudes of higher harmonics decades in the slowest possible
way; we are then led to define an ideal musical instrument by the
condition:
\begin{equation}
    a_{n} \; := \; \frac{a}{n}
\end{equation}
Then:
\begin{equation}
    I ( \omega_{1} , \omega_{2} \, | \, a_{n} ) \; = \; \left\{%
\begin{array}{ll}
     + \infty , & \hbox{if $ \omega_{1} \sim_{ {\mathbb{Q}}} \omega_{2} $  
;} \\
    0, & \hbox{otherwise.} \\
\end{array}%
\right.
\end{equation}
We would like to advise the reader that the  wrong expression $
\frac{n+m}{ n m } $ given in the $ 6^{th} $ chapter of
\cite{Frova-02} for the physical index of consonance among two
commensurable notes with ratio of frequencies $ \frac{m}{n} $ with
$ gcd( m , n) = 1 $ is based on the erroneous procedure of
defining I up to an infinite constant (the harmonic serie).
\newpage
\section{Horizontal rules require the incommensurability of notes}
\label{sec:Horizontal rules require the incommensurability of
notes}

The theory of dissonance resolution inside a fixed tonality and
the theory of modulation ruling how to connect different
tonalities are based on the same central structural element: the
\emph{scales at fixed interval}:

given two positive real numbers $ \omega_{1} \, , \, R \; \in \;
{\mathbb{R}}_{+} $:
\begin{definition}
\end{definition}
SCALE OF $ \omega_{1} $ AT FIXED INTERVAL R:

$ scale ( \omega_{1} \, , \, R ) \; := \; \{ \omega_{n}
\}_{n=1}^{N} $ :
\begin{eqnarray}
  \omega_{n}\; & = & \; R \omega_{n-1} \,  mod ( 2^{{\mathbb{N}}} )  \; \; n 
= 1 , \cdots , N \\
    N & \; := \; & \; \left\{%
\begin{array}{ll}
    \min \{ n \in {\mathbb{N}} \: : \: \omega_{N} \,  =  \, \omega_{1} \, 
mod ( 2^{{\mathbb{N}}} )   \}, & \hbox{if $ \min \{ \cdots \} $  exists} \\
    + \infty , & \hbox{otherwise} \\
\end{array}%
\right.
\end{eqnarray}
\begin{definition}
\end{definition}
$ scale ( \omega_{1} \, , \, R ) $ IS FINITE
\begin{equation}
    N \; < \; + \infty
\end{equation}
We have that:
\begin{theorem} \label{th:finiteness versus the irrationality of the fixed 
interval}
\end{theorem}
FINITENESS VERSUS THE IRRATIONALITY OF THE FIXED INTERVAL:

\begin{equation}
    scale ( \omega_{1} \, , \, R) \text{ is finite} \; \Leftrightarrow \; R 
\notin {\mathbb{Q}}
\end{equation}
\begin{proof}
The thesis follows immediately by the fact that:
\begin{equation}
     \sqrt[n]{2} \, \notin \, {\mathbb{Q}} \; \; \forall n \in
     {\mathbb{N}} \, : \, n \geq 2
\end{equation}
\end{proof}
\begin{corollary}
\end{corollary}
FINITENESS VERSUS THE INCOMMENSURABILITY OF NOTES:
\begin{equation}
     scale ( \omega_{1} \, , \, R) \text{ is finite} \;  \Leftrightarrow
    \; \omega_{i} \nsim_{ {\mathbb{Q}}} \omega_{j} \, i \neq  j < N
\end{equation}
\begin{proof}
It is sufficient to verify that $ \sim_{{\mathbb{Q}}} $ is an
equivalence relation and hence, in particular, is transitive
\end{proof}

Let us now\ observe that the notes in a scale at fixed interval $
scale ( \omega_{1} \, , \, R ) $ don't occur in order of
increasing frequency; given a string or a  sequence of real
numbers $ \{ a_{n} \} $
\begin{definition}
\end{definition}
ORDERING OF $  \{ a_{n} \}  $:
\begin{equation}
  ord [  \{ a_{n}
\}  ] \; = \; \{ \tilde{a}_{n} \}  \; := \;  permutation (  \{
a_{n} \}  ) \, : \,
   i < j  \;  \Rightarrow   \;   \tilde{a}_{i} \, < \,  \tilde{a}_{j}
\end{equation}
\begin{example}
\end{example}
PYTAGORIC SCALE OF $ \omega_{1} $:

The pytagoric scale of a note is constructed following the
\emph{natural cycle of fifths}, i.e. the cycle of fifths assuming
as fifth interval the one specified by the second harmonic of $
\omega_{1}  $ equal to  3:
\begin{equation}
    pytagoric( \omega_{1} ) \; := \;  scale ( \omega_{1} \, , \, 3 )
\end{equation}
By theorem \ref{th:finiteness versus the irrationality of the
fixed interval} it immediately follows that pytagoric scales are
always infinite.
\begin{remark} \label{rem:ordering the pytagoric scales globally or at fifth 
cycles' blocks}
\end{remark}
ORDERING THE PYTAGORIC SCALES GLOBALLY OR AT FIFTH CYCLES BLOCKS

Let us observe that:
\begin{equation}
  \cdot_{n \in {\mathbb{N}}} ord(scale(\nu_{ref} ,3)_{n,n+11}) \; \neq
\; ord(scale( \nu_{ref} , 3))
\end{equation}
a fact, this one, we will see more concretelly in the section
\ref{sec:Horizontal rules require the incommensurability of
notes}.

\bigskip

\begin{example}
\end{example}
N-EQUALLY-TEMPERED SCALES OF  $ \omega_{1} $:

The N-equally-tempered scale of $ \omega_{1} $ is constructed
following the\emph{N-equally-tempered cycle of fifths}, i.e. the
cycle of fifth adopting as basic fifth interval a value different
from the one specified by  Nature through the third harmonic,
chosen in order of obtaining the periodicity of the cycle with
period N.

Indeed, defining:
\begin{equation}
   tempered_{eq}( \omega_{1} , N ) \; := \; scale ( \omega_{1} , 2^{ 
\frac{7}{N}
   })
\end{equation}
by theorem \ref{th:finiteness versus the irrationality of the
fixed interval} it immediately follows that  $ tempered(
\omega_{1} , N ) $ is a finite scale of N notes.

\bigskip

We will now apply these considerations to show how the
incommensurability of notes is required both by the theory of
physical dissonance's resolution inside a fixed tonality and  by
the other ingredient horizontal rules in scores are made of,
namely the theory of modulation between different tonalities.

At this purpose its is necessary, first of all, to introduce a
rigorous mathematical formalization of both Classical Harmony and
some generalizations of its.

Such an objective was pursued by the decennial work of the Guerino
Mazzola's research's group  about \emph{Mathematical Music
Theory}, \cite{Mazzola-01}, \cite{Mazzola-02a},
\cite{Mazzola-02b}.

Some preliminary consideration is required as to the  functorial
formalism in which the whole matter has been finally recasted:

the construction of the abstract mathematical  machinery of
\emph{forms} and \emph{denotators} has been justified by Mazzola
through the following argumentations:
\begin{enumerate}
    \item the temporally evolving nature of the global encyclopedic space
    of human knowledge (called by Mazzola \emph{encyclospace} with
    an explicit referring to the assumption of its informatic
    storing  in a virtual space such as a database embedded in Internet)
    possessing a topological structure endowed in the
    inter-relational links among distinct concepts (cfr. in particular the $
    5^{th}$ chapter "Navigation" and the $ 9^{th} $ chapter
    "Yoneda Perspectives"  of  \cite{Mazzola-02a})
    \item the requirement  to mirror such a structure into the
    structure of the software concretelly computing the
    mathematical concepts of encyclospace
\end{enumerate}
While the latter argumentation, based on the  well known role
played by Category Theory in Theoretical Computer Science
\cite{Asperti-Longo-91}, \cite{Amadio-Curien-98}, is clear
corresponding to the concrete architecture of the program RUBATO,
the former argumentation requires some further analysis:

while in \cite{Mazzola-02b} Mazzola explicitely states that the
existence and size of \emph{form semiotics} is indeed an essential
problem, observing that all the concretelly appearing forms of
Mathematical Music Theory, being regular, can be obtained from
\emph{simple forms} through transfinite recursion, the $ 9^{th} $
chapter "Yoneda Perspectives"  of  \cite{Mazzola-02a} (more or
less explicitely) states that the  passage to the functorial
formalism allows a consistent mathematical adoption of circular
definitions.

Mazzola's claim is, with this regard, trivially false:

the passage from the set-theoretic foundations to the
topos-theoretic foundations of Mathematics, corresponding to a
passage from the Zermelo-Fraenkel formal system augmented with the
Axiom of Choice, ZFC from here and beyond, to a formal system
equiconsistent with the weaker Restricted Zermelo formal system
augmented with the Axiom of Choice, RZC from here and beyond,
doesn't drop the Foundation Axiom banning any kind of circularity
(cfr. the $ 10^{th} $ section "Topos-Theoretic and Set-Theoretic
Foundations" of the $ 6^{th} $ chapter "Topoi and Logic" of
\cite{Mc-Lane-Moerdijk-92}):

the proof of the equiconsistence among the Topos- Theoretic
Foundation of Mathematics and the set-theoretical foundation
furnished by RZC, consists in:
\begin{enumerate}
\item showing that from each model $
{\mathcal{S}} $   of RZC one can construct a well-pointed topos $
{\mathcal{E}} $ (more precisely a model of the first-order-theory
WPT of well-pointed topoi) with choice and natural number object

\item showing, conversely, that for each such well-pointed topos $
{\mathcal{E}} $ one can construct a model of RZC.

\end{enumerate}

The first part of the equi-consistence's proof constructing a
well-pointed  topos $ {\mathcal{E}} $ from a model $ {\mathcal{S}}
$ of set-theory consist simply into the construction of the
category of all sets of $ {\mathcal{S}} $.

The second part of the equi-consistence's proof is, instead, more
complex.

My purpose, with this regard,  is to remark how, in  such a
construction of a model $ {\mathcal{S}} $ of RCZ from a
well-pointed topos with natural numbers object and choice, the
Vicious Circle Principle (adopting the old terminology by Russell
and Poincar\'{e} concerning the exorcism of Russell's paradox for
whose discussion I demand to \cite{Longo-00} and to the section
$4.3$ "The Set-Theoretical Hierarchy" of \cite{Odifreddi-89})
expressed by the following one among RZC's axioms:

\begin{axiom}  \label{ax:axiom of foundation}
\end{axiom}
AXIOM OF FOUNDATION
\begin{equation*}
    ( \exists y ) ( y \in x) \; \rightarrow \; ( \exists y ) [ y \in
    x \, \wedge \, ( \forall z \in y ) ( z \notin x ) ]
\end{equation*}

and stating that any nonempty set has a minimal element with
respect to the set-theoretic membership relation $ \in $ is
granted by the Axiom of Well-Founded Up Trees contained into the
axiomatic definition of a tree in  (the Mitchell-Benabou language
of) an \emph{elementary topos} $ {\mathcal{E}} $ stating that any
non-empty subtree of a tree has a maximal element.

\begin{definition}
\end{definition}
TREE ON $ {\mathcal{E}} $

an object T of  $ {\mathcal{E}} $ endowed with a binary relation $
R \, \mapsto \, T \times T $ such that:
\begin{description}
    \item[Poset:] R is a partial-order relation on T (and will be denoted 
from this reason as $
\leq $ from here and beyond
\item[Root:]
\begin{equation}
    \exists \,  0 \, \in T \; : \; ( 0 \, \leq \, t \; \forall t
\in T )
\end{equation}
  \item[Tree Property:]
\begin{equation}
   \downarrow t  \text{ is linearly ordered by $ \leq $
} \; \forall t \in T
\end{equation}
\item[Well-founded down:]
\begin{equation}
    \forall S \subseteq T \, : \, S \neq \emptyset \; , \exists y
\in S \, : \; y > z \; \Rightarrow \; z \notin S
\end{equation}
(i.e. y is minimal in S w.r.t. $ \leq $ )

\item[Well-founded up:]
\begin{equation}
    \forall S \subseteq T \; : \; ( S \neq \emptyset \; \Rightarrow
\;  \exists w \in S \, : \,  Not (  z > w \, : \, \forall z \in S)
\end{equation}
(i.e. w is maximal in S w.r.t. $ \leq $ )
\item[rigid]
\begin{equation}
    \alpha : T \mapsto T \; \text{ automorfism of T} \; \alpha \,
= \, {\mathbb{I}}
\end{equation}
where $ {\mathbb{I}} $ denotes the identity map on T while an
automorphism of T is defined as a $ leq$ preserving bijection of T

\end{description}
and where:
\begin{definition}
\end{definition}
DOWNWARD CLOSURE OF $ t \in T$:
\begin{equation}
    \downarrow t \; := \; \{ x \, : \, x \, \leq  t \}
\end{equation}

Informally and intuitivelly speaking, the key point consists in
representing a set x as a graph in the following way: x is the
root of the tree;on the first level are the members of x , joined
to x through a  directed edge. Next above there are the members of
the members y of x, each one joined to its y though a directed
edge.

Attempts to build "anti-foundational" axiomatizations of
Mathematics giving up Russell's Vicious Circle Principle (that
should more properly been called the Principle of Censorship of
Vicious Circles) can only arise whether the constraint that the
above graph doesn't contains cycles, and is conseguentially a
tree, is removed.

One possibility is given by the theory of hyper-sets, based on
formal system ZFA obtained by ZFC replacing the Axiom of
Foundation with a suitable Axiom of Anti-foundation.

Hypersets may be reached in many different ways;

Following the $ 2^{th} $ section "Ensembles non-bien fond\'{e}s "
of \cite{Longo-00}  and the chapter 2 "Background on set theory"
an the chapter 10 " Graphs"  of \cite{Barwise-Moss-96} I will
adhere to Aczel's formulations.

Let us assume a \emph{proper class} $ {\mathcal{U}} $ of
ur-elements that are not sets, are not classes and have no
members; given a set s:
\begin{definition}
\end{definition}
s IS TRANSITIVE:
\begin{equation}
    x \, \in  \, s \; \Rightarrow \; x \in {\mathcal{P}} ( s )
\end{equation}
\begin{definition}
\end{definition}
TRANSITIVE CLOSURE OF s:
\begin{equation}
    TC(s) \; := \; \bigcup \, \{ a \, , \, \bigcup a \, , \, \bigcup
\bigcup \, a \, , \, \cdots \}
\end{equation}
\begin{definition}
\end{definition}
SUPPORT OF s:
\begin{equation}
    support(s) \; := \; TC(s) \, \bigcap \, {\mathcal{U}}
\end{equation}
\begin{definition}
\end{definition}
s IS PURE:
\begin{equation}
    support(s) \; = \; \emptyset
\end{equation}
Given an $ A \, \subseteq \, {\mathcal{U}} $:
\begin{equation}
    V_{afa} (A) \; := \; \{ a \, : \, a \text{ is a set } \, and
\, support(a) \, \subseteq  \, A \}
\end{equation}
\begin{definition}
\end{definition}
CLASS OF THE PURE SETS:
\begin{equation}
    V_{afa} \; = \;  V_{afa} (\emptyset)
\end{equation}
Let us now define in a rigorous way the notion of \emph{graph}:
\begin{definition}
\end{definition}
GRAPH:

a couple $ {\mathbf{G}} \; := \; ( G \, , \, \rightarrow_{G} ) $
such that G is a set and $  \rightarrow_{G} $ is a binary relation
over G.

\smallskip

I will denote the proper class of all graphs by GRAPHS.

given $ {\mathbf{G}} = ( G \, , \, \rightarrow_{G} )  \in GRAPHS $
\begin{definition}
\end{definition}
d is  DECORATION OF $ {\mathbf{G}} \; \;  ( \, d = D ( {\mathbf{G}}) \,  ) $

if  d is a function $ d \,  : \,  {\mathbf{G}} \, \rightarrow \,
V_{afa} $ such that:
\begin{equation}
    d(a) \; = \; \{ d(b) \, : \, a \, \rightarrow_{G} b \} \; \;
\forall a , b \in {\mathbf{G}}
\end{equation}
One can than introduce the following:
\begin{axiom} \label{ax:axiom of anti-foundation}
\end{axiom}
AXIOM OF ANTI-FOUNDATION
\begin{equation}
     \forall \, g \,  \in \,GRAPHS  \; \exists ! \,  D(g)
\end{equation}
Let us now consider some formal system:
\begin{definition} \label{def:formal system ZFC-}
\end{definition}
FORMAL SYSTEM $ZFC^{-}$

the formal system with axioms:

\begin{description}
    \item[URELEMENTS]
\begin{equation*}
    (\forall p)(\forall q)[ {\mathcal{U}}(p) \; \rightarrow \;
\neg( q \in p)]
\end{equation*}
\item[EXTENSIONALITY]
\begin{equation*}
    (\forall a) ( \forall b ) [(\forall p)[p \in a \;
\leftrightarrow \; p \in b) \; \rightarrow \; a \, = \, b ]
\end{equation*}
  \item[PAIRING]
\begin{equation*}
  (\forall p)(\forall q) ( \exists a)[ p \in a \wedge q \in a ]
\end{equation*}
\item[UNION]
\begin{equation*}
    (\forall a) ( \exists b)(\forall c \in a )(\forall p \in c) p \in b
\end{equation*}
\item[POWER SET]
\begin{equation*}
   (\forall a) ( \exists b)(\forall c)[c \subseteq a \;
\rightarrow \; c \in b]
\end{equation*}
\item[INFINITY]
\begin{equation*}
    (\exists a)[\emptyset \in a \wedge (\forall b)[b \in a
\;\rightarrow \; (\exists c \in a)c=b \cup \{ b \}]]
\end{equation*}
\item[COLLECTION]
\begin{equation*}
   (\forall a)(\forall p \in a )(\exists q) \varphi (a,p,q) \;
\rightarrow \; (\exists b) (\forall p \in a)(\exists q \in b)
\varphi (a,p,q)
\end{equation*}
\item[SEPARATION]
\begin{equation*}
  (\forall a)(\exists b) (\forall p)[p \in b \; \leftrightarrow \;
p \in a \wedge \Phi (p,a)
\end{equation*}
\item[CHOICE]
\begin{equation*}
    (\forall a) ( \exists r)[ \text{r is a well-order of a}]
\end{equation*}
\item[STRONG PLENITUDE]
\begin{equation}
   (\forall a)(\forall b) [ {\mathcal{U}}(new(a,b)) \; \wedge \;
new (a,b) \notin b \; \wedge \; (\forall c \neq b) [new(a,b) \neq
new(c,b)]]
\end{equation}
\end{description}
\begin{definition} \label{def:formal system ZFC}
\end{definition}
FORMAL SYSTEM ZFC:
\begin{center}
  ZFC = $ ZFC^{-}$ \, + \, AXIOM \ref{ax:axiom of foundation}
\end{center}
\begin{definition} \label{def:formal system ZFA}
\end{definition}
FORMAL SYSTEM ZFA:
\begin{center}
  ZFA = $ ZFC^{-}$ \, + \, AXIOM \ref{ax:axiom of anti-foundation}
\end{center}

The replacement of the axiom\ref{ax:axiom of foundation}
axiom\ref{ax:axiom of anti-foundation} clearly affects the way by
which the comparison of sets is defined so that the Axiom of
Extensionality has to be looked at in a rather different way as it
was for ZFC. As to our considerations, it will be sufficient to
know that ZFA results to be a (consistent ?) formal system  about
which I demand to \cite{Barwise-Moss-96} getting rid of Russell's
paradox in a completelly non-orthodox way.

What it is important here to remark is that  ZFA is not , of
course, equiconsistent with the topos-theoretic foundation of
Mathematics (first order theory of a Well-pointed Topos) and
consequentially with RZC.

Whether Mazzola correctly claims that the passage to the
Topos-theoretical formalization of Mathematics may be seen, from a
philosophical perspective, as a kind of "behavioural revolution"
according to which the identity of a mathematical concept is
completely characterized by its properties as to inter-relation
with other mathematical concepts, his claim that a such a  passage
contains the  abrogation of the Vicious Circle Principle is ,
consequentially, completely false.

The claim that the so called "Yoneda Perspective" gives some kind
of conceptual consistence to postmodern intellectual hoaxes trying
to give an impression of  conceptual consistence to circularities
through semiological arguments or whatsoever (cfr. \cite{Eco-86}
and the $ 2^{th} $ chapter "Dizionario versus Enciclopedia" of
\cite{Eco-84}) is consequentially an intellectual hoax by itself.

One could, indeed, use the  equiconsistence of RZC and the
first-order theory of well-point topoi to obtain a topos-theoretic
version of ZFA. But such an anti-foundational formal system would
be completelly a different thing with respect to the Lawvere's
topos-theoretic foundation of Mathematics where, exactly as in the
set theoretic foundation, circularity is banned.

\bigskip

Since Mazzola's further constructions have as a ground the
formalization of the American school pioneered by Milton Babbit,
the necessary step consists in comparing the information-theoretic
language of \cite{Calude-02} I will adopt and some of the basic
notions of the american tradition as codified by Allen Forte in
\cite{Forte-73}.

Let us consider the scale  $ tempered (C_{2},12) $;

imponing octave-periodicity  it results that the musical alphabet
may be represented as the set of the residue classes modulo twelve
$ {\mathbb{Z}}_{12} $, endowed with its  algebraic structure of a
ring w.r.t. to the operations $ +_{12} $ and $ -_{12} $.

Let us introduce the following maps:
\begin{definition} \label{def:translation's operator on letters}
\end{definition}
TRANSLATION OF $ y \in  {\mathbb{Z}}_{12} $:

the map $ T_{y} \, : \,  {\mathbb{Z}}_{12} \, \mapsto \,
{\mathbb{Z}}_{12} $  such that:
\begin{equation}
    T_{y} (x) \; := \; x \, +_{12} \, y \; \; \forall x \in 
{\mathbb{Z}}_{12}
\end{equation}
and:
\begin{definition} \label{def:inversion's operator on letters}
\end{definition}
INVERSION OPERATOR:

the map $ Inv  \, : \, {\mathbb{Z}}_{12} \, \mapsto \,
{\mathbb{Z}}_{12} $ such that:
\begin{equation}
    Inv(x) \; := \; [0]_{12} \, -_{12} \, x
\end{equation}

It is curious to notice at this point how
a non-prime equal-tempering was coherent with the classical viewpoint on
esthetic consonance, summarized in axiom \ref{ax:axiom of the
naturality of esthetics}, giving rise to  simpler, i.e. non
co-prime ratios among the notes of physically consonant intervals;
from the viewpoint of Atonal Music a prime tempering would have
been, contrary, more reasonable: since 12 is not prime the ring $
{\mathbb{Z}}_{12} $ is not a (Galois) field and hence
multiplication $ \times_{12} $ and division $ \div_{12} $ among
\emph{pitch class sets} cannot be defined so that one cannot
implement the other two symmetry transformation that could,
instead, be implemented on $ {\mathbb{Z}}_{11} $ and $
{\mathbb{Z}}_{13} $.

Of course the musical alphabet $ {\mathbb{Z}}_{12} $ as well as
the basic symmetry transformation are easily implemented by the
following Mathematica expressions from section
\ref{sec:Mathematica's notebook Mathematical Music Toolkit}:
\begin{verbatim}
referencenote= 132;

(*** the  following instructions implements the musical alphabet
Z-12 ***)

letter[n_]:=Mod[n,12]

alphabet=Table[letter[n],{n,0,11}];

FROMletterTOnote[n_] :=
referencenote*Power[2,Power[letter[n],Power[12,-1]]]]

\end{verbatim}

I will adopt from here and beyond also the musical notation:
\begin{eqnarray}
  C \; &:=& \; [0]_{12} \\
  C^{\sharp} \; &:=& \; [1]_{12} \\
  D \; &:=& \; [2]_{12}  \\
  D^{\sharp} \; &:=& \;   [3]_{12}  \\
  E  \; & :=& \; [4]_{12} \\
  F \; & := & \; [5]_{12}    \\
  F^{\sharp} \; & :=  \;  [6]_{12}  \\
  G \; & := \; [7]_{12}   \\
  G^{\sharp} \; & := \;  [8]_{12}  \\
  A \; & := \; [9]_{12}  \\
  A^{\sharp}  \; & := \;  [10]_{12} \; \\
  B \; & := \; [11]_{12} \\
\end{eqnarray}

\medskip

Let us then consider the set $ {\mathbb{Z}}_{12}^{\star} $ of all
the finite strings (that I will also call words from here and
beyond) over $ {\mathbb{Z}}_{12} $, each element of which is
nothing but an \emph{ordered pitch class set} in Forte's
terminology.

Once more we can algorithmically implement the introduced notions
through the following Mathematica expression of the notebook
reported in section\ref{sec:Mathematica's notebook Mathematical
Music Toolkit}:
\begin{verbatim}
FROMletterTOnote[n_] :=
referencenote*Power[2,Power[letter[n],Power[12,-1]]]]

FROMwordTOscale[word_]:=
  Table[FROMletterTOnote[Part[word,i]],{i,1,Length[word]}]

(*** the following instructions introduced the notes without the
Mod-12  constraint ***)

c[1]=FROMletterTOnote[0];

c[n_]:=2*c[n-1]

c\[Sharp][1]=FROMletterTOnote[1];

c\[Sharp][n_]:=2*c\[Sharp][n-1]

d[1]=FROMletterTOnote[2];

d[n_]:=2*d[n-1]

d\[Sharp][1]=FROMletterTOnote[3];

d\[Sharp][n_]:=2*d\[Sharp][n-1]

e[1]=FROMletterTOnote[4];

e[n_]:=2*e[n-1]

f[1]=FROMletterTOFROMletterTOnote[5];

f[n_]:=2*f[n-1]

f\[Sharp][1]=FROMletterTOnote[6];

f\[Sharp][n_]:=2*f\[Sharp][n-1]

g[1]=FROMletterTOnote[7];

g[n_]:=2*g[n-1]

g\[Sharp][1]=FROMletterTOnote[8];

g\[Sharp][n_]:=2*g\[Sharp][n-1]

a[1]=FROMletterTOnote[9];

a[n_]:=2*a[n-1]

a\[Sharp][1]=FROMletterTOnote[10];

a\[Sharp][n_]:=2*a\[Sharp][n-1]

b[1]=FROMletterTOnote[11];

b[n_]:=2*b[n-1]

\end{verbatim}
where the recursive definition of $  x[n] \; x  \, = \, C , \cdots
, B \; n \, > \, 1 $ has being introduced in order to allow
concretelly  to play also higher octaves.

With this regard, if a a monodic piece is represented as a list
each element of which is itself a list of the form $ \{ \, note \,
, \,duration \}$ the involved code is:
\begin{verbatim}
(*** a monodic piece is a list each element of which is itself a list of the
    form {note ,duration}   ***)

playmonodic[piece_]:=
  Do[Play[Sin[piece[[i]][[1]] *2*\[Pi]*t],{t,0,piece[[i]][[2]]}] , {i,1,
      Length[piece]}]
\end{verbatim}

with the durations implemented as:
\begin{verbatim}
(*** insert the "referencetime" ***)

referencetime=4;

semibreve=1*referencetime;

minim= (1/2)*referencetime;

crotchet= (1/4)*referencetime;

quaver=(1/8)*referencetime;

semiquaver=(1/16)*referencetime;

demisemiquaver=(1/32)*referencetime;

hemidemisemiquaver=(1/64)*referencetime;
\end{verbatim}
and  with the following instructions allowing to associate (and to
eliminate) a fixed crothet duration:
\begin{verbatim}
FROMscaleTOpiece[scale_]:=Table[{Part[scale,i],crotchet},{i,1,Length[scale]}]

FROMpieceTOscale[piece_]:=Table[Part[piece,i,1],{i,1,Length[piece]}]
\end{verbatim}

\smallskip

Given a word $  \vec{x} \, \in \,   {\mathbb{Z}}_{12}^{\star} \, -
\, \bigcup_{k=1}^{4}  {\mathbb{Z}}_{12}^{k}  $ and an integer $ i
\, \in \{ 1, \, \cdots \, , | \vec{x} | \} $:
\begin{definition} \label{def:modes of a word}
\end{definition}
MODE OF $ \vec{x} $ of $ i^{th}$ DEGREE :
\begin{equation}
    mode( \vec{x} \, , \, i ) \; := \; {\mathcal{S}}_{cycl}^{i}
( \vec{x} )
\end{equation}
where $ {\mathcal{S}}_{cycl} $ denotes the operator of cyclic
shift.

Given furthermore an integer $ n \in { \mathbb{N}} $
\begin{definition} \label{def:chord of a word at a certain level}
\end{definition}
CHORD OF  $ \vec{x} $ of $ i^{th}$ DEGREE AT LEVEL n :
\begin{equation}
    chord ( \vec{x} \, , \, i \, , \, n ) \; := \cdot_{j=0}^{n+3} mode( 
\vec{x} \, , \, i )_{2j+1}
\end{equation}
where $ \cdot $ denotes the concatenation operator while $  x_{j}
$ denotes the $ j^{th} $ letter of the word $ \vec{x} $.

Let us observe,at this point, that exists a certain maximum level
$ maxlevel( \vec{x} ) $ such that $ chord ( \vec{x} \, , \, i \, , \, n
) \, $ for $ n > maxlevel( \vec{x} ) $ simply adds notes already
contained in the chord.

A simple reasoning allows to infer that:
\begin{equation} \label{eq:maximum level f a word for chords}
    maxlevel ( \vec{x} ) \; = \; 1 \, + \, Int(\frac{ |\vec{x} | -5}{2}) \, 
-
\, \frac{ (-1)^{ |\vec{x} |} - 1 }{2} \dot Int( \frac{ |\vec{x}
|}{2})
\end{equation}

\bigskip

Given our alphabet $ \Sigma := {\mathbb{Z}}_{12} $, a word $
\vec{x} \in \Sigma^{\star} $ and a map $ g \, : \, \Sigma
\rightarrow \Sigma $:
\begin{definition} \label{def:map induced on words by a map on an alphabet}
\end{definition}
MAP INDUCED BY g ON WORDS:

the map $ \hat{g} \, : \, \Sigma^{\star} \, \rightarrow \,
\Sigma^{\star} $ :
\begin{equation}
    \hat{g} ( \vec{x} ) \; = \; \cdot_{i=1}^{| \vec{x} |} \, g (
    x_{i})
\end{equation}

Applying , in particular,  the definition\ref{def:map induced on
words by a map on an alphabet} to the maps of definition
\ref{def:translation's operator on letters} and definition
\ref{def:inversion's operator on letters} one obtains the
following
\begin{definition} \label{def:translation's operator on words}
\end{definition}
TRANSLATION OF $ y \in  {\mathbb{Z}}_{12} $ ON WORDS:

the map $ \hat{T}_{y} \, : \,  {\mathbb{Z}}_{12}^{\star} \,
\mapsto \, {\mathbb{Z}}_{12}^{\star} $

\begin{definition} \label{def:inversion's operator on words}
\end{definition}
INVERSION OPERATOR ON WORDS;

the map $ \hat{Inv}  \, : \, {\mathbb{Z}}_{12}^{\star} \, \mapsto
\, {\mathbb{Z}}_{12}^{\star} $

\bigskip

These operators naturally induce  the following equivalence
relations:
\begin{definition} \label{def:translation equivalence of words}
\end{definition}
$ \vec{x}_{1} \, , \, \vec{x}_{2} \, \in \,
{\mathbb{Z}}_{12}^{\star} $ ARE TRANSLATIONALLY-EQUIVALENT ( $
\vec{x}_{1} \, \sim_{T} \, \vec{x}_{2} $):
\begin{equation}
   \exists y \in {\mathbb{Z}}_{12} \; : \; \vec{x}_{2} \, = T_{y}  
\vec{x}_{1}
\end{equation}
\begin{definition} \label{def:inversion equivalence of words}
\end{definition}
$ \vec{x}_{1} \, , \, \vec{x}_{2} \, \in \,
{\mathbb{Z}}_{12}^{\star} $ ARE INVERSIONALLY-EQUIVALENT ( $
\vec{x}_{1} \, \sim_{Inv} \, \vec{x}_{2} $):
\begin{equation}
    \vec{x}_{2} \, = Inv  \vec{x}_{1}
\end{equation}

that may be managed through the following instructions from
section\ref{sec:Mathematica's notebook Mathematical Music
Toolkit}:

\begin{verbatim}
(*** the following instructions implement the two operation of translation
    and inversion at the basis of atonal music as well as the relative
    tests ***)

translation[word_,n_]:=Table[Mod[Part[word,i]+n,12],{i,1,Length[word]}]

translationequivalenceofwordsQ[w1_,w2_]:=
  Not[Equal[Table[Equal[w2,translation[w1,n]],{n,0,11}],Table[False,{12}]]]

inversion[word_]:=Table[Mod[12-Part[word,i],12],{i,1,Length[word]}]

inversionequivalenceofwordsQ[w1_,w2_]:=Equal[w1,inversion[w2]]

inversioninvarianceofawordQ[w_]:=Equal[w,inversion[w]]
\end{verbatim}

Let us now analyze the musical intervals inside a word: instead of
the complicated Forte's approach I will adopt a  more intuitive
definition of the \emph{interval vector}:

given two musical letters $ x , y \in {\mathbb{Z}}_{12} $:
\begin{definition}
\end{definition}
DISTANCE BETWEEN x AND y :
\begin{equation}
    d ( x \, , \, y) \; := \; x \, -_{12} \, y
\end{equation}
Given a vord $ \vec{x} \in {\mathbb{Z}}_{12}^{ \star } $:
\begin{definition}
\end{definition}
INTERVAL VECTOR OF $ \vec{x} $:
\begin{equation}
   I( \vec{x} ) \; := \; ( d( x_{2} , x_{1} ) \, , \, \cdots \, , d( x_{| 
\vec{x} |} , x_{| \vec{x} | - 1 }
))
\end{equation}
One has obviously that:
\begin{theorem} \label{th:the interval vector determines a translational 
equivalence class}
\end{theorem}
THE INTERVAL VECTOR DETERMINES A TRANSLATIONAL EQUIVALENCE'S CLASS
\begin{equation}
    \vec{x} \, \sim_{T} \, \vec{y} \; \Leftrightarrow \; I(
    \vec{x} ) \, = \, I(
    \vec{y} )
\end{equation}
\begin{proof}
\begin{description}
\item[Proof of the implication $ \Rightarrow $]

By hypothesis:
\begin{equation}
    \exists \, z \, \in {\mathbb{Z}}_{12} \; : \; \vec{y} \, =
\, T_{z} \vec{x}
\end{equation}
One has that:
\begin{multline}
    I( \vec{y}) \;  \; = \;  ( d( x_{2} +_{12} z , x_{1}+_{12} z  ) \, , \, 
\cdots \,
, d( x_{| \vec{x} |} +_{12} z , x_{| \vec{x} | - 1 } +_{12} z ) )
\; = \\
= \, ( d( x_{2} , x_{1} ) \, , \, \cdots \, , d( x_{| \vec{x} |} ,
x_{| \vec{x} | - 1 } ) \; = \; I( \vec{x})
\end{multline}
\item[Proof of the implication $ \Leftarrow $]

By hypothesis:
\begin{equation}
    ( d( x_{2} , x_{1} ) \, , \, \cdots \, , d( x_{| \vec{x}
|} , x_{| \vec{x} | - 1 } ) \; = \; ( d( y_{2} , y_{1} ) \, , \,
\cdots \, , d( y_{| \vec{x} |} , y_{| \vec{y} | - 1 } )
\end{equation}
This implies that:
\begin{multline}
    \exists \, z \, \in \, {\mathbb{Z}}_{12} \; : \; ( d( y_{2} , y_{1} ) \, 
, \,
\cdots \, , d( y_{| \vec{x} |} , y_{| \vec{y} | - 1 } ) \, = \\
= \,  ( d( x_{2} +_{12} z , x_{1}+_{12} z  ) \, , \, \cdots \, ,
d( x_{| \vec{x} |} +_{12} z , x_{| \vec{x} | - 1 } +_{12} z ) )
\end{multline}
implying that:
\begin{equation}
     \vec{y} \, = \, T_{z} \vec{x}
\end{equation}
\end{description}
\end{proof}

Let us now fix a useful notation: words will be represented
enclosed by  $\{ \, \cdots \, \} $ while interval vectors,and
hence translationally-equivalences' classes, will be represented
enclosed by $ (   \, \cdots \, ) $.

So, for example, we will speak about  the diatonic major scale $ (
2 , 2 , 1 , 2 , 2 , 2 , 1 ) $ referring to the translationally
equivalence class:
\begin{equation}
    (2 , 2 , 1 , 2 , 2 , 2 , 1 ) \; = [ {0,2,4,5,7,9,11} ]_{T}
\end{equation}

\medskip
Let us now introduce a particularly important subset of musical
words. Let us consider, at this purpose, a generic finite (i.e. $
card( \Sigma ) < \infty $ ) alphabet $ \Sigma $. We have clearly
that:
\begin{eqnarray}
  card ( \Sigma^{n} )  \; & =& \;  (card(\Sigma))^{n}  \\
  card (\Sigma^{\star}) \; &=&  \; \aleph_{0}  \\
  card( \Sigma^{\infty} ) \; &=& \; \aleph_{1}
\end{eqnarray}
Let us now introduce the following:
\begin{definition} \label{def:nonrepetitivewords of length n}
\end{definition}
NONREPETITIVE WORDS OVER $  \Sigma $ OF LENGTH $ n \in {
\mathbb{N}} $
\begin{equation}
   \Sigma^{n}_{NR} \; := \; \{ \vec{x} \in  \Sigma^{n} \, : \, x_{i} \neq 
x_{j} \; \forall i \neq j \}
\end{equation}
\begin{definition} \label{def:nonrepetitivewords}
\end{definition}
NONREPETITIVE WORDS OVER $  \Sigma $:
\begin{equation}
   \Sigma^{\star}_{NR} \; := \; \{ \vec{x} \in  \Sigma^{\star} \, : \, x_{i} 
\neq x_{j} \; \forall i \neq j \}
\end{equation}
Trivial combinatorial considerations allow to prove the following:
\begin{theorem} \label{th:cardinalities of the sets of nonrepetitive words}
\end{theorem}
CARDINALITIES OF THE SETS OF  NONREPETITIVE WORDS:
\begin{eqnarray}
  card( \Sigma^{n}_{NR} ) \; &=& \; \binom{card(\Sigma)}{n}  \\
  card( \Sigma^{\star}_{NR} )\;  &=& \; \sum_{k=1}^{card(\Sigma)} 
\binom{card(\Sigma)}{k}
\end{eqnarray}
Applying theorem\ref{th:cardinalities of the sets of nonrepetitive
words} to our particular musical alphabet $ {\mathbb{Z}}_{12}$
one obtains the following:
\begin{corollary} \label{th:cardinalities of the sets of nonrepetitive 
musical words}
\end{corollary}
CARDINALITIES OF THE SET OF NONREPETITIVE MUSICAL WORDS
\begin{eqnarray}
  card( {\mathbb{Z}}_{12 \; NR}^{n} ) \; &=& \; \binom{12}{n}  \\
  card( {\mathbb{Z}}_{12 \; NR}^{\star} )\;  &=& \; 
\sum_{k=1}^{card(\Sigma)}
\binom{12}{n} \; = \; 4095
\end{eqnarray}

Mazzola's strong generalization of the notion of tonality is based
on the  formalization of  the 3-voices harmonization of the
degrees of arbitrary nonrepetitive words of length equal to seven
$ \Sigma^{7}_{NR} $.

Gregorian Music \footnote{We will adopt here the standard
denomination of gregorian modes demanding to \cite{Turco-96a} ,
\cite{Turco-96b} and to the fourth chapter "Le scale" of
\cite{Tammaro-03} for any further information concerning their
functionalities in the catholic liturgy and their hystorical
derivation from greek modes} was based on the following strongly
smaller subset of $  \Sigma^{7}_{NR}  $:
\begin{definition} \label{def:gregorian-words}
\end{definition}
GREGORIAN-WORDS
\begin{equation}
   {\mathcal{S}}_{greg} \; := \; \{ \, mode(\vec{x},i) \, \vec{x} \in
\Sigma^{7}_{NR} \; : \;
  I(\vec{x}) = ( 2 , 2 , 1 , 2 , 2 , 2 , 1 ) \, , \, i = 1 \, \, \cdots \, 7 
\}
\end{equation}
One has clearly that:
\begin{equation}
    card (  {\mathcal{S}}_{greg} ) \; = \;  84
\end{equation}
since $  {\mathcal{S}}_{greg} $ contains the 12 version starting
from C to B, of the following modes:
\begin{equation}
    I(s) \;  = \; \left\{%
\begin{array}{ll}
    ( 2 , 2 , 1 , 2 , 2 , 2 , 1 ), & \hbox{ionian (major)} \\
    (2,1,2,2,2,1,2), & \hbox{dorian} \\
    (2,1,2,2,2,1,2), & \hbox{phrigian} \\
    (2,2,2,1,2,2,1), & \hbox{lydian} \\
    (2,2,1,2,2,1,2), & \hbox{myxolydian} \\
    (2,1,2,2,1,2,2), & \hbox{aeolian (minor)} \\
    (1,2,2,1,2,2,2), & \hbox{locrian} \\
\end{array}%
\right.
\end{equation}
i.e. all the elements of $ \Sigma^{7}_{NR} $ whose labels are made
of 2 with the exception of two 1's at distance 4.

Classical Harmony, furthermore, is based on the following smaller
subset of $ {\mathcal{S}}_{greg}  $:
\begin{definition} \label{def:classical-words}
\end{definition}
CLASSICAL-WORDS
\begin{equation}
     {\mathcal{S}}_{cl} \; := \; \{ [i]_{12}-ionian \, , \,
     [i]_{12}-aeolian  \; , \; i \, = \, 1 \, , \,  \cdots  \, 12 \}
\end{equation}
One has clearly that:
\begin{equation}
    card (  {\mathcal{S}}_{cl} ) \; = \;  24
\end{equation}

Among the overwhelming majority of non-gregorian (and hence also
non-classical) nonrepetitive words worth of note are the harmonic
minor scale (2,1,2,2,1,3,1) and its derived modes ( locrian $
\sharp 6 $ $(1,2,2,1,3,1,2) $ , augmented ionic $ (2,2,1,3,1,2,1)
$ , minor lydian $ 7 \flat $ $ (2,1,3,1,2,1,2) $, myxolidian $
\flat 2 \diagup 6 \flat $ $(1,3,1,2,1,2,2)$ , lydian $ \sharp 2 $
$(3,1,2,1,2,2,1)$  , super-locrian 7 diminuished
$(1,2,1,2,2,1,3)$), the melodic minor scale $ (2,1,2,2,2,2,1) $
and its derived modes (dorian $ \flat 2 $ $ (1,2,2,2,2,1,2) $,
augmented lydian $ (2,2,2,2,1,2,1)$, dominating lydian $
(2,2,2,1,2,1,2) $, myxolidian $ \flat 6$ $(2,2,1,2,1,2,2)$,
locrian $ \sharp 2$ $(2,1,2,1,2,2,2)$, superlocrian
$(1,2,1,2,2,2,2)$)  as well as  the 12-equally tempered
approcimation of many ethnic scales cfr. the voices "Modi
(Modalit\'{a}"  and "Scale musicali antiche e moderne" of
\cite{Righini-80} and \cite{Castro-98} \footnote{It should be
observed that it is only  by 12-equally tempered approximation of
the makam "Ahavah Rabbah"  that it collapes to the fifth's mode of
the harmonic minor scale. These brief and very schematic
observations could be considered, anyway, as no more than a taste
about the complexity of all the studies concerning Jewish Music,
its identity, its esthetics (and consequentially, according to
those insisting on the link music versus poietry , cfr. the
subsection "Musica e Poesia" of the $ 2^{th} $ chapter
"L'Occidente cristiano e l'idea di Musica" of \cite{Fubini-95}),
its poetics \cite{Schoekel-89}), the bi-directional flow of
information with Christian music in its hystorical develop
\cite{Levi-02}. I demand to \cite{Fubini-94},
\cite{Bayer-Avenari-Boehm-72}, \cite{Adler-03} and
\cite{Seroussi-Braun-Schleifer-Sharvit-Manasseh-Kaufman-Shelemay-Hirshberg-Bohlman-Werb-Feldman-Harran-Knapp-Bloch-Thwaite-Levin-Yating-01}
for any further information} such as the major tzigan
$(1,3,2,1,1,3,1) $ and minor tzigan scale $ (2,1,3,1,1,3,1) $, the
indian scale $ (1,3,1,2,1,3,1) $ , the ind\'{u} scale
$(2,2,1,2,1,2,2)$, the hungarian scale $ ( 3,1,2,1,2,1,2) $  and
the minor napoletan scale $ (1,2,2,2,1,3,1) $.

We are now ready to introduce the following notions:

\begin{definition} \label{def:Mazzola tonality}
\end{definition}
MAZZOLA TONALITY:

a couple  $  ( \vec{x} \, ,  \, chord^{(3)} ) $ where:
\begin{enumerate}
    \item
\begin{equation}
    \vec{x} \; \in \; \Sigma^{7}_{NR}
\end{equation}
    \item $ chord^{(3)} \, : \, \{1,2,3,4,5,6,7 \} \, \mapsto  \, 
{\mathbb{Z}}_{12}^{3} $ is the map associating to
    to the $ i^{th} $-degree the chord:
\begin{equation}
    chord^{(3)} (i)  \; = \; chord ( \vec{x} \, , \, i \, , \, 1 )
\end{equation}
\end{enumerate}
Since a Mazzola-tonality is determined by its underlying
non-repetitive eptatonic words, the set $ {\mathcal{T}}_{Maz} $ of
all the possible Mazzola-tonalities is bijective to $
\Sigma^{7}_{NR} $ , so that, clearly, one have that:
\begin{equation}
    card ( {\mathcal{T}}_{Maz} ) \; = \; 792
\end{equation}
Similarly, introduced the following notions:
\begin{definition} \label{def:gregorian tonalities}
\end{definition}
GREGORIAN TONALITIES:
\begin{equation}
    {\mathcal{T}}_{greg} \; := \; \{ ( \vec{x} \, , \, chord^{(3)} ) \, :
    \,  \vec{x} \in {\mathcal{S}}_{greg} \}
\end{equation}
\begin{definition} \label{def:classical tonalities}
\end{definition}
CLASSICAL TONALITIES:
\begin{equation}
    {\mathcal{T}}_{cl} \; := \; \{ ( \vec{x} \, , \, chord^{(3)} ) \, :
    \,  \vec{x} \in {\mathcal{S}}_{cl} \}
\end{equation}
the bijectivity of $ {\mathcal{S}}_{greg} $ and  $
{\mathcal{T}}_{greg} $ and the bijectivity of  $
{\mathcal{S}}_{cl} $ and  $ {\mathcal{T}}_{cl} $ imply that:
\begin{equation}
    card ( {\mathcal{T}}_{greg} ) \; = \; 84
\end{equation}
\begin{equation}
    card ( {\mathcal{T}}_{cl} ) \; = \; 24
\end{equation}
Let us now observe that, as to the possibility of performing the
harmonization of the degrees on the scale, there is no reason to
restrict ourselves to $ \Sigma^{7}_{NR} $; let us observe, with
this regard  \cite{Righini-80}, \cite{Castro-98}, that the popular
music of many cultures adopts commonly  the pentatonic scale
$(2,2,3,2,4)$ and its modes, that the 12-equally tempered
approximation of some ethnic scale such as the japanese Hira Joski
$ (1,4,1,4,2)$ or the japanese pelog $ ( 1,4,1,4,2 ) $ are also
pentatonic, that the blues scale $ ( 3,2,1,1,3 ) $ is made of six
notes and the many be-bop scales such as the bebop major $
(2,2,1,2,1,1,2,1) $ or the be-bop dominant $ (2,2,1,2,2,1,1,1) $
and the relative mode are made of 8 notes and so on.

Looking at  Jazz Harmony \cite{Levine-95}, \cite{Strunk-00},
\cite{Gramaglia-92}, \cite{Mariani-02} furthermore, it appears
natural to give up the restriction of considering only first-level
chords.

One arrives, consequentially, to the following definition: given a
word $ \vec{x} \,  \in \,  \Sigma^{\star}_{NR} - \bigcup_{k=1}^{4}
\Sigma^{k}_{NR} $ and a number $ n \in {\mathbb{N}}_{+} $ such
that $ n \; < \;  maxlevel( \vec{x} ) $
\begin{definition} \label{def:tonality of a non-repetitive word at a certain 
level}
\end{definition}
TONALITY OF THE WORD $ \vec{x} \, \in \, \Sigma^{\star}_{NR} \, -
\, \bigcup_{k=1}^{4} \Sigma^{k}_{NR} $ AT THE LEVEL n (n-TONALITY
OF $ \vec{x} $) :

$ tonality( \vec{x} \, , \, n) \; := \;    ( \vec{x} \, , \,
chord^{(n)} ) $, where $ chord^{(n)} \, : \, \{ 1  , \cdots ,
|\vec{x} | \} \, \mapsto \, {\mathbb{Z}}_{12}^{\star} $ is the map
such that:
\begin{equation}
    chord^{(n)} (i)  \; = \; chord ( \vec{x} \, , \, i \, , \, n )
\end{equation}

I will denote the set of all the n-tonalities by
${\mathcal{T}}_{n} $ while I will denote by $ {\mathcal{T} } $ the
set of all the tonalities at any level. One has clearly that:
\begin{multline}
   card ( {\mathcal{T} } ) \; = \; \sum_{\vec{x} \, \in \, 
\Sigma^{\star}_{NR} \, -
\, \bigcup_{k=1}^{4} \Sigma^{k}_{NR}} \,
\sum_{n=1}^{maxlevel(\vec{x})} 1   \;
= \; \sum_{k=5}^{12} card ( \Sigma_{NR}^{k} ) \, \sum_{n=1}^{1 \, + \, 
Int(\frac{ k -5}{2}) \, - \, \frac{ (-1)^{
k} - 1 }{2} \dot Int( \frac{k}{2})} 1 \\ = \; 10100
\end{multline}

Given an n-tonality $ t \in {\mathcal{T}}_{n} $:
\begin{definition}
\end{definition}
HARMONIC WORDS OF t:
\begin{equation}
    {\mathcal{HW}} (t) \; := \; (Range( chord^{(n)} (t)))^{\star}
\end{equation}

Given two tonalities $ t_{1} , t_{2} \; \in \; {\mathcal{T}} $:
\begin{definition}
\end{definition}
$ t_{1} \; AND \; t_{2} $ ARE TRANSLATIONALLY-EQUIVALENT ( $ t_{1}
\, \sim_{T} \, t_{2} )$:
\begin{equation}
    \exists z \in {\mathbb{Z}}_{12} \; : \; chord( t_{1} \; i ) \, = \, 
T_{z}  chord( t_{2} \; i
) \; \; \forall i
\end{equation}

Once again all these mathematical concepts are easily implemented
trough the following expressions of the section
\ref{sec:Mathematica's notebook Mathematical Music Toolkit}:
\begin{verbatim}
(*** mode[word, i])  gives the i - th mode of a word
on the musical alphabet   Z-12 ***)

mode[word_,i_]:=RotateLeft[word,i-1];

(***  chord[word,i,level] gives the chord of the i-th degree of a
word at a chosen level ***)

chord[word_,i_,level_]:= Table[Part[mode[word,i],
      
If[Mod[2n+1,Length[word]]==0,Length[word],Mod[2n+1,Length[word]]]],{n,0,
      1+level}]

tonality[word_,level_]:=Table[chord[word,i,level],{i,1,Length[word]}]

harmonicwords[t_,n_]:=Strings[t,n]

harmonicwordsupto[t_,n_]:=
  If[n==1,harmonicwords[t,1],
    Join[harmonicwordsupto[t,n-1],harmonicwords[t,n]]]

harmonicwordintermofdegrees[t_,listofdegrees_]:=
  Table[Part[t,Part[listofdegrees,i]],{i,1,Length[listofdegrees]}]

degreeofachordinatonality[w_,t_]:=Part[Flatten[Position[t,w]],1]

FROMharmonicwordTOphysicalchord[hw_]:=
  Table[FROMwordTOscale[Part[hw,i]],{i,1,Length[hw]}]

tonalitymembershipQ[hw_,t_]:=MemberQ[harmonicwords[t,Length[hw]],hw]
\end{verbatim}

It may be useful to report some example:

\begin{example}
\end{example}
TONALITIES OF THE MOST COMMON WORDS AT DIFFERENT LEVELS

\begin{verbatim}
tonality[majorword[0],1] =
{{0,4,7},{2,5,9},{4,7,11},{5,9,0},{7,11,2},{9,0,4},{11,2,5}}

tonality[majorword[0],2] =
{{0,4,7,11},{2,5,9,0},{4,7,11,2},{5,9,0,4},{7,11,2,5},{9,0,4,7},{11,2,5,9}}

tonality[majorword[0],3] =
{{0,4,7,11,2},{2,5,9,0,4},{4,7,11,2,5},{5,9,0,4,7},{7,11,2,5,9},{9,0,4,7,11},{
    11,2,5,9,0}}

tonality[majorword[0],4] =
{{0,4,7,11,2,5},{2,5,9,0,4,7},{4,7,11,2,5,9},{5,9,0,4,7,11},{7,11,2,5,9,0},{9,
    0,4,7,11,2},{11,2,5,9,0,4}}

tonality[majorword[0],5] =
{{0,4,7,11,2,5,9},{2,5,9,0,4,7,11},{4,7,11,2,5,9,0},{5,9,0,4,7,11,2},{7,11,2,
    5,9,0,4},{9,0,4,7,11,2,5},{11,2,5,9,0,4,7}}

tonality[minorword[0],1] =
{{0,3,7},{2,5,8},{3,7,10},{5,8,0},{7,10,2},{8,0,3},{10,2,5}}

tonality[minorword[0],2] =
{{0,3,7,10},{2,5,8,0},{3,7,10,2},{5,8,0,3},{7,10,2,5},{8,0,3,7},{10,2,5,8}}

tonality[minorword[0],3] =
{{0,3,7,10,2},{2,5,8,0,3},{3,7,10,2,5},{5,8,0,3,7},{7,10,2,5,8},{8,0,3,7,10},{
    10,2,5,8,0}}

tonality[minorword[0],4] =
{{0,3,7,10,2,5},{2,5,8,0,3,7},{3,7,10,2,5,8},{5,8,0,3,7,10},{7,10,2,5,8,0},{8,
    0,3,7,10,2},{10,2,5,8,0,3}}

tonality[minorword[0],5] =
{{0,3,7,10,2,5,8},{2,5,8,0,3,7,10},{3,7,10,2,5,8,0},{5,8,0,3,7,10,2},{7,10,2,
    5,8,0,3},{8,0,3,7,10,2,5},{10,2,5,8,0,3,7}}

tonality[harmonicminorword[0],1] =
{{0,3,7},{2,5,8},{3,7,11},{5,8,0},{7,11,2},{8,0,3},{11,2,5}}

tonality[harmonicminorword[0],2] =
{{0,3,7,11},{2,5,8,0},{3,7,11,2},{5,8,0,3},{7,11,2,5},{8,0,3,7},{11,2,5,8}}

tonality[harmonicminorword[0],3] =
{{0,3,7,11,2},{2,5,8,0,3},{3,7,11,2,5},{5,8,0,3,7},{7,11,2,5,8},{8,0,3,7,11},{
    11,2,5,8,0}}

tonality[harmonicminorword[0],4] =
{{0,3,7,11,2,5},{2,5,8,0,3,7},{3,7,11,2,5,8},{5,8,0,3,7,11},{7,11,2,5,8,0},{8,
    0,3,7,11,2},{11,2,5,8,0,3}}

tonality[harmonicminorword[0],5] =
{{0,3,7,11,2,5,8},{2,5,8,0,3,7,11},{3,7,11,2,5,8,0},{5,8,0,3,7,11,2},{7,11,2,
    5,8,0,3},{8,0,3,7,11,2,5},{11,2,5,8,0,3,7}}

tonality[dorianword[0],1] =
{{0,3,7},{2,5,9},{3,7,10},{5,9,0},{7,10,2},{9,0,3},{10,2,5}}

tonality[dorianword[0],2] =
{{0,3,7,10},{2,5,9,0},{3,7,10,2},{5,9,0,3},{7,10,2,5},{9,0,3,7},{10,2,5,9}}

tonality[dorianword[0],3] =
{{0,3,7,10,2},{2,5,9,0,3},{3,7,10,2,5},{5,9,0,3,7},{7,10,2,5,9},{9,0,3,7,10},{
    10,2,5,9,0}}

tonality[dorianword[0],4] =
{{0,3,7,10,2,5},{2,5,9,0,3,7},{3,7,10,2,5,9},{5,9,0,3,7,10},{7,10,2,5,9,0},{9,
    0,3,7,10,2},{10,2,5,9,0,3}}

tonality[dorianword[0],5] =
{{0,3,7,10,2,5,9},{2,5,9,0,3,7,10},{3,7,10,2,5,9,0},{5,9,0,3,7,10,2},{7,10,2,
    5,9,0,3},{9,0,3,7,10,2,5},{10,2,5,9,0,3,7}}

tonality[phrigianword[0],1] =
{{0,3,7},{1,5,8},{3,7,10},{5,8,0},{7,10,1},{8,0,3},{10,1,5}}

tonality[phrigianword[0],2] =
{{0,3,7,10},{1,5,8,0},{3,7,10,1},{5,8,0,3},{7,10,1,5},{8,0,3,7},{10,1,5,8}}

tonality[phrigianword[0],3] =
{{0,3,7,10,1},{1,5,8,0,3},{3,7,10,1,5},{5,8,0,3,7},{7,10,1,5,8},{8,0,3,7,10},{
    10,1,5,8,0}}

tonality[phrigianword[0],4] =
{{0,3,7,10,1,5},{1,5,8,0,3,7},{3,7,10,1,5,8},{5,8,0,3,7,10},{7,10,1,5,8,0},{8,
    0,3,7,10,1},{10,1,5,8,0,3}}

tonality[phrigianword[0],5] =
{{0,3,7,10,1,5,8},{1,5,8,0,3,7,10},{3,7,10,1,5,8,0},{5,8,0,3,7,10,1},{7,10,1,
    5,8,0,3},{8,0,3,7,10,1,5},{10,1,5,8,0,3,7}}

tonality[lydianword[0],1] =
{{0,4,7},{2,6,9},{4,7,11},{6,9,0},{7,11,2},{9,0,4},{11,2,6}}

tonality[lydianword[0],2] =
{{0,4,7,11},{2,6,9,0},{4,7,11,2},{6,9,0,4},{7,11,2,6},{9,0,4,7},{11,2,6,9}}

tonality[lydianword[0],3] =
{{0,4,7,11,2},{2,6,9,0,4},{4,7,11,2,6},{6,9,0,4,7},{7,11,2,6,9},{9,0,4,7,11},{
    11,2,6,9,0}}

tonality[lydianword[0],4] =
{{0,4,7,11,2,6},{2,6,9,0,4,7},{4,7,11,2,6,9},{6,9,0,4,7,11},{7,11,2,6,9,0},{9,
    0,4,7,11,2},{11,2,6,9,0,4}}

tonality[lydianword[0],5] =
{{0,4,7,11,2,6,9},{2,6,9,0,4,7,11},{4,7,11,2,6,9,0},{6,9,0,4,7,11,2},{7,11,2,
    6,9,0,4},{9,0,4,7,11,2,6},{11,2,6,9,0,4,7}}

tonality[mixolydianword[0],1] =
{{0,4,7},{2,5,9},{4,7,10},{5,9,0},{7,10,2},{9,0,4},{10,2,5}}

tonality[mixolydianword[0],2] =
{{0,4,7,10},{2,5,9,0},{4,7,10,2},{5,9,0,4},{7,10,2,5},{9,0,4,7},{10,2,5,9}}

tonality[mixolydianword[0],3] =
{{0,4,7,10,2},{2,5,9,0,4},{4,7,10,2,5},{5,9,0,4,7},{7,10,2,5,9},{9,0,4,7,10},{
    10,2,5,9,0}}

tonality[mixolydianword[0],4] =
{{0,4,7,10,2,5},{2,5,9,0,4,7},{4,7,10,2,5,9},{5,9,0,4,7,10},{7,10,2,5,9,0},{9,
    0,4,7,10,2},{10,2,5,9,0,4}}

tonality[mixolydianword[0],5] =
{{0,4,7,10,2,5,9},{2,5,9,0,4,7,10},{4,7,10,2,5,9,0},{5,9,0,4,7,10,2},{7,10,2,
    5,9,0,4},{9,0,4,7,10,2,5},{10,2,5,9,0,4,7}}

tonality[locrianword[0],1] =
{{0,3,6},{1,5,8},{3,6,10},{5,8,0},{6,10,1},{8,0,3},{10,1,5}}

tonality[locrianword[0],2] =
{{0,3,6,10},{1,5,8,0},{3,6,10,1},{5,8,0,3},{6,10,1,5},{8,0,3,6},{10,1,5,8}}

tonality[locrianword[0],3] =
{{0,3,6,10,1},{1,5,8,0,3},{3,6,10,1,5},{5,8,0,3,6},{6,10,1,5,8},{8,0,3,6,10},{
    10,1,5,8,0}}

tonality[locrianword[0],4] =
{{0,3,6,10,1,5},{1,5,8,0,3,6},{3,6,10,1,5,8},{5,8,0,3,6,10},{6,10,1,5,8,0},{8,
    0,3,6,10,1},{10,1,5,8,0,3}}

tonality[locrianword[0],5] =
{{0,3,6,10,1,5,8},{1,5,8,0,3,6,10},{3,6,10,1,5,8,0},{5,8,0,3,6,10,1},{6,10,1,
    5,8,0,3},{8,0,3,6,10,1,5},{10,1,5,8,0,3,6}}

tonality[tziganminorword[0],1] =
{{0,3,7},{2,6,8},{3,7,11},{6,8,0},{7,11,2},{8,0,3},{11,2,6}}

tonality[tziganminorword[0],2] =
{{0,3,7,11},{2,6,8,0},{3,7,11,2},{6,8,0,3},{7,11,2,6},{8,0,3,7},{11,2,6,8}}

tonality[tziganminorword[0],3] =
{{0,3,7,11,2},{2,6,8,0,3},{3,7,11,2,6},{6,8,0,3,7},{7,11,2,6,8},{8,0,3,7,11},{
    11,2,6,8,0}}

tonality[tziganminorword[0],4] =
{{0,3,7,11,2,6},{2,6,8,0,3,7},{3,7,11,2,6,8},{6,8,0,3,7,11},{7,11,2,6,8,0},{8,
    0,3,7,11,2},{11,2,6,8,0,3}}

tonality[tziganminorword[0],5] =
{{0,3,7,11,2,6,8},{2,6,8,0,3,7,11},{3,7,11,2,6,8,0},{6,8,0,3,7,11,2},{7,11,2,
    6,8,0,3},{8,0,3,7,11,2,6},{11,2,6,8,0,3,7}}

tonality[jewishword[0],1] =
{{0,4,7},{1,5,8},{4,7,10},{5,8,0},{7,10,1},{8,0,4},{10,1,5}}

tonality[jewishword[0],2] =
{{0,4,7,10},{1,5,8,0},{4,7,10,1},{5,8,0,4},{7,10,1,5},{8,0,4,7},{10,1,5,8}}

tonality[jewishword[0],3] =
{{0,4,7,10,1},{1,5,8,0,4},{4,7,10,1,5},{5,8,0,4,7},{7,10,1,5,8},{8,0,4,7,10},{
    10,1,5,8,0}}

tonality[jewishword[0],4] =
{{0,4,7,10,1,5},{1,5,8,0,4,7},{4,7,10,1,5,8},{5,8,0,4,7,10},{7,10,1,5,8,0},{8,
    0,4,7,10,1},{10,1,5,8,0,4}}

tonality[jewishword[0],5] =
{{0,4,7,10,1,5,8},{1,5,8,0,4,7,10},{4,7,10,1,5,8,0},{5,8,0,4,7,10,1},{7,10,1,
    5,8,0,4},{8,0,4,7,10,1,5},{10,1,5,8,0,4,7}}

tonality[majorpentatonicword[0],1] =
{{0,4,9},{2,7,0},{4,9,2},{7,0,4},{9,2,7}}

tonality[majorpentatonicword[0],2] =
{{0,4,9,2},{2,7,0,4},{4,9,2,7},{7,0,4,9},{9,2,7,0}}

tonality[majorpentatonicword[0],3] =
{{0,4,9,2,7},{2,7,0,4,9},{4,9,2,7,0},{7,0,4,9,2},{9,2,7,0,4}}

tonality[minorpentatonicword[0],1] =
{{3,7,0},{5,10,3},{7,0,5},{10,3,7},{0,5,10}}

tonality[minorpentatonicword[0],2] =
{{3,7,0,5},{5,10,3,7},{7,0,5,10},{10,3,7,0},{0,5,10,3}}

tonality[minorpentatonicword[0],3] =
{{3,7,0,5,10},{5,10,3,7,0},{7,0,5,10,3},{10,3,7,0,5},{0,5,10,3,7}}

tonality[bluesword[0],1] =
{{0,5,7},{3,6,10},{5,7,0},{6,10,3},{7,0,5},{10,3,6}}

tonality[esatonalword[0],1] =
{{0,4,8},{2,6,10},{4,8,0},{6,10,2},{8,0,4},{10,2,6}}

tonality[augmentedword[0],1] =
{{0,4,8},{3,7,11},{4,8,0},{7,11,3},{8,0,4},{11,3,7}}

tonality[halfwholediminishedword[0],1] =
{{0,3,6},{1,4,7},{3,6,9},{4,7,10},{6,9,0},{7,10,1},{9,0,3},{10,1,4}}

tonality[halfwholediminishedword[0],2] =
{{0,3,6,9},{1,4,7,10},{3,6,9,0},{4,7,10,1},{6,9,0,3},{7,10,1,4},{9,0,3,6},{10,
    1,4,7}}

tonality[bebopmajorword[0],1] =
{{0,4,7},{2,5,8},{4,7,9},{5,8,11},{7,9,0},{8,11,2},{9,0,4},{11,2,5}}

tonality[bebopmajorword[0],2] =
{{0,4,7,9},{2,5,8,11},{4,7,9,0},{5,8,11,2},{7,9,0,4},{8,11,2,5},{9,0,4,7},{11,
    2,5,8}}

tonality[bebopdominant[0],1] =
{{0,4,7},{2,5,9},{4,7,10},{5,9,11},{7,10,0},{9,11,2},{10,0,4},{11,2,5}}

tonality[bebopdominant[0],2] =
{{0,4,7,10},{2,5,9,11},{4,7,10,0},{5,9,11,2},{7,10,0,4},{9,11,2,5},{10,0,4,
    7},{11,2,5,9}}

tonality[chromaticword[0],1] =
{{0,2,4},{1,3,5},{2,4,6},{3,5,7},{4,6,8},{5,7,9},{6,8,10},{7,9,11},{8,10,0},{
    9,11,1},{10,0,2},{11,1,3}}

tonality[chromaticword[0],2] =
{{0,2,4,6},{1,3,5,7},{2,4,6,8},{3,5,7,9},{4,6,8,10},{5,7,9,11},{6,8,10,0},{7,
    9,11,1},{8,10,0,2},{9,11,1,3},{10,0,2,4},{11,1,3,5}}

tonality[chromaticword[0],3] =
{{0,2,4,6,8},{1,3,5,7,9},{2,4,6,8,10},{3,5,7,9,11},{4,6,8,10,0},{5,7,9,11,1},{
    
6,8,10,0,2},{7,9,11,1,3},{8,10,0,2,4},{9,11,1,3,5},{10,0,2,4,6},{11,1,3,5,
    7}}

tonality[chromaticword[0],4] =
{{0,2,4,6,8,10},{1,3,5,7,9,11},{2,4,6,8,10,0},{3,5,7,9,11,1},{4,6,8,10,0,2},{
    5,7,9,11,1,3},{6,8,10,0,2,4},{7,9,11,1,3,5},{8,10,0,2,4,6},{9,11,1,3,5,
    7},{10,0,2,4,6,8},{11,1,3,5,7,9}}
\end{verbatim}

\bigskip
\begin{example}
\end{example}
THE II-V-I PROGRESSION

Let us consider the II-V-I progression, the most usual progression
in Jazz \cite{Levine-95}:
\begin{verbatim}

harmonicwordintermofdegrees[tonality[majorword[0],1],{2,5,1}] =
{{2,5,9},{7,11,2},{0,4,7}}

harmonicwordintermofdegrees[tonality[majorword[0],2],{2,5,1}] =
{{2,5,9,0},{7,11,2,5},{0,4,7,11}}

harmonicwordintermofdegrees[tonality[majorword[0],3],{2,5,1}] =
{{2,5,9,0,4},{7,11,2,5,9},{0,4,7,11,2}}

harmonicwordintermofdegrees[tonality[majorword[0],4],{2,5,1}] =
{{2,5,9,0,4,7},{7,11,2,5,9,0},{0,4,7,11,2,5}}

harmonicwordintermofdegrees[tonality[majorword[0],5],{2,5,1}] =
{{2,5,9,0,4,7,11},{7,11,2,5,9,0,4},{0,4,7,11,2,5,9}}

harmonicwordintermofdegrees[tonality[minorword[0],1],{2,5,1}] =
{{2,5,8},{7,10,2},{0,3,7}}

harmonicwordintermofdegrees[tonality[minorword[0],2],{2,5,1}] =
{{2,5,8,0},{7,10,2,5},{0,3,7,10}}

harmonicwordintermofdegrees[tonality[minorword[0],3],{2,5,1}] =
{{2,5,8,0,3},{7,10,2,5,8},{0,3,7,10,2}}

harmonicwordintermofdegrees[tonality[minorword[0],4],{2,5,1}] =
{{2,5,8,0,3,7},{7,10,2,5,8,0},{0,3,7,10,2,5}}

harmonicwordintermofdegrees[tonality[minorword[0],5],{2,5,1}] =
{{2,5,8,0,3,7,10},{7,10,2,5,8,0,3},{0,3,7,10,2,5,8}}

harmonicwordintermofdegrees[tonality[harmonicminorword[0],1],{2,5,1}]
= {{2,5,8},{7,11,2},{0,3,7}}

harmonicwordintermofdegrees[tonality[harmonicminorword[0],2],{2,5,1}]
= {{2,5,8,0},{7,11,2,5},{0,3,7,11}}

harmonicwordintermofdegrees[tonality[harmonicminorword[0],3],{2,5,1}]
= {{2,5,8,0,3},{7,11,2,5,8},{0,3,7,11,2}}

harmonicwordintermofdegrees[tonality[harmonicminorword[0],4],{2,5,1}]
= {{2,5,8,0,3,7},{7,11,2,5,8,0},{0,3,7,11,2,5}}

harmonicwordintermofdegrees[tonality[harmonicminorword[0],5],{2,5,1}]
= {{2,5,8,0,3,7,11},{7,11,2,5,8,0,3},{0,3,7,11,2,5,8}}

harmonicwordintermofdegrees[tonality[dorianword[0],1],{2,5,1}] =
{{2,5,9},{7,10,2},{0,3,7}}

harmonicwordintermofdegrees[tonality[dorianword[0],2],{2,5,1}] =
{{2,5,9,0},{7,10,2,5},{0,3,7,10}}

harmonicwordintermofdegrees[tonality[dorianword[0],3],{2,5,1}] =
{{2,5,9,0,3},{7,10,2,5,9},{0,3,7,10,2}}

harmonicwordintermofdegrees[tonality[dorianword[0],4],{2,5,1}] =
{{2,5,9,0,3,7},{7,10,2,5,9,0},{0,3,7,10,2,5}}

harmonicwordintermofdegrees[tonality[dorianword[0],5],{2,5,1}] =
{{2,5,9,0,3,7,10},{7,10,2,5,9,0,3},{0,3,7,10,2,5,9}}

harmonicwordintermofdegrees[tonality[phrigianword[0],1],{2,5,1}] =
{{1,5,8},{7,10,1},{0,3,7}}

harmonicwordintermofdegrees[tonality[phrigianword[0],2],{2,5,1}] =
{{1,5,8,0},{7,10,1,5},{0,3,7,10}}

harmonicwordintermofdegrees[tonality[phrigianword[0],3],{2,5,1}] =
{{1,5,8,0,3},{7,10,1,5,8},{0,3,7,10,1}}

harmonicwordintermofdegrees[tonality[phrigianword[0],4],{2,5,1}] =
{{1,5,8,0,3,7},{7,10,1,5,8,0},{0,3,7,10,1,5}}

harmonicwordintermofdegrees[tonality[phrigianword[0],5],{2,5,1}] =
{{1,5,8,0,3,7,10},{7,10,1,5,8,0,3},{0,3,7,10,1,5,8}}

harmonicwordintermofdegrees[tonality[lydianword[0],1],{2,5,1}] =
{{2,6,9},{7,11,2},{0,4,7}}

harmonicwordintermofdegrees[tonality[lydianword[0],2],{2,5,1}] =
{{2,6,9,0},{7,11,2,6},{0,4,7,11}}

harmonicwordintermofdegrees[tonality[lydianword[0],3],{2,5,1}] =
{{2,6,9,0,4},{7,11,2,6,9},{0,4,7,11,2}}

harmonicwordintermofdegrees[tonality[lydianword[0],4],{2,5,1}] =
{{2,6,9,0,4,7},{7,11,2,6,9,0},{0,4,7,11,2,6}}

harmonicwordintermofdegrees[tonality[lydianword[0],5],{2,5,1}] =
{{2,6,9,0,4,7,11},{7,11,2,6,9,0,4},{0,4,7,11,2,6,9}}

harmonicwordintermofdegrees[tonality[mixolydianword[0],1],{2,5,1}]
= {{2,5,9},{7,10,2},{0,4,7}}

harmonicwordintermofdegrees[tonality[mixolydianword[0],2],{2,5,1}]
= {{2,5,9,0},{7,10,2,5},{0,4,7,10}}

harmonicwordintermofdegrees[tonality[mixolydianword[0],3],{2,5,1}]
= {{2,5,9,0,4},{7,10,2,5,9},{0,4,7,10,2}}

harmonicwordintermofdegrees[tonality[mixolydianword[0],4],{2,5,1}]
= {{2,5,9,0,4,7},{7,10,2,5,9,0},{0,4,7,10,2,5}}

harmonicwordintermofdegrees[tonality[mixolydianword[0],5],{2,5,1}]
= {{2,5,9,0,4,7,10},{7,10,2,5,9,0,4},{0,4,7,10,2,5,9}}

harmonicwordintermofdegrees[tonality[locrianword[0],1],{2,5,1}] =
{{1,5,8},{6,10,1},{0,3,6}}

harmonicwordintermofdegrees[tonality[locrianword[0],2],{2,5,1}] =
{{1,5,8,0},{6,10,1,5},{0,3,6,10}}

harmonicwordintermofdegrees[tonality[locrianword[0],3],{2,5,1}] =
{{1,5,8,0,3},{6,10,1,5,8},{0,3,6,10,1}}

harmonicwordintermofdegrees[tonality[locrianword[0],4],{2,5,1}] =
{{1,5,8,0,3,6},{6,10,1,5,8,0},{0,3,6,10,1,5}}

harmonicwordintermofdegrees[tonality[locrianword[0],5],{2,5,1}] =
{{1,5,8,0,3,6,10},{6,10,1,5,8,0,3},{0,3,6,10,1,5,8}}

harmonicwordintermofdegrees[tonality[tziganminorword[0],1],{2,5,1}]
= {{2,6,8},{7,11,2},{0,3,7}}

harmonicwordintermofdegrees[tonality[tziganminorword[0],2],{2,5,1}]
= {{2,6,8,0},{7,11,2,6},{0,3,7,11}}

harmonicwordintermofdegrees[tonality[tziganminorword[0],3],{2,5,1}]
= {{2,6,8,0,3},{7,11,2,6,8},{0,3,7,11,2}}

harmonicwordintermofdegrees[tonality[tziganminorword[0],4],{2,5,1}]
= {{2,6,8,0,3,7},{7,11,2,6,8,0},{0,3,7,11,2,6}}

harmonicwordintermofdegrees[tonality[tziganminorword[0],5],{2,5,1}]
= {{2,6,8,0,3,7,11},{7,11,2,6,8,0,3},{0,3,7,11,2,6,8}}

harmonicwordintermofdegrees[tonality[jewishword[0],1],{2,5,1}] =
{{1,5,8},{7,10,1},{0,4,7}}

harmonicwordintermofdegrees[tonality[jewishword[0],2],{2,5,1}] =
{{1,5,8,0},{7,10,1,5},{0,4,7,10}}

harmonicwordintermofdegrees[tonality[jewishword[0],3],{2,5,1}] =
{{1,5,8,0,4},{7,10,1,5,8},{0,4,7,10,1}}

harmonicwordintermofdegrees[tonality[jewishword[0],4],{2,5,1}] =
{{1,5,8,0,4,7},{7,10,1,5,8,0},{0,4,7,10,1,5}}

harmonicwordintermofdegrees[tonality[jewishword[0],5],{2,5,1}] =
{{1,5,8,0,4,7,10},{7,10,1,5,8,0,4},{0,4,7,10,1,5,8}}

harmonicwordintermofdegrees[tonality[majorpentatonicword[0],1],{2,5,1}]
= {{2,7,0},{9,2,7},{0,4,9}}

harmonicwordintermofdegrees[tonality[majorpentatonicword[0],2],{2,5,1}]
= {{2,7,0,4},{9,2,7,0},{0,4,9,2}}

harmonicwordintermofdegrees[tonality[majorpentatonicword[0],3],{2,5,1}]
= {{2,7,0,4,9},{9,2,7,0,4},{0,4,9,2,7}}

harmonicwordintermofdegrees[tonality[minorpentatonicword[0],1],{2,5,1}]
= {{5,10,3},{0,5,10},{3,7,0}}

harmonicwordintermofdegrees[tonality[minorpentatonicword[0],2],{2,5,1}]
= {{5,10,3,7},{0,5,10,3},{3,7,0,5}}

harmonicwordintermofdegrees[tonality[minorpentatonicword[0],3],{2,5,1}]
= {{5,10,3,7,0},{0,5,10,3,7},{3,7,0,5,10}}

harmonicwordintermofdegrees[tonality[bluesword[0],1],{2,5,1}] =
{{3,6,10},{7,0,5},{0,5,7}}

harmonicwordintermofdegrees[tonality[esatonalword[0],1],{2,5,1}] =
{{2,6,10},{8,0,4},{0,4,8}}

harmonicwordintermofdegrees[tonality[augmentedword[0],1],{2,5,1}]
= {{3,7,11},{8,0,4},{0,4,8}}

harmonicwordintermofdegrees[tonality[halfwholediminishedword[0],1],{2,5,1}]
= {{1,4,7},{6,9,0},{0,3,6}}

harmonicwordintermofdegrees[tonality[halfwholediminishedword[0],2],{2,5,1}]
= {{1,4,7,10},{6,9,0,3},{0,3,6,9}}

harmonicwordintermofdegrees[tonality[wholehalfdiminishedword[0],1],{2,5,1}]
= {{2,5,8},{6,9,0},{0,3,6}}

harmonicwordintermofdegrees[tonality[wholehalfdiminishedword[0],2],{2,5,1}]
= {{2,5,8,11},{6,9,0,3},{0,3,6,9}}

harmonicwordintermofdegrees[tonality[wholetonediminishedword[0],1],{2,5,1}]
= {{1,4,8},{6,10,1},{0,3,6}}

harmonicwordintermofdegrees[tonality[wholetonediminishedword[0],2],{2,5,1}]
= {{1,4,8,0},{6,10,1,4},{0,3,6,10}}

harmonicwordintermofdegrees[tonality[wholetonediminishedword[0],3],{2,5,1}]
= {{1,4,8,0,3},{6,10,1,4,8},{0,3,6,10,1}}

harmonicwordintermofdegrees[tonality[wholetonediminishedword[0],4],{2,5,1}]
= {{1,4,8,0,3,6},{6,10,1,4,8,0},{0,3,6,10,1,4}}

harmonicwordintermofdegrees[tonality[wholetonediminishedword[0],5],{2,5,1}]
= {{1,4,8,0,3,6,10},{6,10,1,4,8,0,3},{0,3,6,10,1,4,8}}

harmonicwordintermofdegrees[tonality[bebopmajorword[0],1],{2,5,1}]
= {{2,5,8},{7,9,0},{0,4,7}}

harmonicwordintermofdegrees[tonality[bebopmajorword[0],2],{2,5,1}]
= {{2,5,8,11},{7,9,0,4},{0,4,7,9}}

harmonicwordintermofdegrees[tonality[bebopdominant[0],1],{2,5,1}]
= {{2,5,9},{7,10,0},{0,4,7}}

harmonicwordintermofdegrees[tonality[bebopdominant[0],2],{2,5,1}]
= {{2,5,9,11},{7,10,0,4},{0,4,7,10}}

harmonicwordintermofdegrees[tonality[chromaticword[0],1],{2,5,1}]
= {{1,3,5},{4,6,8},{0,2,4}}

harmonicwordintermofdegrees[tonality[chromaticword[0],2],{2,5,1}]
= {{1,3,5,7},{4,6,8,10},{0,2,4,6}}

harmonicwordintermofdegrees[tonality[chromaticword[0],3],{2,5,1}]
= {{1,3,5,7,9},{4,6,8,10,0},{0,2,4,6,8}}

harmonicwordintermofdegrees[tonality[chromaticword[0],4],{2,5,1}]
= {{1,3,5,7,9,11},{4,6,8,10,0,2},{0,2,4,6,8,10}}

\end{verbatim}

\bigskip

Given two tonalities $ t_{1} \, , \,  t_{2} \; \in \;
{\mathcal{T}}  $:
\begin{definition}
\end{definition}
PIVOTAL DEGREES OF $ t_{1} $ AND $ t_{2} $:
\begin{equation}
    {\mathcal{P}}(  t_{1} \, , \,   t_{2} ) \; := \; Range(
    chord^{(level(t_{1})} ( t_{1} )) \, \bigcap \,  Range(
    chord^{(level(t_{2}))} ( t_{2} ))
\end{equation}
Such a concept, essential as to Modulation Theory, may be easily
computed through the following Mathematica expression from the
notebook of section\ref{sec:Mathematica's notebook Mathematical
Music Toolkit}:

\begin{verbatim}

pivotaldegrees[t1_,t2_]:=Intersection[t1,t2]

(*** "pivotaldegreesintermofdegrees[t1,t2]" give a list of two
elements,  the first being the pivotal degrees of "t1" and "t2"
expressed as degrees of "t1", the second being the pivotal degrees
of "t1" and "t2" expressed as degrees of t2" ***)

pivotaldegreesintermofdegrees[t1_,t2_]:={
    Table[degreeofachordinatonality[Part[pivotaldegrees[t1,t2],i],t1],{i,1,
        Length[pivotaldegrees[t1,t2]]}],
    Table[degreeofachordinatonality[Part[pivotaldegrees[t1,t2],i],t2],{i,1,
        Length[pivotaldegrees[t1,t2]]}]}

\end{verbatim}
It is important to remark that, obviously, the pivotal degrees
depends crucially from the level at which the involved tonalities
are built:

\begin{example} \label{ex:pivotal degrees of major tonalities at different 
levels along a clockwise step in fifths' cycle}
\end{example}
PIVOTAL DEGREES OF MAJOR TONALITIES AT DIFFERENT LEVELS ALONG A
CLOCKWISE STEP IN THE FIFTHS' CYCLE
\begin{multline}
  pivotaldegreesintermofdegrees[tonality[majorword[0],1],
  tonality[majorword[7],1]] \; = \\
   \{\{1,3,5,6\},\{4,6,1,2\}\}
\end{multline}
\begin{multline}
   pivotaldegreesintermofdegrees[tonality[majorword[0],2],
  tonality[majorword[7],2]] \; = \\
   \{\{1,3,6\},\{4,6,2\}\}
\end{multline}
\begin{multline}
   pivotaldegreesintermofdegrees[tonality[majorword[0],3],
  tonality[majorword[7],3]] \; = \\
   \{\{1,6\},\{4,2\}\}
\end{multline}
\begin{multline}
   pivotaldegreesintermofdegrees[tonality[majorword[0],4],
  tonality[majorword[7],4]]\; = \\
\{\{6\},\{2\}\}
\end{multline}
\begin{multline}
   pivotaldegreesintermofdegrees[tonality[majorword[0],5],
  tonality[majorword[7],5]]\; = \\
  \emptyset
\end{multline}
so that, as to Classical Harmony, the $ I_{C} \, = \, IV_{G} \, =
\, C $, the $ III_{C} \, = \, VI_{G} \, = \, E m $, the $ V_{C} \,
= \, I_{G} \, = \, G7 $, the $ VI_{C} \, = \, II_{G} \, = \, A m $
are four different streets one can follow to modulate from C Major
to G Major, a thing almost  forbidden in that context.

In Jazz Harmony, where parallel motion is allowed, tonalities are
defined at a level $ > 1 $ with the number of possible streets
lowering with the considered level: at second level one still has
that \footnote{As to the nomenclature for chords I will adopt that
of American Jazz Theory  exposed in the appendix L.2 "Third Chain
Classes" of \cite{Mazzola-02a}} $ I_{C}^{7} \, = \, IV_{G}^{7} \,
= C 7 + $, $ III_{C}^{7} \, = \, VI_{G}^{7} \, = \, E m 7 $, at
third level remains only that $ I_{C}^{9} \, = \, IV_{G}^{9} \, =
\, C 7 + / 9 $ and  that $ VI_{C}^{9} \, = \, II_{G}^{9} \, = \, A
m 9 $, at fourth level remains only that $  VI_{C}^{11} \, = \,
II_{G}^{11} \, = \, A m 11 $.

\smallskip

\begin{example} \label{ex:pivotal degrees of jewish tonalities at different 
levels along a clockwise step in fifths' cycle}
\end{example}
PIVOTAL DEGREES OF JEWISH TONALITIES AT DIFFERENT LEVELS ALONG A
CLOCKWISE STEP IN THE FIFTHS' CYCLE

\begin{verbatim}
pivotaldegreesintermofdegrees[tonality[jewishword[0],1],
  tonality[jewishword[7],1]] = {{4},{7}}

pivotaldegreesintermofdegrees[tonality[jewishword[0],2],
  tonality[jewishword[7],2]] = {{},{}}

\end{verbatim}
so that parallel motion of a clockwise step in the fifths' cycle
has strongly less chances for jewish tonalities than for the major
tonalities discussed in the example \ref{ex:pivotal degrees of
major tonalities at different levels along a clockwise step in
fifths' cycle}: at the first level one has the only pivotal degree
$ IV_{C} \; = \; VII_{G} \; = \; Fm $ while no pivotal degree
exists already at the second level.

\smallskip

Let us now pass to analyze the concept of cadence: with an eye to
its application as to Modulation Theory, Mazzola (cfr. the $
26^{th} $ chapter "Cadences" of \cite{Mazzola-02a}) has correctly
catched the structural peculiarity of such harmonic words: they
identify a tonality among a suitable class of tonalities.

What is important, at this point, to stress is how such a concept
critically depends by a context of tonalities among which the
cadential one has to identify its tonality.

Mazzola's definition of cadence may be easily translated also as
to our more general definition of tonality, namely the
definition\ref{def:tonality of a non-repetitive word at a certain
level}, in the following way:

given a tonality $ t \in {\mathcal{T}}$ and a context $
{\mathcal{T}}_{context} \subseteq {\mathcal{T}} $
\begin{definition} \label{cadences of a tonality w.r.t. to a context}
\end{definition}
CADENCES OF $  t $ W.R.T. THE CONTEXT $ {\mathcal{T}}_{context} $
\begin{equation}
     { \mathcal{C}}(t \, , \, {\mathcal{T}}_{context} ) \; := \; \{ c  \, 
\in \, {\mathcal{HW}} (t) \; :
     \;
      ( \,  c  \, \in \, {\mathcal{HW}} (t) \; \Rightarrow \;  t \, = \,  u 
\,
      ) \, \, \forall u \in {\mathcal{T}}_{context} \}
\end{equation}
To understand how crucial is the role of the context in the
definition\ref{cadences of a tonality w.r.t. to a context} let us
observe that:
\begin{theorem} \label{th:no-go theorem for cadences in a multi-mode 
context}
\end{theorem}
NO-GO THEOREM  FOR CADENCES IN A MULTI-MODE CONTEXT

\begin{hypothesis}
\end{hypothesis}
\begin{equation}
    \vec{x}_{1} \, , \, \vec{x}_{2} \, \in \, {\mathbb{Z}}_{12}^{\star}
\end{equation}
\begin{equation}
    {\mathcal{T}}_{context} \subset {\mathcal{T}} \, : \,
    tonality(  \vec{x}_{1} , n ) , tonality(  \vec{x}_{2} ,n ) \in
     {\mathcal{T}}_{context}
\end{equation}
\begin{equation}
    \exists i \, : \vec{x}_{2} = mode(  \vec{x}_{1} , i)
\end{equation}
\begin{thesis}
\end{thesis}
\begin{equation}
  { \mathcal{C}}[tonality(\vec{x}_{1}\, , \,n) \, , \, 
{\mathcal{T}}_{context} ] \; = \;
  \emptyset
\end{equation}
\begin{proof}
Since:
\begin{equation}
  chord ( \vec{x}_{2} , j , n ) \, = \, chord ( \vec{x}_{2} , 
Mod[j+i,|\vec{x_{1}}|]    , n
  ) \; \; \forall j = 1 , \cdots , |\vec{x_{2}}|
\end{equation}
one has that:
\begin{equation}
    hw \in tonality[ \vec{x_{1}} , n ] \; \Rightarrow \; hw \in tonality[ 
\vec{x_{2}} , n
    ] \; \Rightarrow \; hw \notin { \mathcal{C}}(tonality[\vec{x}_{1}\, , 
\,n] \, , \, {\mathcal{T}}_{context} )
\end{equation}
\end{proof}

An immediate consequence of the theorem \ref{th:no-go theorem for
cadences in a multi-mode context} is the following:
\begin{corollary}
\end{corollary}
\begin{eqnarray}
  { \mathcal{C}}(t \, , \, {\mathcal{T}}_{greg} ) \; & = &\;
  \emptyset \; \; \forall t \in {\mathcal{T}}_{greg}  \\
{ \mathcal{C}}(t \, , \, {\mathcal{T}}_{cl} ) \;& = &\;
  \emptyset \; \; \forall t \in {\mathcal{T}}_{cl}
\end{eqnarray}

A someway natural choice of the context respect to which to
analyze whether an harmonic word is a tonality are translational
equivalence classes, a fact this one we will analyze more deeply
when we will discuss Modulation Theory.

Such an idea may be formalized, anyway, in the following way:

given a word $ \vec{x} \in \Sigma^{\star} $ a natural number $ n
\in {\mathbb{N}} \, : \, n < maxlevel( \vec{x} ) $
\begin{definition} \label{def:natural context of a word at a certain level}
\end{definition}
NATURAL CONTEXT OF $ tonality[ \vec{x} , n ] $
\begin{equation}
    {\mathcal{T}}_{n.c.} ( \vec{x} , n ) \; := \; \{
    tonality(\vec{y} \, , \, n) \: : \: \vec{x} \, \sim_{T}
    \vec{y} \}
\end{equation}

Given  a tonality $ t \in {\mathcal{T}} $ and a context $
{\mathcal{T}}_{context} \subseteq {\mathcal{T}} $
\begin{definition}
\end{definition}
MINIMAL CADENCES OF $  t $ W.R.T. THE CONTEXT $
{\mathcal{T}}_{context} $
\begin{equation}
    {\mathcal{MC}} ( t \, , \, {\mathcal{T}}_{context} ) \; := \; \{ c_{1}  
\, \in \, {
\mathcal{C}}(t, {\mathcal{T}}_{context}) \, : \,  \nexists c_{2}
\in { \mathcal{C}}(t,{\mathcal{T}}_{context}) , c_{2} \, <_{p} \,
c_{1} \}
\end{equation}

Cadences and minimal cadences with respect to suitable contexts
may be easily computed through the following expressions from the
Mathematica notebook of section \ref{sec:Mathematica's notebook
Mathematical Music Toolkit}:

\begin{verbatim}

setofthemajortonalities[n_]:=Table[tonality[majorword[i],n],{i,0,11}];

setoftheminortonalities[n_]:=Table[tonality[minorword[i],n],{i,0,11}];

setoftheclassicaltonalities[n_]:=
  Union[setofthemajortonalities[n],setoftheminortonalities[n]];

setofthegregoriantonalities[n_]:=
  Flatten[Table[tonality[mode[majorword[i],j],n],{j,1,7},{i,0,11}],1];

setofthemazzolatonalities[n_]:=
  
Table[tonality[Part[mazzolawords,i],n],{i,1,Length[nonrepetitivewords[7]]}];

setofthejewishtonalities[n_]:=Table[tonality[jewishword[i],n],{i,0,11}];

truthQ[x_]:=Equal[x,True]

cadenceQ[hw_,setoftonalities_]:= If[Length[ Select[Table[
            tonalitymembershipQ[hw,Part[setoftonalities,i]],{i,1,
              Length[setoftonalities]}],truthQ]]==1,True,False]

subharmonicwords[hw_]:=LexicographicSubsets[hw]

cadences[t_,n_,setoftonalities_]:=
  generalizedselect[harmonicwordsupto[t,n],cadenceQ,setoftonalities]

cadencesintermofdegrees[t_,n_,setoftonalities_]:=
  Table[degreeofachordinatonality[
      Part[Part[cadences[t,n,setoftonalities],i],j],t],{i,1,
      Length[cadences[t,n,setoftonalities] ]},{j,1,
      Length[ Part[cadences[t,n,setoftonalities],i] ]}]

minimalcadenceQ[hw_,setoftonalities_]:=And[cadenceQ[hw,setoftonalities],
        
Equal[Table[cadenceQ[prefix[hw,i],setoftonalities],{i,1,Length[hw]-1}],
      Table[False,{Length[hw]-1}]]]

minimalcadences[t_,n_,setoftonalities_]:=
  generalizedselect[harmonicwords[t,n],minimalcadenceQ,setoftonalities]

minimalcadencesintermofdegrees[t_,n_,setoftonalities_]:=
  Table[degreeofachordinatonality[
      Part[Part[minimalcadences[t,n,setoftonalities],i],j],t],{i,1,
      Length[minimalcadences[t,n,setoftonalities] ]},{j,1,
      Length[ Part[minimalcadences[t,n,setoftonalities],i] ]}]
\end{verbatim}

\begin{remark} \label{rem:the perfect and plagal harmonic progressions}
\end{remark}
THE PERFECT AND PLAGAL HARMONIC PROGRESSIONS:

Though expressing the idea that the role of a cadence is, in
certain contexts, that of declaring the underlying tonality,
hystorical tradition has led most of the more common Harmony's
manuals \cite{Dubois-01}, \cite{Piston-89}, \cite{Schoenberg-02}
to speak about the so-called perfect cadence $ V-I $ and plagal
cadence $ IV -I$. That such harmonic words are not cadences even
with respect to the little class of major tonalities may be easily
verified observing that, for example, the harmonic word $ ( V_{C}
\, , \, I_{C}) \, = \, {{7,11,2},{0,4,7}} \, = \, ( G \, , \, C )
$ may be seen also as $ ( I_{G} \, , \, IV_{G} )$ and that the
harmonic word $ ( IV_{C} \, , \, I_{C} ) \, = \, {{5,9,0},{0,4,7}}
\, = \; (F \, , \, C) $ may be seen also as $ (I_{F} \, , \, V_{F}
)$. Such a situation, anyway, changes if one takes in
consideration major tonalities at a livel greater than the first.

\smallskip

\begin{example}
\end{example}
CADENCES AND MINIMAL CADENCES OF THE C MAJOR TONALITY W.R.T. THE
CONTEXT OF MAJOR TONALITIES AT DIFFERENT LEVELS

Obviously:
\begin{verbatim}
cadencesintermofdegrees[tonality[majorword[0],1],1,setoftheclassicaltonalities[1]]
=  {}
\end{verbatim}
since any one-letter harmonic word in a major tonality x may be
also seen as a one letter harmonic word in the relative minor
tonality.

Restricting the set of the tonalities to the major ones one finds
that:
\begin{verbatim}
cadencesintermofdegrees[tonality[majorword[0],1],1,setofthemajortonalities[1]]
= {{7}}
\end{verbatim}
showing that $ VII_{C} $ is the only one-letter harmonic word to
be a cadence w.r.t. the set of major tonalities at first level
since it is the only diminished triad $( \{11,2,5\} \, = \,  B m 5
- ) $ in tonality[majorword[0],1], the other degrees being major
triads $( \{0,4,7\} = C \, , \,  \{5,9,0\} = F \,  , \, \{7,11,2\}
= G) $ or minor triads $  (\{2,5,9\} = Dm \,  , \, \{4,7,11\} = Em
\, , \, \{9,0,4\}= Am) $.

Raising the level of the tonalities, anyway, the situation
changes:
\begin{verbatim}
cadencesintermofdegrees[tonality[majorword[0],2],1,setofthemajortonalities[2]]
{{5},{7}}
\end{verbatim}
since at the second level $ V_{C}^{(7)} \, = \,  \{7,11,2,5\} \, =
\, G 7 \, = \, $ is the only first kind major seventh chord, $
I_{C}^{(7)} \, = \, \{0,4,7,11\}\, = \, C 7 + $ and $ IV_{C}^{(7)}
\, = \, \{5,9,0,4\} \, = \, F 7 + $ being fourth kind seventh
chords.

Passing to the third level:
\begin{verbatim}
cadencesintermofdegrees[tonality[majorword[0],3],1,setofthemajortonalities[3]]
= {{3},{5},{7}}
\end{verbatim}
one obtains a new one-letter cadential set, i.e. $ III_{C}^{(9)}
\, = \, \{4,7,11,2,5\} \, = \, E m 9 - $, that is a ninth chord of
different kind w.r.t. $ II_{C}^{(9)} \, = \,  \{2,5,9,0,4\} \, =
\, D m 9 $ and $ VI_{C} ^{(9)} \, = \, \{9,0,4,7,11\} \, = \, Am9
$.

Passing to the fourth level:
\begin{verbatim}
cadencesintermofdegrees[tonality[majorword[0],4],1,setofthemajortonalities[4]]
= {{1},{3},{4},{5},{7}}
\end{verbatim}
one obtains two new one-letter cadences, i.e. $ I_{C}^{(11)} \, =
\,  \{0,4,7,11,2,5\} \, = \,C7+/11 $  and $ IV_{C}^{(11)} \, = \,
\{5,9,0,4,7,11\} \, = \, C7+ / 11+ $ that are $ 11^{th} $ chords
of different type.

Passing, finally, to the fifth level:
\begin{verbatim}
cadencesintermofdegrees[tonality[majorword[0],5],1,setofthemajortonalities[5]]
= {{1},{2},{3},{4},{5},{6},{7}}
\end{verbatim}
one can that all the one-letter harmonic words become cadential.

Let us now pass to analyze harmonic words of length two,
restricting the analysis to the minimal cadences.

At the $ 1^{th} $ level:
\begin{verbatim}
minimalcadencesintermofdegrees[tonality[majorword[0],1],2,
  setofthemajortonalities[1]] =
{{1,7},{2,3},{2,5},{2,7},{3,2},{3,4},{3,7},{4,3},{4,5},{4,7},{5,2},{5,4},{5,
    7},{6,7}}
\end{verbatim}
giving us confirmation of what we saw in the remark\ref{rem:the
perfect and plagal harmonic progressions}.

At the  $ 2^{th} $ level:
\begin{verbatim}
minimalcadencesintermofdegrees[tonality[majorword[0],2],2,
  setofthemajortonalities[2]] =
{{1,2},{1,4},{1,5},{1,7},{2,1},{2,3},{2,5},{2,7},{3,2},{3,4},{3,5},{3,7},{4,
    1},{4,3},{4,5},{4,7},{6,5},{6,7}}
\end{verbatim}
one sees that the plagal harmonic word $ ( IV_{C}^{(7)} ,
I_{C}^{7} ) \, = \, \{\{5,9,0,4\},\{0,4,7,11\}\} \, = \, ( F7+ ,
C7+) $ becomes a cadence.

At the  $ 3^{th} $ level:
\begin{verbatim}
minimalcadencesintermofdegrees[tonality[majorword[0],3],2,
  setofthemajortonalities[3]] =
{{1,2},{1,3},{1,4},{1,5},{1,7},{2,1},{2,3},{2,5},{2,6},{2,7},{4,1},{4,3},{4,
    5},{4,6},{4,7},{6,2},{6,3},{6,4},{6,5},{6,7}}
\end{verbatim}
this is still true, while at the $ 4^{th} $ level:
\begin{verbatim}
minimalcadencesintermofdegrees[tonality[majorword[0],4],2,
  setofthemajortonalities[4]] =
{{2,1},{2,3},{2,4},{2,5},{2,6},{2,7},{6,1},{6,2},{6,3},{6,4},{6,5},{6,7}}
\end{verbatim}
this is no longer true, since, though remaining a cadence:
\begin{verbatim}
cadenceQ[harmonicwordintermofdegrees[tonality[majorword[0],4],{4,1}],
  setofthemajortonalities[4]] = True
\end{verbatim}
the plagal cadence is no more minimal since we saw previously
that, at this level $ VI_{C}^{7} $ is a cadence by itself. For the
same reason:
\begin{verbatim}
minimalcadencesintermofdegrees[tonality[majorword[0],5],2,
  setofthemajortonalities[5]] = {}
\end{verbatim}
no cadential set at the $ 5^{th} $ level can be minimal.

\medskip

\begin{example}
\end{example}
CADENCES AND MINIMAL CADENCES OF THE C JEWISH TONALITY  W.R.T. THE
CONTEXT OF JEWISH TONALITIES AT DIFFERENT LEVELS

At the $ 1^{th} $ level one has that:
\begin{verbatim}
cadencesintermofdegrees[tonality[jewishword[0],1],1,setofthejewishtonalities[1]]
= {{6}}
\end{verbatim}
so that the only one-letter cadential set of the C Jewish tonality
is $ VI_{C} \, = \, \{8,0,4\} \, = \, G^{\sharp} 5 + $.

It is sufficient, anyway, to pass to the $ 2^{th} $ level:
\begin{verbatim}
cadencesintermofdegrees[tonality[jewishword[0],2],1,
  setofthejewishtonalities[2]] = {{1},{2},{3},{4},{5},{6},{7}}
\end{verbatim}
that any one-letter harmonic-word becomes a cadence, such a
property continuing, clearly, to hold at higher levels.

Passing to harmonic words of two letters one has that:
\begin{verbatim}
minimalcadencesintermofdegrees[tonality[jewishword[0],1],2,
  setofthejewishtonalities[1]] =
{{1,2},{1,3},{1,4},{1,5},{1,6},{1,7},{2,1},{2,3},{2,4},{2,5},{2,6},{2,7},{3,1},
{3,2},{3,4},{3,5},{3,6},{3,7},{4,1},{4,2},{4,3},{4,5},{4,6},{4,7},{5,1},
{5,2},{5,3},{5,4},{5,6},{5,7},{7,1},{7,2},{7,3},{7,4},{7,5},{7,6}}
\end{verbatim}
It is interesting to note that, in particular,  that both the
perfect harmonic word $( V_{C} , I_{C} ) \, = \, ( \{7,10,1 \} ,
\{ 0,4,7 \} ) \, = \, ( G m 5 - , C ) $ and the plagal harmonic
word $( IV_{C} , I_{C} ) \, = \, ( \{5,8,0\} , \{0,4,7 \} \, = \,
( Fm , C ) $ are cadential.

\medskip

Let us now analyze the Muzzulini - Mazzola's theory of
modulations in the framework of Symmetry's Theory
\cite{Muzzulini-95}, \cite{Mazzola-02a}.

Given an arbitrary finite alphabet $ \Sigma $, a word $ \vec{x}
\in \Sigma^{\star} $ and a map $ g \, : \, \Sigma \rightarrow
\Sigma  $:
\begin{definition} \label{def;symmetry of a word}
\end{definition}
$ g $ IS A SYMMETRY OF $ \vec{x} $:
\begin{equation}
    \hat{g} ( \vec{x} ) \; = \; \vec{x}
\end{equation}
The symmetries of a word $ \vec{x} $ constitute a group I will
denote as $ SYM( \vec{x} ) $. For a systematic catalogation of the
symmetry groups of all the elements of $ \Sigma_{NR}^{\star} $ I
demand to the appendix L.1 "Chords and Third Chain Classes" of
\cite{Mazzola-02a}.

Given two  tonalities $ t_{1} \, , \, t_{2} \; \in \;
{\mathcal{T}} $ such that $ t_{2} \in {\mathcal{T}}_{n.c.} (
t_{1)} $, let's say $ chord^{-1} ( t_{2} ) \; = \hat{T}_{z} \;
chord^{-1} ( t_{2} ) $ :
\begin{definition} \label{def;Mazzola modulator}
\end{definition}
MAZZOLA MODULATOR FROM $ t_{1} $ TO  $ t_{2} $:

a map   $ g \, : \, chord^{-1} ( t_{1} )  \, \rightarrow \,
chord^{-1} ( t_{2}) $ of the form:
\begin{equation}
    g \; = \; T_{z} \circ h \; z \in {\mathbb{Z}}_{12} \, , \, h
\in SYM( chord^{-1} ( t_{1} ) )
\end{equation}

\begin{definition} \label{def:Mazzola modulation}
\end{definition}
MAZZOLA MODULATION FROM $ t_{1} $ TO $ t_{2} $:

a couple $ ( g \, , \, c ) $ such that:
\begin{enumerate}
    \item g is a modulator from $ t_{1} $ to $ t_{2} $
    \item
\begin{equation}
    c \; \in \;  {\mathcal{C}} [ t_{2} \, , \, {\mathcal{T}}_{n.c.} ( t_{2}
    )]
\end{equation}
\end{enumerate}
where $ {\mathcal{T}}_{n.c.} ( t_{2} )$ is shortcut notation
w.r.t. the definition \ref{def:natural context of a word at a
certain level} of evident meaning.

\begin{definition}
\end{definition}
MAZZOLA TONAL MUSICAL PIECES:
\begin{multline}
  {\mathcal{MP}}_{Maz} \; := \; \{ (hw_{1} \, m( t_{1} , t_{2} )  \,  hw_{2} 
\,   \cdots \,  hw_{n-1} \, m( t_{n-1} , t_{n} )  \,  hw_{n}  ) \; : \\
   t_{i} \,  \in  \,  {\mathcal{T}}_{Maz} \, , hw_{i} \in {\mathcal{HW}}( 
t_{i})  \, ,   m( t_{i-1} , t_{i} ) \,  \in \,  {\mathcal{M}}_{Maz}( t_{i-1} 
, t_{i} ) \, i = 1 , \cdots , n \, \\
     n \in {\mathbb{N}} \}
\end{multline}

The Muzzulini-Mazzola's theory, introduces, at this point the
notion of \emph{modulation quantum} based on the analogy between
the quanta mediating physical interaction in Particle Physics.

Such an analogy is, anyway, rather superficial since it  is based
simply on the visualization of Feynman diagram's in the
perturbative expansions of Quantum Yang-Mills' Field Theories
while the notion of symmetry, at the heart of the
Mazzola-Muzzulini's modulation is in no way linked with some kind
of "gauging" of a Lie group \cite{Weinberg-95} \cite{Weinberg-96},
resulting in the symmetry group of the modulation quantum.

\smallskip

Let us observe, at this point, that, though explicitly thought to
be the mathematical formalization of the three-stages' process
through which Sch\"{o}nberg, in the $ 9^{th} $ chapter of
\cite{Schoenberg-02}, codifying  the passage from a tonality $
t_{1} $ to a tonality $ t_{2} $ as:
\begin{enumerate}
    \item an harmonic word in a tonality $ t_{1} $
    \item an harmonic word belonging both to $ t_{1} $ and to $
    t_{2}$
    \item a $ t_{2} $-cadence determining the new tonality $ t_{2} $
\end{enumerate}
the Muzzulini-Mazzola's approach based on Symmetry Theory and
r5esultin in the definition\ref{def:Mazzola modulation} imposes
very stronger constraints.

Given two tonalities $ t_{1}, t_{2} \in {\mathcal{T}} $:
\begin{definition}
\end{definition}
MODULATIONS FROM $  t_{1} $ to $  t_{2} $:
\begin{equation}
  {\mathcal{M}} ( t_{1}, t_{2} ) \; = \; \{ \, ( p , c ) \, : \, p
  \in {\mathcal{P}}(  t_{1} \, , \,   t_{2} ) \, , \, c \in {
  \mathcal{C}}[t_{2} \, , \, {\mathcal{T}}_{n.c.} ( t_{2} ) ] \, \}
\end{equation}
\begin{definition}
\end{definition}
TONAL MUSICAL PIECES:
\begin{multline}
  {\mathcal{MP}} \; := \; \{ (hw_{1} \, m( t_{1} , t_{2} )  \,  hw_{2} \,   
\cdots \,  hw_{n-1} \, m( t_{n-1} , t_{n} )  \,  hw_{n}) \; : \\
   t_{i} \,  \in  \,  {\mathcal{T}} \, , hw_{i} \in {\mathcal{HW}}( t_{i})  
\,  m( t_{i-1} , t_{i} ) \,  \in \,  {\mathcal{M}}( t_{i-1} , t_{i} ) \, i = 
1 , \cdots , n \, \\
     n \in {\mathbb{N}} \}
\end{multline}

\begin{example}
\end{example}
ITERATING THE II-V-I PROGRESSION MOVING ALONG THE FIFTH'S CYCLE

Though disliked in Classical Music, parallel motion of perfect
fifth, and hence modulation of one step forward along the fifth's
cycles, has been reconsidered in all kinds of Modern Music
\cite{Piston-89}, \cite{Levine-95}, \cite{Strunk-00}. Such a kind
of modulations are, in some way, the simplest ones since the
tonalities connected by modulation have the highest number of
pivotal degrees, though, as we have seen in the
section\ref{ex:pivotal degrees of major tonalities at different
levels along a clockwise step in fifths' cycle}, such a number
deacreases when the level at which tonalities are considered
increases.

We have seen therein that the only pivotal degree persisting when
the level of the major tonalities is raised from 1 to 4 is the $
6^{th} $ degree of the departure-tonality coinciding with the $
2^{th} $ degree of the arrival tonality. We shall therefore adopt
such  pivotal degree. As cadence we will adopt the celebrated
cadential set (w.r.t. the natural context of the major tonalities)
II-V-I progression.

We obtain conseguentially the following a tonal musical piece
constisting in the parallel motion of the II-V-I progression along
the fifths' cycle.

Adopting the expression implemented in the Mathematica notebook of
section \ref{sec:Mathematica's notebook Mathematical Music
Toolkit}:

\begin{verbatim}
stepalongfifthscycle[n_]:=Mod[n+7,12]

fifthscycle=NestList[stepalongfifthscycle,0,12];
\end{verbatim}

the required musical piece may be immediately computed in the
following way:

\begin{verbatim}
musicalpiece[n_,l_]:= Join[harmonicwordintermofdegrees[
      tonality[majorword[Part[fifthscycle,n]],l],{2,5,1}]]

Flatten[Table[musicalpiece[n,1],{n,1,13}],1] =
{{2,5,9},{7,11,2},{0,4,7},{9,0,4},{2,6,9},{7,11,2},{4,7,11},{9,1,4},{2,6,9},{
    
11,2,6},{4,8,11},{9,1,4},{6,9,1},{11,3,6},{4,8,11},{1,4,8},{6,10,1},{11,3,
    
6},{8,11,3},{1,5,8},{6,10,1},{3,6,10},{8,0,3},{1,5,8},{10,1,5},{3,7,10},{
    8,0,3},{5,8,0},{10,2,5},{3,7,10},{0,3,7},{5,9,0},{10,2,5},{7,10,2},{0,4,
    7},{5,9,0},{2,5,9},{7,11,2},{0,4,7}}

Flatten[Table[musicalpiece[n,2],{n,1,13}],1] =
{{2,5,9,0},{7,11,2,5},{0,4,7,11},{9,0,4,7},{2,6,9,0},{7,11,2,6},{4,7,11,2},{9,
    
1,4,7},{2,6,9,1},{11,2,6,9},{4,8,11,2},{9,1,4,8},{6,9,1,4},{11,3,6,9},{4,
    8,11,3},{1,4,8,11},{6,10,1,4},{11,3,6,10},{8,11,3,6},{1,5,8,11},{6,10,1,
    
5},{3,6,10,1},{8,0,3,6},{1,5,8,0},{10,1,5,8},{3,7,10,1},{8,0,3,7},{5,8,0,
    
3},{10,2,5,8},{3,7,10,2},{0,3,7,10},{5,9,0,3},{10,2,5,9},{7,10,2,5},{0,4,
    7,10},{5,9,0,4},{2,5,9,0},{7,11,2,5},{0,4,7,11}}

Flatten[Table[musicalpiece[n,3],{n,1,13}],1] =
{{2,5,9,0,4},{7,11,2,5,9},{0,4,7,11,2},{9,0,4,7,11},{2,6,9,0,4},{7,11,2,6,9},{
    4,7,11,2,6},{9,1,4,7,11},{2,6,9,1,4},{11,2,6,9,1},{4,8,11,2,6},{9,1,4,8,
    
11},{6,9,1,4,8},{11,3,6,9,1},{4,8,11,3,6},{1,4,8,11,3},{6,10,1,4,8},{11,3,
    6,10,1},{8,11,3,6,10},{1,5,8,11,3},{6,10,1,5,8},{3,6,10,1,5},{8,0,3,6,
    
10},{1,5,8,0,3},{10,1,5,8,0},{3,7,10,1,5},{8,0,3,7,10},{5,8,0,3,7},{10,2,
    
5,8,0},{3,7,10,2,5},{0,3,7,10,2},{5,9,0,3,7},{10,2,5,9,0},{7,10,2,5,9},{0,
    4,7,10,2},{5,9,0,4,7},{2,5,9,0,4},{7,11,2,5,9},{0,4,7,11,2}}

Flatten[Table[musicalpiece[n,4],{n,1,13}],1] =
{{2,5,9,0,4,7},{7,11,2,5,9,0},{0,4,7,11,2,5},{9,0,4,7,11,2},{2,6,9,0,4,7},{7,
    
11,2,6,9,0},{4,7,11,2,6,9},{9,1,4,7,11,2},{2,6,9,1,4,7},{11,2,6,9,1,4},{4,
    
8,11,2,6,9},{9,1,4,8,11,2},{6,9,1,4,8,11},{11,3,6,9,1,4},{4,8,11,3,6,9},{
    
1,4,8,11,3,6},{6,10,1,4,8,11},{11,3,6,10,1,4},{8,11,3,6,10,1},{1,5,8,11,3,
    
6},{6,10,1,5,8,11},{3,6,10,1,5,8},{8,0,3,6,10,1},{1,5,8,0,3,6},{10,1,5,8,
    
0,3},{3,7,10,1,5,8},{8,0,3,7,10,1},{5,8,0,3,7,10},{10,2,5,8,0,3},{3,7,10,
    
2,5,8},{0,3,7,10,2,5},{5,9,0,3,7,10},{10,2,5,9,0,3},{7,10,2,5,9,0},{0,4,7,
    10,2,5},{5,9,0,4,7,10},{2,5,9,0,4,7},{7,11,2,5,9,0},{0,4,7,11,2,5}}
\end{verbatim}

\bigskip

Let us now analyze how the mathematical structure of Harmony is
modified by giving up the constraint of the incommensurability of
notes by replacing, as underlying starting scale, the
equally-tempered chromatic scale $ tempered_{eq}( C_{2} , 12 ) $
with the pytagoric scale  $ pytagoric(  C_{2} ) $.

I will here radically depart from Mazzola's way of treating both
\emph{pytagoric tuning} and \emph{just-intonation tuning} whose
philosophy could be described as hiding a problem by taking the
quotient over that: the introduction of the quotient module  $
J_{Kt} $ of the section 13.4.2.2 "Just Triadic Degree
Interpretations"  as well as that of the enharmonic projection enh
of the section24.1.3  of \cite{Mazzola-02a} are therein functional
to the concealment of a structrural problem inside an elegant
mathematical abstraction.

But, please, tell a musician he have to play an equivalence class
of notes differing by a comma and hear his answer $ \cdots ... $.

\smallskip

The net effect  of replacing the
equally-tempered chromatic scale $ tempered_{eq}( C_{2} , 12 ) $
with the pytagoric scale  $ pytagoric(  C_{2} ) $ may be
formalized as the following ansatz concerning the adopted musical
alphabet:

\begin{definition}
\end{definition}
PYTAGORIC ANSATZ:
\begin{equation}
  {\mathbb{Z}}_{12} \; \mapsto \; {\mathbb{Z}}_{12}^{\star}
\end{equation}
where:
\begin{equation}
    {\mathbb{Z}}_{12}^{\star} \; := \; \{ [i]_{12}^{(n)} \, , \, i = 1 , 
\cdots \ 12 \, , \, n \in {\mathbb{N}}  \}
\end{equation}
\begin{equation}
[i]_{12}^{(n)} \; := \; (ord((scale(\nu_{ref},3))_{12*n+1,12*n+11})_{i+1}
\end{equation}
where the note  $ [i]_{12}^{(n)} $ is the new note obtained by $
[i]_{12} $ after n turns of the \emph{fifths' cycle} and where
attention has to be kept as to the remark\ref{rem:ordering the
pytagoric scales globally or at fifth cycles' blocks}

The pytagoric ansatz may be also represented  by the following
change of  musical notation :
\begin{eqnarray}
  C^{(n)} \; &:=& \; [0]_{12}^{n} \\
  C^{\sharp \, (n)} \; &:=& \;  [1]_{12}^{n} \\
  D^{(n)} \; &:=& \; [2]_{12}^{(n)}  \\
  D^{\sharp \, (n)} \; &:=& \;  [3]_{12}^{(n)}  \\
  E^{(n)}  \; & :=& \; [4]_{12}^{(n)} \\
  F^{(n)} \; & := & \; [5]_{12}^{(n)}    \\
  F^{\sharp \, (n)} \; & := & [6]_{12}^{(n)}  \\
  G^{(n)} \; & := \; [7]_{12}^{(n)}   \\
  G^{\sharp \, (n)} \; & := \;  [8]_{12}^{(n)}  \\
  A^{(n)} \; & := \; [9]_{12}^{(n)} \\
  A^{\sharp \, (n)}  \; & := \; [10]_{12}^{(n)}  \\
  B^{(n)} \; & := \; [11]_{12}^{(n)}  \\
\end{eqnarray}

One has by construction that:
\begin{theorem}
\end{theorem}
\begin{equation}
   {\mathbf{p}} ( [i]_{12}^{(n+1)} ) \; = \; {\mathbf{p}} ( [i]_{12}^{(n)}
   ) \, + \,  {\mathbf{Kf}} \; \; \forall n \in {\mathbb{N}} \, ,
   \, i = 0 , \cdots , 11
\end{equation}
Introduced the  notation $ \Xi := {\mathbb{Z}}_{12}^{\star} $:
\begin{definition}
\end{definition}
FIFTH CYCLES' ORDERING RELATION OVER $ \Xi $:
\begin{equation}
    [n_{1}]_{12}^{(m)} \, <_{fc} \, [n_{2}]_{12}^{(m)}  \; \Leftrightarrow
    \; n_{1} \, < \,  n_{2}
\end{equation}
\begin{equation}
    [n_{1}]_{12}^{(m_{1})} \, <_{fc} \, [n_{2}]_{12}^{(m_{2})} \; 
\Leftrightarrow
   \; m_{1} \, < \,  m_{2}
\end{equation}
One has clearly that:
\begin{eqnarray}
   card ( \Xi) \, = \, card (  {\mathbb{Z}}_{12}^{\star}  ) \; = \;
   \aleph_{0} \\
   card ( \Xi^{\star} ) \; = \; card( ({\mathbb{Z}}_{12}^{\star}) ^{\star}  
) \; = \; \aleph_{1}
\end{eqnarray}

Given the  pytagoric alphabet $ \Xi $, a word $ \tilde{x} \in
\Xi^{\star} $ and a map $ g \, : \, \Xi \rightarrow \Xi $:
\begin{definition} \label{def:map induced on pytagoric words by
a map on the pytagoric alphabet}
\end{definition}
MAP INDUCED BY g ON WORDS:

the map $ \hat{g} \, : \, \Xi^{\star} \, \rightarrow \,
\Xi^{\star} $ :
\begin{equation}
    \hat{g} ( \tilde{x} ) \; = \; \cdot_{i=1}^{| \tilde{x} |} g
    (x_{i})
\end{equation}

Given a letter $ z \, := \, [i]_{12}^{(n)} \, \in \, \Xi $

\begin{definition}
\end{definition}
PYTAGORIC TRANSLATION BY z :

the map $ T_{z} \,  :  \, \Xi \, \mapsto \, \Xi $:
\begin{equation}
    T_{z} ( [j]_{12}^{(m)} ) \; = \; [i +_{12} j]^{(m)} \; \; \forall m
    \in {\mathbb{N}} \, , \,  \forall j \in \{ 0,1, \cdots, 11 \}
\end{equation}

\smallskip

\begin{definition}
\end{definition}
PYTAGORIC INVERSION:

the map $ Inv^{Pyt}  \,  :  \, \Xi \, \mapsto \, \Xi $:
\begin{equation}
    Inv ( [j]_{12}^{(m)} ) \; = \; [ 0 \, -_{12} \, j ]_{12}^{(m)}
    \; \; \forall m
    \in {\mathbb{N}} \, , \,  \forall j \in \{ 0,1, \cdots, 11 \}
\end{equation}

\smallskip

\begin{definition} \label{def:operator of fifth cycle's raising}
\end{definition}
OPERATOR OF FIFTH CYCLE'S RAISING:

the map $ C_{+} \,  :  \, \Xi \, \mapsto \, \Xi $:
\begin{equation}
    C_{+} ( [j]_{12}^{(m)} ) \; := \; [j]_{12}^{(m+1)}   \; \; \forall m
    \in {\mathbb{N}} \, , \,  \forall j \in \{ 0,1, \cdots, 11 \}
\end{equation}

\smallskip

\begin{definition} \label{def:operator of fifth cycle's lowering}
\end{definition}
OPERATOR OF FIFTH CYCLE'S LOWERING:

the map $ C_{-} \,  :  \, \Xi \, \mapsto \, \Xi $:
\begin{equation}
    C_{-} ( [j]_{12}^{(m)} ) \; := \; [j]_{12}^{(m-1)}   \; \; \forall m
    \in {\mathbb{N}} \, , \,  \forall j \in \{ 0,1, \cdots, 11 \}
\end{equation}

Once again the effect of the pytagoric ansatz may be concretely
analyzed adopting  the following Mathematica expressions from the
notebook of section \ref{sec:Mathematica's notebook Mathematical
Music Toolkit}:

\begin{verbatim}

<<DiscreteMath`Combinatorica`;

<<DiscreteMath`Permutations`

putinorder[listofnumbers_]:=Permute[listofnumbers,Ordering[listofnumbers]]

pytagoricletter[n_,m_]:={letter[n],m}

pytagoricalphabetuptofifthcycles[m_]:=
  Flatten[Table[pytagoricletter[n,j],{j,0,m},{n,0,11}],1]

FROMpytagoricletterTOnote[x_]:=
  Part[putinorder[
      Take[scaleatfixedinterval[referencenote,3,12*Part[x,2]+11],-12]],
    Part[x,1]+1]

pytagoricscaleuptofifthcycles[m_]:=
  Map[FROMpytagoricletterTOnote,pytagoricalphabetuptofifthcycles[m]]

pytagoricword[w_,n_]:=Table[pytagoricletter[Part[w,i],n],{i,1,Length[w]}]

\end{verbatim}
\begin{example}
\end{example}
THE PART OF THE PYTAGORIC MUSICAL ALPHABET GENERATED BY THE FIRTS
10 FIFTHS' CYCLES


(where we have shown to take such an high number of decimals to
remark that only rational numbers are involved) that, graphically
represented, appear as:

\begin{verbatim}
ListPlot[N[pytagoricscaleuptofifthcycles[10],225]] =


\end{verbatim}

\smallskip

\includegraphics{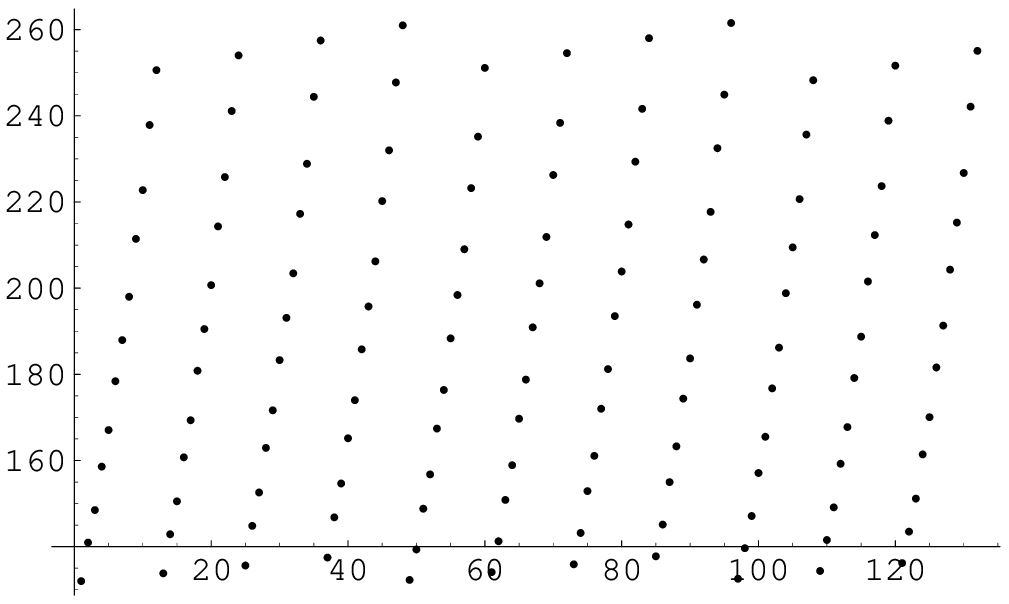}

\smallskip
revealing, also in a pictorial way:

\begin{verbatim}

ListPlot[Table[
    
FROMpytagoricletterTOnote[Part[pytagoricalphabetuptofifthcycles[100],i]],{
      i,1,Length[pytagoricalphabetuptofifthcycles[100]]}],
  PlotJoined\[Rule]True] =

\end{verbatim}

\smallskip

\includegraphics{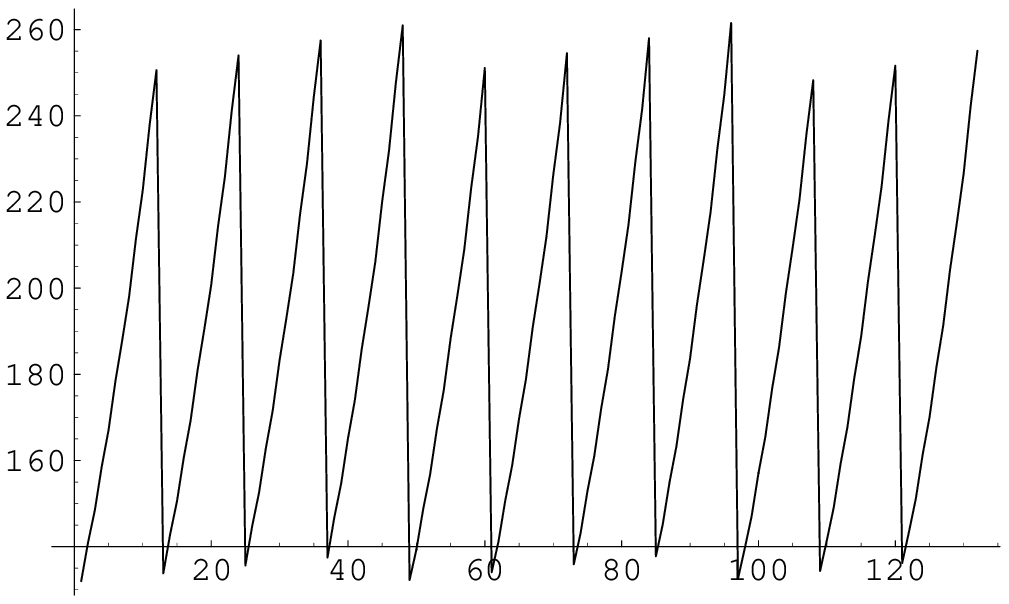}

\smallskip

\begin{verbatim}

ListPlot3D[Table[FROMpytagoricletterTOnote[{i,n}],{i,0,11},{n,0,10}]]
=
\end{verbatim}

\smallskip

\includegraphics{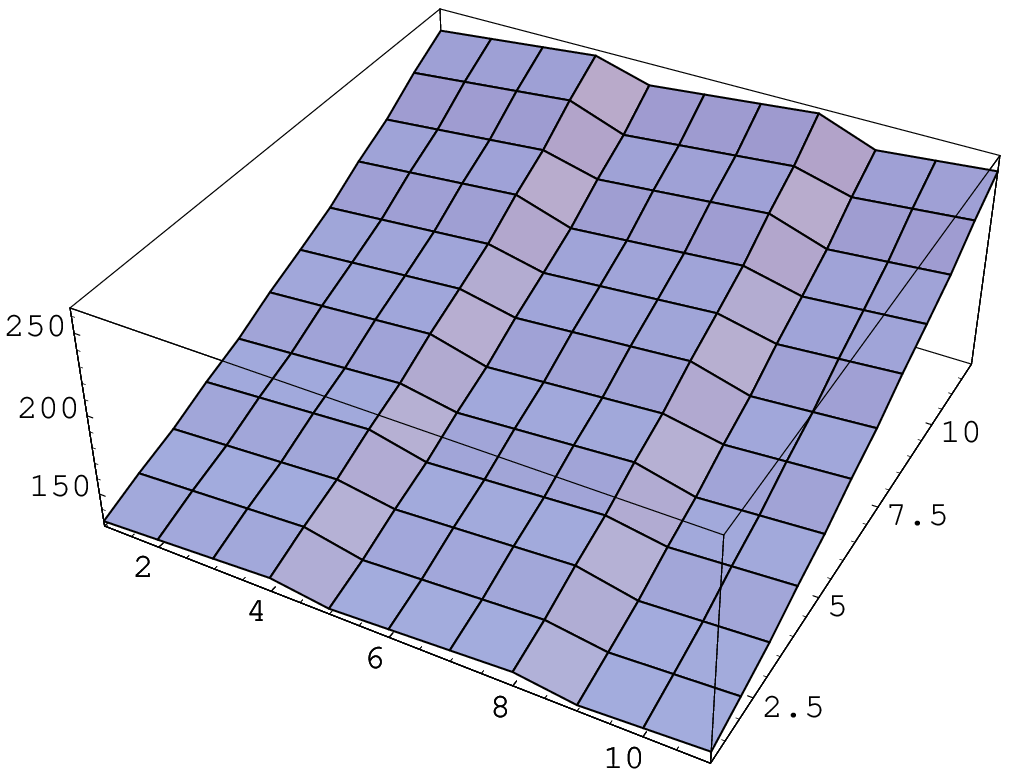}

\smallskip

the obvious fact that the ordering relation $ <_{fc} $ over $
\Xi^{\star}_{NR} $ doesn't correspond to the natural ordering
relation $ < $ among numbers

\medskip

It is then
interesting to see how all the previous formalization of Harmony
is affected by the \emph{pytagoric ansatz};

given a word $ \vec{x} = (x_{1} \, , \, \cdots \, , x_{| \vec{x}
|} ) \in {\mathbb{Z}}_{12}^{\star} $ and a natural number $ n \in
{\mathbb{N}} $:
\begin{definition}
\end{definition}
PYTAGORIC WORD OF $ \vec{x} $ AT CYCLE n:
\begin{equation}
\tilde{ pw} ( \vec{x} , n ) \; := (x_{1}^{(n)}  \, , \, \cdots \, , x_{| 
\vec{x}
|}^{(n)} )
\end{equation}

Given a  pytagoric word $ \tilde{ x} \, \in \,   \Xi^{\star} \, -
\, \bigcup_{k=1}^{4}  \Xi^{k}  $ and an integer $ i \, \in \{ 1,
\, \cdots \, , | \vec{x} | \} $:
\begin{definition} \label{def:modes of a  pytagoric word}
\end{definition}
MODE OF $ \tilde{x} $ of $ i^{th}$ DEGREE :
\begin{equation}
    mode( \tilde{x} \, , \, i ) \; := \; {\mathcal{S}}_{cycl}^{i}
( \tilde{x} )
\end{equation}
where $ {\mathcal{S}}_{cycl} $ denotes again the operator of
cyclic shift.

Given furthermore an integer $ n \in { \mathbb{N}} $
\begin{definition} \label{def:chord of a pytagoric word at a certain level}
\end{definition}
CHORD OF  $ \tilde{x} $ of $ i^{th}$ DEGREE AT LEVEL n :
\begin{equation}
    chord ( \tilde{x} \, , \, i \, , \, n ) \; := \cdot_{j=0}^{n+3} mode( 
\tilde{x} \, , \, i )_{2j+1}
\end{equation}
where $ \cdot $ denotes the concatenation operator while $  x_{j}
$ denotes the $ j^{th} $ letter of the pytagoric word $ \tilde{x}
$.

Exactly as without the pytagoric ansatz, there exists a certain
maximum level
$ maxlevel( \tilde{x} ) $ such that $ chord ( \tilde{x} \, , \, i \, , \, n
) \, $ for $ n > maxlevel( \tilde{x} ) $ simply adds notes already
contained in the chord, again given by:
\begin{equation} \label{eq:maximum level of a  pytagoric word for chords}
    maxlevel ( \tilde{x} ) \; = \; 1 \, + \, Int(\frac{ |\tilde{x} | -5}{2}) 
\, -
\, \frac{ (-1)^{ |\tilde{x} |} - 1 }{2} \dot Int( \frac{
|\tilde{x} |}{2})
\end{equation}

\medskip

Let us now analyze how the previously discussed mathematical
formalization of Harmony is affected by the pytagoric ansatz.

Given a pytagoric word $ \tilde{x} \,  \in \,  \Xi^{\star}_{NR} -
\bigcup_{k=1}^{4} \Xi^{k}_{NR} $ and a number $ n \in
{\mathbb{N}}_{+} \, : \,  n < maxlevel ( \vec{x} )  $:
\begin{definition} \label{def:pytagoric tonality of a non-repetitive word at 
a certain level}
\end{definition}
PYTAGORIC TONALITY OF THE PYTAGORIC WORD $ \tilde{x} \, \in \,
\Xi^{\star}_{NR} \, - \, \bigcup_{k=1}^{4} \Xi^{k}_{NR} $ AT THE
LEVEL n (n-TONALITY OF $ \vec{x} $) :

$ pytagorictonality[\tilde{x} \, , \, n ] \; := \;  ( \tilde{x} \, , \, 
chord^{(n)} ) $, where $ chord^{(n)} \, :
\,  \{ 1  , \cdots , |\tilde{x} | \}  \, \mapsto \, \Xi^{\star} $
is the map such that:
\begin{equation}
    chord^{(n)} (i)  \; = \; chord ( \tilde{x} \, , \, i \, , \, n )
\end{equation}

I will denote the set of all the pytagoric n-tonalities by
${\mathcal{T}}_{n}^{Pyt} $ while I will denote by $ {\mathcal{T}
}^{Pyt} $ the set of all the pytagoric tonalities at any level.

One has clearly that:
\begin{equation}
    card( {\mathcal{T}}_{n}^{Pyt} ) \; = card ( \Xi^{\star}_{NR} \, - \,
\bigcup_{k=1}^{4} \Xi^{k}_{NR} ) \; = \; \aleph_{1}
\end{equation}
\begin{equation}
    card( {\mathcal{T}} ) \; = \; card ( \bigcup_{n \in {\mathbb{N}}}
    {\mathcal{T}}_{n}^{Pyt} ) \; = \; \aleph_{1}
\end{equation}
Given a pytagoric n-tonality $ t \in {\mathcal{T}}_{n}^{Pyt} $:
\begin{definition}
\end{definition}
PYTAGORIC HARMONIC WORDS OF t:
\begin{equation}
    {\mathcal{HW}}^{Pyt} (t) \; := \; (Range( chord^{(n)} (t))^{\star}
\end{equation}
Given two pytagoric tonalities $ t_{1} \, , \,  t_{2} \; \in \;
{\mathcal{T}}^{Pyt}  $:
\begin{definition}
\end{definition}
PYTAGORIC PIVOTAL DEGREES OF $ t_{1} $ AND $ t_{2} $:
\begin{equation}
    {\mathcal{P}}^{Pyt}(  t_{1} \, , \,   t_{2} ) \; := \; Range(
    chord^{(level(t_{1})} ( t_{1} )) \, \bigcap \,  Range(
    chord^{(level(t_{2}))} ( t_{2} ))
\end{equation}

Given a pytagoric tonality $ t \in {\mathcal{T}}^{Pyt}$ and a
context $ {\mathcal{T}}^{Pyt}_{context} \subseteq
{\mathcal{T}}^{Pyt} $
\begin{definition} \label{pytagoric cadences of a tonality w.r.t. to a 
context}
\end{definition}
PYTAGORIC CADENCES OF $  t $ W.R.T. THE CONTEXT $
{\mathcal{T}}_{context} $
\begin{equation}
     { \mathcal{C}}^{Pyt}(t \, , \, {\mathcal{T}}_{context} ) \; := \; \{ c  
\, \in \, {\mathcal{HW}}^{Pyt} (t) \; :  \\
      ( \,  c  \, \in \, {\mathcal{HW}^{Pyt}} (t) \; \Rightarrow \;  t \, = 
\,  u \,
      ) \, \, \forall u \in {\mathcal{T}}_{context} \}
\end{equation}

Given a pytagoric word $ \tilde{x} \in \Xi^{\star} $ and a
natural number $ n \in {\mathbb{N}} \, : \, n < maxlevel( \vec{x}
) $
\begin{definition} \label{def:pytagoric natural context of a pytagoric word 
at a certain level}
\end{definition}
    PYTAGORIC NATURAL CONTEXT OF $ pytagorictonality[ \tilde{x} , n ] $
\begin{equation}
    {\mathcal{T}}^{Pyt}_{n.c.} ( \tilde{x} , n ) \; := \; \{
    pytagorictonality(\tilde{y} \, , \, n) \: : \: \tilde{x} \, 
\sim_{T}^{Pyt}
    \tilde{y} \}
\end{equation}
where $ \sim_{T}^{Pyt} $ is the pytagorical translational
equivalence as induced by the map $ \hat{T}^{Pyt} $ obtained
applying  the definition\ref{def:map induced on pytagoric words by
a map on the pytagoric alphabet} to the alphabet $ \Xi $.

\begin{definition}
\end{definition}
PYTAGORIC MINIMAL CADENCES OF $  t $ W.R.T. THE CONTEXT $
{\mathcal{T}}^{Pyt}_{context} $
\begin{equation}
    {\mathcal{MC}} ( t \, , \, {\mathcal{T}}^{Pyt}_{context} ) \; := \; \{ 
c_{1}  \, \in \, {
\mathcal{C}}(t, {\mathcal{T}}_{context}) \, : \,  \nexists c_{2}
\in { \mathcal{C}}(t,{\mathcal{T}}^{Pyt}_{context}) , c_{2} \,
<_{p} \, c_{1} \}
\end{equation}

Given two pytagoric tonalities $ t_{1}, t_{2} \in
{\mathcal{T}}^{Pyt} $:
\begin{definition} \label{def:pytagoric modulations}
\end{definition}
PYTAGORIC MODULATIONS FROM $  t_{1} $ to $  t_{2} $:
\begin{equation}
  {\mathcal{M}}^{Pyt} ( t_{1}, t_{2} ) \; = \; \{ \, ( p , c ) \, : \, p
  \in {\mathcal{P}}^{Pyt}(  t_{1} \, , \,   t_{2} ) \, , \, c \in {
  \mathcal{C}}^{Pyt}(t_{2}) \, \}
\end{equation}
\begin{definition} \label{def:pytagoric tonal musical pieces}
\end{definition}
PYTAGORIC TONAL MUSICAL PIECES:
\begin{multline}
  {\mathcal{MP}}^{Pyt} \; := \; \{ (hw_{1} \, m( t_{1} , t_{2} )  \,  hw_{2} 
\,   \cdots \,  hw_{n-1} \, m( t_{n-1} , t_{n})   \,  hw_{n} ) \; : \\
   t_{i} \,  \in  \,  {\mathcal{T}}^{Pyt} \, , hw_{i} \in 
{\mathcal{HW}}^{Pyt}( t_{i})  \,  m( t_{i-1} , t_{i} ) \,  \in \,  
{\mathcal{M}}^{Pyt}( t_{i-1} , t_{i} ) \, i = 1 , \cdots , n \, \\
     n \in {\mathbb{N}} \}
\end{multline}

\begin{remark}
\end{remark}
PYTAGORICALLY-TUNED AND JUST-TUNED EULER NOTES VERSUS PYTAGORIC
AND JUST SCALES

It is important, at this point to analyze the relation existing
among the \emph{pytagorical Euler notes} of
definition\ref{def:pytagorically-tuned Euler notes} and the
\emph{pytagoric scale} $ pytagoric ( \nu_{ref} , 3 ) $ as well as
the relation existing among the \emph{just-tuned Euler notes} of
definition\ref{def:just-tuned Euler notes} and the celebrated
\emph{just-intonation scale} \footnote{often called the
\emph{simple-ratios' scale} or the \emph{Zarlinian scale}, this
last name being, anyway, based to an erroneous hystorical
attribution to Gioseffo Zarlino (cfr. the voice "Scale Musicali
(le)  antiche e moderne" of \cite{Righini-80}} we are going to
introduce.

Looking at the pitches obtained by pytagoric Euler points lying in
a cube of side = 10 having as center the origin (i.e. the Euler
point of $ \nu_{ref} $:

\begin{verbatim}
N[Map[FROMeulerpointTOpitch,Flatten[Table[x*octaveepoint+y*fifthpoint,

{x,-5,5},{y,-5,5}],1]]] =

{-15509.8,-13607.8,-11705.9,-9803.91,-7901.96,-6000.,-4098.04,-2196.09,-294.\
135,1607.82,
  
3509.78,-14309.8,-12407.8,-10505.9,-8603.91,-6701.96,-4800.,-2898.04,-996.\
09,905.865,2807.82,
  
4709.78,-13109.8,-11207.8,-9305.87,-7403.91,-5501.96,-3600.,-1698.04,203.91,
  2105.87,4007.82,
  5909.78,-11909.8,-10007.8,-8105.87,-6203.91,-4301.96,-2400.,-498.045,
  1403.91,3305.87,5207.82,
  
7109.78,-10709.8,-8807.82,-6905.87,-5003.91,-3101.96,-1200.,701.955,2603.91,
  4505.87,6407.82,8309.78,-9509.78,-7607.82,-5705.87,-3803.91,-1901.96,0,
  1901.96,3803.91,5705.87,7607.82,
  
9509.78,-8309.78,-6407.82,-4505.87,-2603.91,-701.955,1200.,3101.96,5003.91,
  6905.87,8807.82,10709.8,-7109.78,-5207.82,-3305.87,-1403.91,498.045,2400.,
  
4301.96,6203.91,8105.87,10007.8,11909.8,-5909.78,-4007.82,-2105.87,-203.91,
  1698.04,3600.,5501.96,7403.91,9305.87,11207.8,
  13109.8,-4709.78,-2807.82,-905.865,996.09,2898.04,4800.,6701.96,8603.91,
  10505.9,12407.8,14309.8,-3509.78,-1607.82,294.135,2196.09,4098.04,6000.,
  7901.96,9803.91,11705.9,13607.8,15509.8}

\end{verbatim}

and comparing it with  the first 100 pitchs of the pytagoric
scale:

\begin{verbatim}
N[Map[FROMnoteTOpitch,scaleatfixedinterval[132,3,100]]] =

{-15509.8,-13607.8,-11705.9,-9803.91,-7901.96,-6000.,-4098.04,-2196.09,-294.\
135,1607.82,
  
3509.78,-14309.8,-12407.8,-10505.9,-8603.91,-6701.96,-4800.,-2898.04,-996.\
09,905.865,2807.82,
  
4709.78,-13109.8,-11207.8,-9305.87,-7403.91,-5501.96,-3600.,-1698.04,203.91,
  2105.87,4007.82,
  5909.78,-11909.8,-10007.8,-8105.87,-6203.91,-4301.96,-2400.,-498.045,
  1403.91,3305.87,5207.82,
  
7109.78,-10709.8,-8807.82,-6905.87,-5003.91,-3101.96,-1200.,701.955,2603.91,
  4505.87,6407.82,8309.78,-9509.78,-7607.82,-5705.87,-3803.91,-1901.96,0,
  1901.96,3803.91,5705.87,7607.82,
  
9509.78,-8309.78,-6407.82,-4505.87,-2603.91,-701.955,1200.,3101.96,5003.91,
  6905.87,8807.82,10709.8,-7109.78,-5207.82,-3305.87,-1403.91,498.045,2400.,
  
4301.96,6203.91,8105.87,10007.8,11909.8,-5909.78,-4007.82,-2105.87,-203.91,
  1698.04,3600.,5501.96,7403.91,9305.87,11207.8,
  13109.8,-4709.78,-2807.82,-905.865,996.09,2898.04,4800.,6701.96,8603.91,
  10505.9,12407.8,14309.8,-3509.78,-1607.82,294.135,2196.09,4098.04,6000.,
  7901.96,9803.91,11705.9,13607.8,15509.8}


\end{verbatim}
one is euristically led to discover that the \emph{pytagoric Euler
notes} and the \emph{notes of the pytagoric scales} are different
concepts:
\begin{verbatim}
  Intersection[
  N[Map[FROMeulerpointTOpitch,Flatten[Table[x*octaveepoint+y*fifthpoint
  ,{x,-5,5},{y,-5,5}],1]]],
N[Map[FROMnoteTOpitch,scaleatfixedinterval[132,3,100]]] =

  {0}
\end{verbatim}

\smallskip

The situation is subtler as to \emph{just intonation scales}:
\begin{definition} \label{def:diatonic just intonation scale of C major}
\end{definition}
DIATONIC JUST INTONATION SCALE OF C MAJOR:

the scale $ \{ \nu_{1} \, , \, \cdots \, , \, \nu_{7} \} $
specified by the following table:

\begin{tabular}{|c|c|c|}
  \hline
  $\nu \, Hz $ & $ {\mathbf{mi}}(\nu  , \nu_{ref}) $  & $ {\mathbf{p}}(\nu) 
\, Cents $ \\
  \hline
   132 & 1 & 0 \\
  148.5 & $ \frac{9}{8} $ & 203.91 \\
  165 & $ \frac{5}{4} $ & 386.314 \\
  176 & $ \frac{4}{3} $ & 498.045 \\
  198 & $ \frac{3}{2} $ & 701.955 \\
  220 & $ \frac{5}{3 }$ & 719.354 \\
  247.5 & $ \frac{15}{8} $ & 1088.27 \\ \hline
\end{tabular}

\smallskip

Behind the introduction of the scale of
definition\ref{def:diatonic just intonation scale of C major} laid
a great musical innovation (that as many ones had to face the
opposition of the Catholic Church)  , namely the acceptance of the
third-major interval as a consonant one.

Tonal Harmony was codified assuming the axiom\ref{ax:axiom of the
naturality of esthetics} such an innovation assumed as cornerstone
the answer\ref{an:the naife answer} to the question
\ref{qu:question on the more natural musical messages}.

Since:
\begin{itemize}
    \item  $ harmonic(  \nu_{{ref}} \, , \, 1 ) \, = \, 2 \nu_{ref} $
had been already used as a basis for the assumption of the
axiom\ref{ax:axiom of perception of repetition for the first
harmonic} with the consequential constraint of the rescaling to
$scale-range ( \nu_{ref} )$
    \item $ harmonic( \nu_{ref} \, , \, 2 ) \,
= \, 3 \nu_{ref} $ had been used in the construction of the
pytagoric scale to fix the perfect fifth interval to the value $
{\mathbf{mi}} ( \nu_{perfect \,  fifth} \: , \: \nu_{ref} ) \, =
\, \frac{3}{2} $
    \item $ harmonic( \nu_{ref} \, , \, 3 ) \, =
\, 4 \nu_{ref}  $  couldn't add anything new w.r.t. the
consequences induced by $ harmonic(  \nu_{{ref}} \, , \, 1 ) \, =
\, 2 \nu_{ref}  $
\end{itemize}
it was clear that the next protagonist had to be $   harmonic(
\nu_{ref} \, , \, 2 ) \, = \, 5 \nu_{ref} $.

One could, at this point, thought to introduce $ scale( \nu_{ref}
\, 5 ) $ but this would have signified to throw away the
consequences induced by  $ harmonic( \nu_{ref} \, , \, 3 ) $.

\smallskip

Behind the introduction of the \emph{just intonation diatonic
scales} laid an other idea, namely the pytagorical creed  that the
musical consonant intervals correspond to simple ratios of integer
numbers.

Given a rational number $ r \,  = \,  \frac{n}{m} \, n , m \in
{\mathbb{N}} \, , \, m \neq 0 $ measures of its "simplicity" such
as the following:
\begin{definition}
\end{definition}
EMPIRICAL MEASURE OF SIMPLICITY OF $ r$:
\begin{equation}
    {\mathcal{ESM}} (r ) \; := \; \frac{1}{ gcd(n,m) } \frac{m+n}{m \cdot n}
\end{equation}
(whose sometimes claimed link with physical consonance we confuted
in the example\ref{ex:physical consonance of notes w.r.t. an ideal
instrument}) gave to the pytagorical major third intervals too low
values:
\begin{multline}
    {\mathcal{ESM}}  ( {\mathbf{mi} } ( E^{(0)} \, , \,
    C^{(0)} )) \; = \; {\mathcal{ESM}}  ( {\mathbf{mi} } ( A^{(0)} \, , \,
    F^{(0)} )) \; \\
\; {\mathcal{ESM}}  ( {\mathbf{mi} } ( B^{(0)} \, , \,
    G^{(0)} )) \; = \;  {\mathcal{ESM}} ( \frac{81}{64} ) \;
    = \frac{145}{5184} \; = \; 0.0279707
\end{multline}
specially if compared with:
\begin{equation}
     {\mathcal{ESM}} ( \frac{5}{4} ) \; = \; \frac{9}{20} \; = \;
     0.45
\end{equation}
The definition\ref{def:diatonic just intonation scale of C major}
was assumed as a sort of compromize between the necessity of
implementing the new ratio $ \frac{5}{4} $ for the third major's
interval trying to conserve as much as possible the ratio $
\frac{3}{2} $ for the perfect fifth's interval, a compromise
lacking of the mathematical elegance of the pytagoric scale.

So, while preserving the ratio $ \frac{3}{2} $ for the $ C - G $,
the $ E - B $ , $ F-C' $ perfect fifths' intervals it gives a
ratio of $ \frac{40}{27} $ for the $ D - A $ fifth interval for
which:
\begin{equation}
    {\mathcal{ESM}} ( \frac{40}{27} ) \; = \; \frac{67}{1080} \; =
    \; 0.062037
\end{equation}

\smallskip

Given the empirical nature of such a compromise, it doesn't
astonish that the attempt of completing the diatonic
just-intonation scale adding the 5 lacking notes necessary to
obtain a putative chromatic just intonation scale following the
same principle of looking for simple ratios, was performed by
different authors in many ways.

Let us observe, at this point, that:
\begin{theorem} \label{th:the diatonic major scale is made of just-Euler 
notes}
\end{theorem}
THE DIATONIC JUST-INTONATION MAJOR SCALE IS MADE OF JUST-EULER
NOTES:
\begin{eqnarray}
  {\mathcal{P}}_{Euler} ( \nu_{ref} ) \; &=& \; ( \, 0 \,  , \,  0 \,  , \,  
0 \, ) \\
  {\mathcal{P}}_{Euler} ( \frac{9}{8} \nu_{ref} ) \; &=& \; ( \, -3 \,  , \, 
  2 \,  , \,  0 \,  ) \\
  {\mathcal{P}}_{Euler} ( \frac{5}{4} \nu_{ref} ) \; &=& \; ( \, -2 \,  , \, 
  0 \,  , \, 1  \,) \\
  {\mathcal{P}}_{Euler} ( \frac{4}{3} \nu_{ref} ) \; &=& \; ( \, 2 \,  , \,  
-1 \,  , \,  0 \, ) \\
  {\mathcal{P}}_{Euler} ( \frac{3}{2 }\nu_{ref} ) \; &=& \; ( \, -1 \,  , \, 
  1 \,  , \,  0 \, ) \\
  {\mathcal{P}}_{Euler} ( \frac{5}{3} \nu_{ref} ) \; &=& \; ( \, 0 \,  , \,  
-1 \,  , \,  1 \, ) \\
  {\mathcal{P}}_{Euler} ( \frac{15}{8} \nu_{ref} ) \; &=& \; ( \, -3 \,  , 
\,  1 \,  , \,  1 \, ) \\
\end{eqnarray}

Theorem\ref{th:the diatonic major scale is made of just-Euler
notes} suggests to choice the other 5 notes among the
simple-ratio's just Euler notes; one obtains in this way  the
following (cfr. the appendix K "Just and Well-Tempered Tuning" of
\cite{Mazzola-02a}):
\begin{definition}
\end{definition}
VOGEL'S CHROMATIC JUST-INTONATION SCALE OF C MAJOR

\begin{eqnarray}
  C^{(0)}_{just} \; & = & \;  \nu_{ref} \\
  C^{\sharp \, (0)}_{just} \; &  = & \; \frac{16}{15} \nu_{ref}  \\
  D^{(0)}_{just} \; &=& \; \frac{9}{8} \nu_{ref}   \\
  D^{\sharp \, (0)}_{just} \; & = & \; \frac{6}{5} \nu_{ref}   \\
  E^{(0)}_{just} \; & = & \; \frac{5}{4} \nu_{ref}  \\
  F^{(0)}_{just} \; & = & \; \frac{4}{3} \nu_{ref}   \\
  F^{\sharp \, (0)}_{just} \; &  = & \; \frac{45}{32} \nu_{ref} \\
  G^{(0)}_{just} \; &=& \; \frac{3}{2} \nu_{ref}   \\
  G^{\sharp \, (0)}_{just} \; &=& \; \frac{8}{5}  \nu_{ref}   \\
  A^{(0)}_{just} \; &=& \;  \frac{5}{3} \nu_{ref}  \\
  A^{\sharp \, (0)}_{just} \; & = & \; \frac{16}{9} \nu_{ref}  \\
  B ^{(0)}_{just} \; &=& \; \frac{15}{8} \nu_{ref} \\
\end{eqnarray}
where:
\begin{eqnarray*}
  {\mathcal{P}}_{Euler} ( C^{\sharp \, (0)}_{just} ) \; &=& \; ( \, 4 \, , 
\, -1 \, , \, -1 \, )  \\
  {\mathcal{P}}_{Euler} ( D^{\sharp \, (0)}_{just} ) \; &=& \; ( \, 1 \, , 
\, 1 \, , \, -1 \, ) \\
  {\mathcal{P}}_{Euler} ( F^{\sharp \, (0)}_{just} ) \; &=& \; ( \, -5 \, , 
\, 2 \, , \, 1 \, )  \\
  {\mathcal{P}}_{Euler} ( G^{\sharp \, (0)}_{just} ) \; &=& \; ( \, 3 \, , 
\, 0 \, , \, -1 \, ) \\
  {\mathcal{P}}_{Euler} ( A^{\sharp \, (0)}_{just} ) \; &=& \; ( \, 4 \, , 
\, -2 \, , \, 0 \, )
\end{eqnarray*}

\smallskip

Visualized in Euler's space Vogel's chromatic just-intonation
scale of C major appears as:

\begin{verbatim}
<<Graphics`Graphics3D`

vogelchromaticjustlistofeulerpoints:={{0,0,0},{4,-1,-1},{-3,2,0},{1,1,-1},
{-2,0,1},{2,-1,0},{-5,2,1},{-1,1,0},{3,0,-1},{0,-1,1},{4,-2,0}{-3,1,1}}

ScatterPlot3D[vogelchromaticjustlistofeulerpoints]
\end{verbatim}

\smallskip

\includegraphics{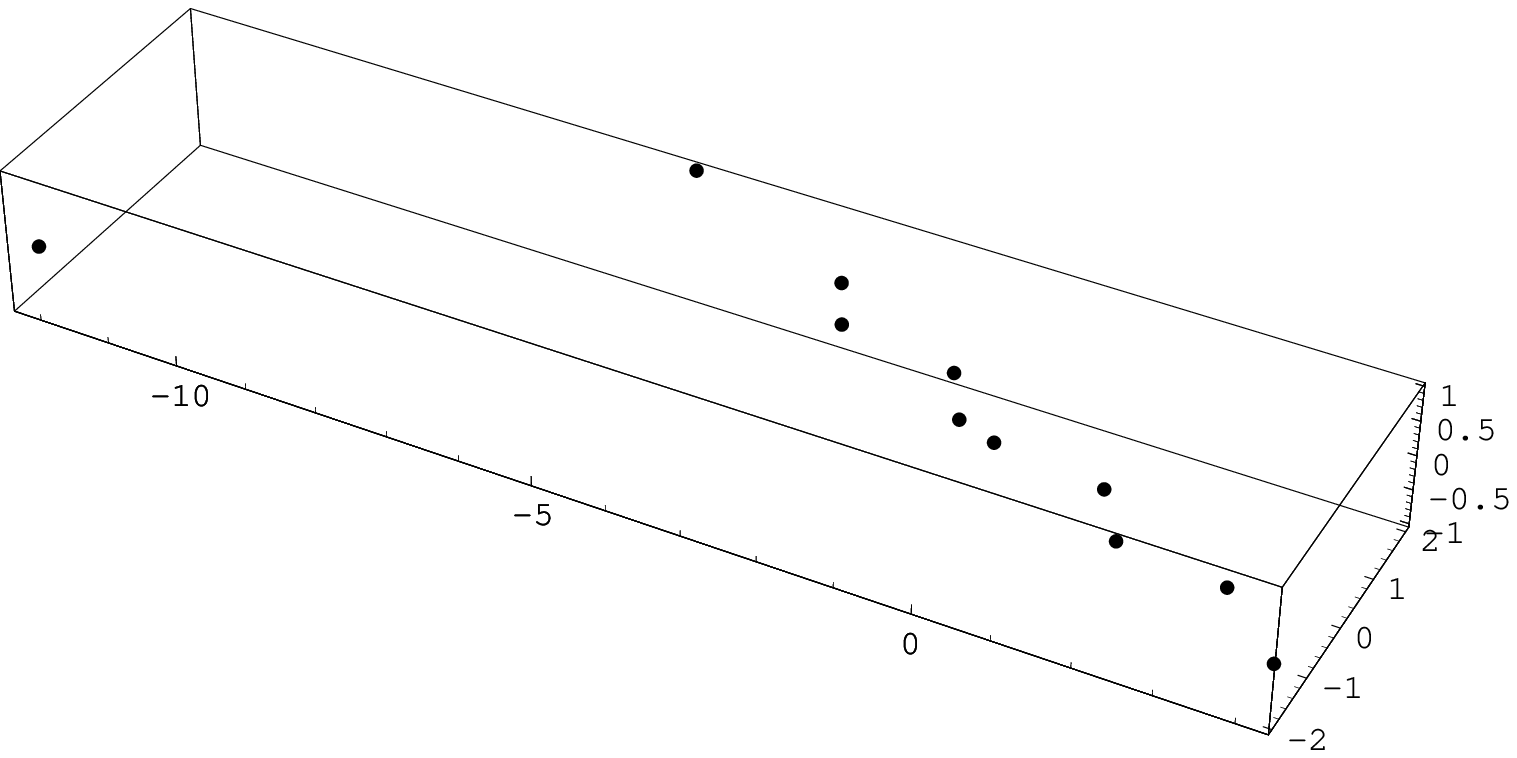}

\bigskip

Let us, now, observe that whichever choice one effects as to the 5
chromatic notes one adopts to complete the diatonic just
intonation scale of C major one have that, owing the the
\emph{third-comma}, fifth cycles don't close so that, a situation
anbalogous to the one seen as to pytagoric scale occurs.

It follows that once again the passage from the
12-equally-tempered scale to the just-intonation scale may be
formalized through the following:

\begin{definition}
\end{definition}
JUST INTONATION ANSATZ:
\begin{equation}
  {\mathbb{Z}}_{12} \; \mapsto \; {\mathbb{Z}}_{12}^{\star \, (just)}
\end{equation}
where:
\begin{equation}
    {\mathbb{Z}}_{12}^{\star \, (just)} \; := \; \{ [i]_{12}^{(n) \, (just} 
\, , \, i = 1 , \cdots \ 12 \, , \, n \in {\mathbb{N}}  \}
\end{equation}
\begin{equation}
[i]_{12}^{(n)} \; := \; (ord((scale(\nu_{ref},3))_{12*n+1,12*n+11})_{i+1}
\end{equation}
where the note  $ [i]_{12}^{(n) \, (just)} $ is the new note
obtained by $ [i]_{12} $ after n turns of the \emph{fifths'
cycle}, defined by the conditions:
\begin{itemize}
    \item
\begin{equation}
     {\mathbf{p}} ( [i]_{12}^{(n+1) \, (just)} ) \; := \; {\mathbf{p}} ( 
[i]_{12}^{(n) \, (just)}
   ) \, + \,  {\mathbf{Kt}} \; \; \forall n \in {\mathbb{N}} \, ,
   \, i = 0 , \cdots , 11
\end{equation}
    \item
\begin{equation}
  {\mathcal{P}}_{Euler} ( [i]_{12}^{(n) \, (just) } ) \; \in \;  
{\mathcal{N}}_{Euler}^{just-tuned}
\end{equation}

\end{itemize}

\smallskip

The \emph{just intonation ansatz} may be also represented  by the
following change of musical notation :
\begin{eqnarray}
  C^{(n) \, just} \; &:=& \; [0]_{12}^{(n) \, just} \\
  C^{\sharp \, (n) \, just} \; &:=& \;  [1]_{12}^{(n) \, just} \\
  D^{(n) \, just} \; &:=& \; [2]_{12}^{(n) \, just}  \\
  D^{\sharp \, (n)  \, just } \; &:=& \;  [3]_{12}^{(n)  \, just}  \\
  E^{(n)  \, just}  \; & :=& \; [4]_{12}^{(n)  \, just} \\
  F^{(n)  \, just} \; & := & \; [5]_{12}^{(n)  \, just}    \\
  F^{\sharp \, (n)  \, just} \; & := & [6]_{12}^{(n)  \, just}  \\
  G^{(n)  \, just} \; & := \; [7]_{12}^{(n)  \, just}   \\
  G^{\sharp \, (n)  \, just} \; & := \;  [8]_{12}^{(n)  \, just}  \\
  A^{(n)  \, just} \; & := \; [9]_{12}^{(n)  \, just} \\
  A^{\sharp \, (n)  \, just}  \; & := \; [10]_{12}^{(n)  \, just}  \\
  B^{(n)  \, just} \; & := \; [11]_{12}^{ (n)  \, just }  \\
\end{eqnarray}

The effect of the \emph{just-intonation ansatz} in the
formalization of Harmony is structurally the same that the effect
of the \emph{pytagoric ansatz}. So it will be sufficient to repeat
step by step the pytagoric footsteps to get the correct
just-intonation notions in a straightforward way.

As to our goal, namely to show how these ansatz, required by the
Vertical Rules of Tonal Harmony in order to obtain the physical
consonance of notes, lead to an inconsistence in the formalization
of Harmony, I will limit myself to take into account simply the
pytagoric case, every consideration performed being trivially
translated to the just-intonation case by replacing $
{\mathbf{Kf}} $  with $ {\mathbf{Kt}} $ etc.

In terms of the introduced mathematical staff our goal may be
rephrased as the illustration of the deficiencies hidden inside $
{\mathcal{MP}}^{Pyt} $ owing to the peculiarity of pytagoric
Modulation's Theory

\smallskip

A first  obvious consideration is codified by the following:
\begin{theorem} \label{th:on the impossibility of modulating from one cycle 
to another}
\end{theorem}
ON THE IMPOSSIBILITY OF MODULATING FROM ONE CYCLE TO ANOTHER

\begin{hypothesis}
\end{hypothesis}

\begin{equation*}
    \vec{x} \, , \, \vec{y} \in {\mathbb{Z}}_{12}^{\star} \; : \;
    {\mathcal{M}} ( tonality[ \vec{x} , n ] , tonality[ \vec{y} , n
    ] \; \neq \; \emptyset \; \; n \in {\mathbb{N}} \, : n \, \leq
    \min ( maxlevel(\vec{x} ) ,  maxlevel(\vec{y} ) )
\end{equation*}

\begin{thesis}
\end{thesis}
\begin{equation}
   {\mathcal{M}}^{Pyt} \{ tonality[ \tilde{pw} ( \vec{x} , i ) , n  ] , 
tonality[ \tilde{pw} ( \vec{y} , j ) , n
   ] \} \; = \; \emptyset  \forall i \neq j
\end{equation}

\begin{proof}

Since:
\begin{equation}
  \tilde{pw} ( \vec{x} , i ) \, \bigcap \, \tilde{pw} ( \vec{x} , j
  ) \; = \; \emptyset
\end{equation}
it follows that:
\begin{equation}
    {\mathcal{P}}^{Pyt} ( \{ tonality[ \tilde{pw} ( \vec{x} , i ) , n  ] , 
tonality[ \tilde{pw} ( \vec{y} , j ) , n
   ] \} \; = \; \emptyset
\end{equation}
applying the definition\ref{def:pytagoric modulations} the thesis
immediately follows

\end{proof}

\begin{example}
\end{example}
IMPOSSIBILITY OF MODULATING AMONG PYTAGORIC C MAJOR SCALES AT
DIFFERENT CYCLES

\begin{verbatim}
  tonality[pytagoricword[majorword[0],2],1] =

{{{0,2},{4,2},{7,2}},{{2,2},{5,2},{9,2}},{{4,2},{7,2},{11,2}},{{5,2},{9,2},{0,
      2}},{{7,2},{11,2},{2,2}},{{9,2},{0,2},{4,2}},{{11,2},{2,2},{5,2}}}

tonality[pytagoricword[majorword[0],3],1] =

{{{0,3},{4,3},{7,3}},{{2,3},{5,3},{9,3}},{{4,3},{7,3},{11,3}},{{5,3},{9,3},{0,
      3}},{{7,3},{11,3},{2,3}},{{9,3},{0,3},{4,3}},{{11,3},{2,3},{5,3}}}

Intersection[{{{0,2},{4,2},{7,2}},{{2,2},{5,2},{9,2}},{{4,2},{7,2},{11,2}},{{
5,2},{9,2},{0,2}},{{7,2},{11,2},{2,2}},{{9,2},{0,2},{4,2}},{{11,2},{2,
2},{5,2}}}
,{{{0,3},{4,3},{7,3}},{{2,3},{5,3},{9,3}},{{4,3},{7,3},{11,
3}},{{5,3},{9,3},{0,3}},{{7,3},{11,3},{2,3}},{{9,3},{0,3},{4,3}},{{11,
3},{2,3},{5,3}}} ] = {}

\end{verbatim}

Let us then return to the Hindemith's citation presented in
section\ref{sec:Introduction}, i.e.:

\begin{center}
  \emph{"Anyone who has ever tasted the delights of pure intonation by the 
continual displacement of the comma in string-quartet
  playing , must come to the conclusion that there can be no such
  thing as atonal music, in which the existence of tone
  relationship is denied" (cited from chapter 4 "Harmony", section 10 
"Atonality and Polytonality" of
  \cite{Hindemith-42})}
\end{center}

and let us formalize Hindemith's "continual displacement of the
comma" inside Pytagoric Modulation Theory.

Our key tools will be the fifth cycle's raising $ C_{+} $ and
lowering $ C_{-} $ operators of, respectively, definition\ref{def:operator 
of fifth
cycle's raising} and  definition\ref{def:operator of fifth cycle's
lowering} using them, first of all, to state the following
obvious:
\begin{corollary}
\end{corollary}
\begin{equation}
  {\mathcal{M}}^{Pyt} \{ tonality[ \tilde{x} , n  ] , tonality[ \hat{C}_{+} 
( \tilde{x} ) , n ) ] \} \; = \; \emptyset
\; \; \forall  \tilde{x} \in \Xi^{\star}
\end{equation}
\begin{proof}
It is a particular case of theorem\ref{th:on the impossibility of
modulating from one cycle to another}
\end{proof}

Raising and lowering the letters of a pytagoric word in different
ways one can result to pytagoric words whose tonality may be
obtained by modulation by the starting one; I will call this kind
of modulations \emph{comma displacement modulations}.

To introduce them it is useful, first of all, to enlarge the
mathematical toolbox in the following way: given a pytagoric word
$ \tilde{x} = ( x_{1} \, , \, \cdots \, , x_{| \vec{x} | }) \in
\Xi^{\star} $

\begin{definition} \label{def:cycle moving operator applied to pytagoric 
words}
\end{definition}
\begin{equation}
    C_{\pm} ( \tilde{x}, i ) \; := \; \tilde{y} =  ( y_{1} \, ,  \,  \cdots 
\, , y_{| \vec{x} |
    }) \, : \, y_{n} = \left\{%
\begin{array}{ll}
    C_{\pm} x_{n} , & \hbox{if $ n = i $;} \\
    x_{n}, & \hbox{otherwise.} \\
\end{array}%
\right.
\end{equation}
Obviously the operators of definition\ref{def:cycle moving
operator applied to pytagoric words} easily implemented on
computer through the following expressions from the notebook of
section\ref{sec:Mathematica's notebook Mathematical Music
Toolkit}:

\begin{verbatim}
cycleraising[pytagoricword_,i_]:=
  Table[If[n==i,{Part[pytagoricword,n,1],Part[pytagoricword,n,2]+1},
      Part[pytagoricword,n]],{n,1,Length[pytagoricword]}]

cyclelowering[pytagoricword_,i_]:=
  Table[If[n==i,{Part[pytagoricword,n,1],Part[pytagoricword,n,2]-1},
      Part[pytagoricword,n]],{n,1,Length[pytagoricword]}]
\end{verbatim}

Given two pytagoric words $ \tilde{x} \, , \,  \tilde{y} \in
\Xi^{\star}$:
\begin{definition} \label{def:comma-displacement equivalence relation among 
words}
\end{definition}
$ \tilde{x} $ IS A COMMA-DISPLACEMENT OF $  \tilde{y} \; (
\tilde{x} \,  \sim_{K} \,  \tilde{y} ) $

$ \tilde{x} $ may be obtained from $   \tilde{y} $ applying to it
a suitable composition of the operators of
definition\ref{def:cycle moving operator applied to pytagoric
words}

\smallskip

Given two pytagoric tonalities $ t_{1} \, ,
\, t_{2} \, \in \, {\mathcal{T}} $:

\begin{definition} \label{def:comma-displacement equivalence relation among 
tonalities}
\end{definition}
$ t_{1} $ IS A COMMA-DISPLACEMENT OF $ t_{2} \; ( t_{1} \,
\sim_{K} \, t_{2} ) $:
\begin{equation}
  \exists  \tilde{x}_{1} , \tilde{x}_{2} \in \Xi^{\star} \ :
  t_{1} = tonality[ \tilde{x}_{1} , n ] \,  , \,   t_{2} = tonality[ 
\tilde{x}_{2} ,
  n] \, ,
  \tilde{x}_{1} \sim_{K}  \tilde{x}_{2}
\end{equation}
The notation adopted in definition \ref{def:comma-displacement
equivalence relation among words} and
definition\ref{def:comma-displacement equivalence relation among
tonalities} is justified by the fact that in both cases the fact
that comma displacement is an equivalence relation may be easily
proved.

\begin{definition}
\end{definition}
COMMA DISPLACEMENT'S PYTAGORIC MODULATIONS:
\begin{equation}
    {\mathcal{M}}^{Pyt}_{K} \; := \; \{ m \in
    {\mathcal{M}}^{Pyt} ( t_{1} ,  t_{2} ) \, : \, t_{1} ,  t_{2}
    \in {\mathcal{T}}^{Pyt} \, , \,  t_{1} \sim_{K}  t_{2} \}
\end{equation}

\begin{example} \label{ex:a comma displacement modulation of C major}
\end{example}
A COMMA DISPLACEMENT MODULATION OF C MAJOR

Given the C major pytagoric word $ \tilde{pw}_{1} $ of level zero:
\begin{verbatim}
pw1 = pytagoricword[majorword[0],0] =

{{0,0},{2,0},{4,0},{5,0},{7,0},{9,0},{11,0}}
\end{verbatim}
and that, let's call it  $ \tilde{pw}_{2} $ obtained by  $
\tilde{pw}_{1} $ raising the $ 5^{th} $ letter:
\begin{verbatim}
pw2=cycleraising[pw1,5] =

{{0,0},{2,0},{4,0},{5,0},{7,1},{9,0},{11,0}}
\end{verbatim}
let us introduce their tonalities of first level:
\begin{verbatim}
pt1=tonality[pw1,1] =

{{{0,0},{4,0},{7,0}},{{2,0},{5,0},{9,0}},{{4,0},{7,0},{11,0}},{{5,0},{9,0},{0,
      0}},{{7,0},{11,0},{2,0}},{{9,0},{0,0},{4,0}},{{11,0},{2,0},{5,0}}}

pt2=tonality[pw2,1] =

{{{0,0},{4,0},{7,1}},{{2,0},{5,0},{9,0}},{{4,0},{7,1},{11,0}},{{5,0},{9,0},{0,
      0}},{{7,1},{11,0},{2,0}},{{9,0},{0,0},{4,0}},{{11,0},{2,0},{5,0}}}
\end{verbatim}
One has that:
\begin{verbatim}
pivotaldegreesintermofdegrees[pt1,pt2] =

{{2,4,6,7},{2,4,6,7}}

\end{verbatim}
so, as to the passage in a neutral environment common to both the
involved tonalities, there is no problem.

As to cadences we have, following the definition
\ref{def:pytagoric modulations}, to consider as  $
{\mathcal{T}}_{context} $:
\begin{verbatim}
context =
Table[tonality[pytagoricword[majorword[n],0],1],{n,0,11}] =

{{{{0,0},{4,0},{7,0}},{{2,0},{5,0},{9,0}},{{4,0},{7,0},{11,0}},{{5,0},{9,0},{
        0,0}},{{7,0},{11,0},{2,0}},{{9,0},{0,0},{4,0}},{{11,0},{2,0},{5,
        
0}}},{{{1,0},{5,0},{8,0}},{{3,0},{6,0},{10,0}},{{5,0},{8,0},{0,0}},{{
        
6,0},{10,0},{1,0}},{{8,0},{0,0},{3,0}},{{10,0},{1,0},{5,0}},{{0,0},{3,
        
0},{6,0}}},{{{2,0},{6,0},{9,0}},{{4,0},{7,0},{11,0}},{{6,0},{9,0},{1,
        
0}},{{7,0},{11,0},{2,0}},{{9,0},{1,0},{4,0}},{{11,0},{2,0},{6,0}},{{1,
        
0},{4,0},{7,0}}},{{{3,0},{7,0},{10,0}},{{5,0},{8,0},{0,0}},{{7,0},{10,
        0},{2,0}},{{8,0},{0,0},{3,0}},{{10,0},{2,0},{5,0}},{{0,0},{3,0},{7,
        
0}},{{2,0},{5,0},{8,0}}},{{{4,0},{8,0},{11,0}},{{6,0},{9,0},{1,0}},{{
        
8,0},{11,0},{3,0}},{{9,0},{1,0},{4,0}},{{11,0},{3,0},{6,0}},{{1,0},{4,
        
0},{8,0}},{{3,0},{6,0},{9,0}}},{{{5,0},{9,0},{0,0}},{{7,0},{10,0},{2,
        
0}},{{9,0},{0,0},{4,0}},{{10,0},{2,0},{5,0}},{{0,0},{4,0},{7,0}},{{2,
        0},{5,0},{9,0}},{{4,0},{7,0},{10,0}}},{{{6,0},{10,0},{1,0}},{{8,0},{
        
11,0},{3,0}},{{10,0},{1,0},{5,0}},{{11,0},{3,0},{6,0}},{{1,0},{5,0},{
        8,0}},{{3,0},{6,0},{10,0}},{{5,0},{8,0},{11,0}}},{{{7,0},{11,0},{2,
        
0}},{{9,0},{0,0},{4,0}},{{11,0},{2,0},{6,0}},{{0,0},{4,0},{7,0}},{{2,
        
0},{6,0},{9,0}},{{4,0},{7,0},{11,0}},{{6,0},{9,0},{0,0}}},{{{8,0},{0,
        0},{3,0}},{{10,0},{1,0},{5,0}},{{0,0},{3,0},{7,0}},{{1,0},{5,0},{8,
        
0}},{{3,0},{7,0},{10,0}},{{5,0},{8,0},{0,0}},{{7,0},{10,0},{1,0}}},{{{
        
9,0},{1,0},{4,0}},{{11,0},{2,0},{6,0}},{{1,0},{4,0},{8,0}},{{2,0},{6,
        0},{9,0}},{{4,0},{8,0},{11,0}},{{6,0},{9,0},{1,0}},{{8,0},{11,0},{2,
        
0}}},{{{10,0},{2,0},{5,0}},{{0,0},{3,0},{7,0}},{{2,0},{5,0},{9,0}},{{
        
3,0},{7,0},{10,0}},{{5,0},{9,0},{0,0}},{{7,0},{10,0},{2,0}},{{9,0},{0,
        
0},{3,0}}},{{{11,0},{3,0},{6,0}},{{1,0},{4,0},{8,0}},{{3,0},{6,0},{10,
        
0}},{{4,0},{8,0},{11,0}},{{6,0},{10,0},{1,0}},{{8,0},{11,0},{3,0}},{{
        10,0},{1,0},{4,0}}}}
\end{verbatim}
respect to which the pytagoric tonality $ pt_{2} $ seems not to
have has cadences:
\begin{verbatim}
cadencesintermofdegrees[pt2,3,context] =

{}
\end{verbatim}
So is seems to be impossible to comma-modulate from pt1 to pt2.

An analogous situation occurs if $ \tilde{pw_{2}} $ is defined as:
\begin{equation}
    \tilde{pw_{2}} \; := \; C_{+} ( \tilde{pw_{1}} \, , \, i ) \;
    \; i =1,2,3,4
\end{equation}

\smallskip

The situation is, anyway, different for $ i \, = \, 6$:
\begin{verbatim}

pw2=cycleraising[pw1,6] =
{{0,0},{2,0},{4,0},{5,0},{7,0},{9,1},{11,0}}

pt2 =tonality[pw2,1] =
{{{0,0},{4,0},{7,0}},{{2,0},{5,0},{9,1}},{{4,0},{7,0},{11,0}},{{5,0},{9,1},{0,
      0}},{{7,0},{11,0},{2,0}},{{9,1},{0,0},{4,0}},{{11,0},{2,0},{5,0}}}

pivotaldegreesintermofdegrees[pt1,pt2] = {{1,3,5,7},{1,3,5,7}}

cadencesintermofdegrees[pt2,2,context] =
{{7},{1,7},{3,7},{5,7},{7,1},{7,3},{7,5},{7,7}}

\end{verbatim}
so that $
\{\{\{11,0\},\{2,0\},\{5,0\}\},\{\{11,0\},\{2,0\},\{5,0\}\}\} \,
\in \,  {\mathcal{M}}^{Pyt}_{K} $.

For i=7, anyway comma-modulation return to be impossible:
\begin{verbatim}
pw2=cycleraising[pw1,7] =
{{0,0},{2,0},{4,0},{5,0},{7,0},{9,0},{11,1}}

pivotaldegreesintermofdegrees[pt1,pt2] =
{{1,2,4,6},{1,2,4,6}}

cadences[pt2,2,context] = {}
\end{verbatim}

\smallskip

The possibility of performing the Hindemith's comma-displacement
not only in terms of the melody but at an harmonically
structurally consistent way occurs, as we have empirically
evidentiated, with very low frequency.

\smallskip

\begin{remark}
\end{remark}
THE STRUCTURAL INCONSISTENCE OF TONAL HARMONY VERSUS
MEYER-EPPLER'S VALENCE THEORY

Since $ {\mathbb{Q}} $ is dense in $ {\mathbb{R}} $:
\begin{equation}
card( \{ x \in ( {\mathbb{R}} - {\mathbb{Q}}) \, : \, r_{1} < x <
r_{2} \} \; = \; \aleph_{1} \; \; \forall  r_{1} \, , \, r_{2} \in
{\mathbb{Q}} , r_{1} < r_{2}
\end{equation}
one could be led to think that the structural problem shown in
this work is silly at on a physical ground in that, concretelly,
Phychoacustics have tought to us that human hear cannot
distinguish notes $ \omega $ and $ \omega \, + \, \epsilon $
nearer that a certain bound $ epsilon $.

Such a phenomenon has been formalized by Werner Meyer-Eppler
though his Valence Theory based on the following (cfr. the section
B.2. " Discriminating Tones: Werner Meyer Eppler's Valence Theory"
of \cite{Mazzola-02a}):

given two \emph{sonor signals} $ s_{1} $ and $ s_{2} $ and a
predicate P about sounds:

\begin{definition}
\end{definition}
$ s_{1} $ IS METAMERE TO $ s_{2} $ W.R.T. P ( $ s_{1} \, \sim_{P}
\, s_{2} ) $

$ s_{1} $ cannot be distinguished by $ s_{2} $ w.r.t. to P by
human listeners.

\begin{definition}
\end{definition}
VALENCE OF $ s_{1} $  W.R.T. P (P-VALENCE OF  $ s_{1} $)
\begin{equation}
    \sim_{P}  s_{1} \,  \; := \; \{ s_{2} \, : \, s_{1} \, \sim_{P}
\, s_{2} \}
\end{equation}

Identified a suitable set of predicates $ \{ P _{i} \} $ relevant
as to sonor signal distinguishing, one considers as basic
theoretical objects of Acoustics the valences w.r.t. to $
\bigwedge_{i}  P_{i} $

\medskip

In particular it may be proved that infinite just-intonation Euler
notes belong to the same $ \bigwedge_{i}  P_{i} $-valence.

As correctly stressed by Mazzola, the key bug of Meyer-Eppler's
Valence Theory consists in that the relation $ \sim_{P} $ is not
an equivalence relation lacking both of the reflexive and, the
more important point, of the transitive property.

In a less general context in which only \emph{sounds} are taken
into account and considering the valence of them with respect to
the predicate
\begin{equation}
    P ( \nu ) \,  := \,  << \, \nu \in ( \nu - \frac{\epsilon}{2} \, ,
    \, \nu + \frac{\epsilon}{2} ) \, >>
\end{equation}
one could think that, making in the definition\ref{def:physical
index of consonance of two sounds}, definition\ref{def:physical
index of consonance of N sounds} and definition \ref{def:physical
index of consonance of N notes w.r.t. a musical instrument}, the
following:
\begin{definition}
\end{definition}
P-VALENCE ANSATZ:
\begin{equation}
    \delta_{Kronecker} ( \sum_{i=1}^{N} n_{i} \omega_{i}  ) \;
    \rightarrow \;  \delta_{Kronecker}^{\epsilon} ( \sum_{i=1}^{N} n_{i} 
\omega_{i}  )
\end{equation}
where:
\begin{equation}
  \delta_{Kronecker}^{\epsilon} ( \omega ) \; := \; \left\{%
\begin{array}{ll}
    1, & \hbox{if $ | \omega | \leq \epsilon $} \\
    0, & \hbox{otherwise} \\
\end{array}%
\right.
\end{equation}
any claimed formal inconsistence of Tonal Harmony disappears.

This is not true since the structural problems concerning the
extreme difficulty of making pytagoric modulations ,and
specifically comma-modulation since any pytagoric modulation may
be though as the composition of  an "ordinary" modulation followed
by a comma-modulation), that I have shown in the simple example
\ref{ex:a comma displacement modulation of C major} are the effect
of a structural problem.

And iterating the comma-displacements, for the level of degree
that it is possible not only at a melodic but at an harmonic
level, differences lying inside the same P-valence accumulate
producing disasters as the (eventual) readers can concretelly
compute using the notebook of \ref{sec:Mathematica's notebook
Mathematical Music Toolkit}

\newpage
\appendix
\section{The Mathematica's notebook \emph{Mathematical Music Toolkit}} 
\label{sec:Mathematica's
notebook Mathematical Music Toolkit}

\begin{verbatim}

Mathematical Music Toolkit (\[Alpha]-version - 8/2/2004 at 18:10)
In[1]:= Off[General::spell1]; In[2]:=
<<DiscreteMath`Combinatorica`; In[3]:=
<<DiscreteMath`Permutations` In[4]:=
<<DiscreteMath`KroneckerDelta` In[5]:=
putinorder[listofnumbers_]:=Permute[listofnumbers,Ordering[listofnumbers]]
In[6]:= \!\(\( (*\(**\ "\<generalizedselect[l,predicate,y]\>"\)\
picks\ out\ all\
      elements\ e\_i\ of\ the\ list\ l\ for\ which\ the\ binary\ predicate\
      predicate[e\_i, y]\ is\ True\ ***) \n (*
    \(**\ N\) . B .
      \( : \ since\ the\ instruction\ is\ implemented\ throwing\ away\
        elements\ of\ pattern\ String\ \ the\ list\ l\ must\ not\ to\ 
contain
        \ \(\(strings\ !!\)!\)\)\ ***) \ \)\)
In[7]:= generalizedselect[l_,predicate_,y_]:=
  DeleteCases[
    Table[If[predicate[Part[l,i],y],Part[l,i],"throw 
away"],{i,1,Length[l]}],
    x_String]
In[8]:= (*** harmonic[nu,n] gives the n-th harmonic of the note nu
***) In[9]:= harmonic[nu_,n_]:=(n+1)*nu In[10]:= (***
scalerange[nu] gives the scale-range of the note nu ***) In[11]:=
scalerange[nu_]:=Interval[{nu,harmonic[nu,1]}] In[12]:= (***
rescalingtororange[nu,mu] rescales the note mu to the scale-
    range of the note nu ***)
In[13]:= \!\(rescalingtorange[nu_, mu_] := \
    If[mu\  < \ nu, \n\ \ \ \
      mu*\ 2\^IntegerPart[\(Log[nu] - Log[mu]\)\/Log[2] + 1], \
      If[mu > 2  nu, \nmu*2\^\(-IntegerPart[\(Log[mu] - 
Log[nu]\)\/Log[2]]\),
        \n"\<undefined\>"]]\)
In[14]:= multiplication[nu_,ratio_]:=ratio*nu In[15]:=
multiplicationfornesting[list_]:={multiplication[Part[list,1],Part[list,2]],
    Part[list,2]}
In[16]:= scaleatfixedinterval[nu_,ratio_,n_]:=
  Prepend[Table[
      rescalingtorange[nu,
        Part[Nest[multiplicationfornesting,{nu,ratio},i],1]],{i,1,n}],nu]
In[17]:= orderedscaleatfixedinterval[nu_,ratio_,n_]:=
  putinorder[scaleatfixedinterval[nu,ratio,n]]
In[18]:= (*** referencenote is the absolute reference note as
given by a suitable
    diapason ***)
In[19]:= referencenote= 132; In[20]:= \!\( (*\(**\ the\)\
following\ instructions\ implements\ the\ musical\
      alphabet\ Z\_12\ ***) \)
In[21]:= letter[n_]:=Mod[n,12] In[22]:=
alphabet=Table[letter[n],{n,0,11}]; In[23]:= \!\(\( (*\(**\ the\)\
following\ instructions\ introduce\ the\ conversion\
        among\ the\ the\ musical\ alphabet\ Z\_12\ and\ the\ notes\ of\ the\
        12 - equally\ temmpered\ chromatic\ scale\ of\ referencenote; \
    words\ in\ Z\_12\ corresponds\ to\ Allen\ Forte' s\ notion\ of\ pitch -
      class - set\ ***) \ \ \)\)
In[24]:= \!\(FROMletterTOnote[n_] :=
referencenote*2\^\(letter[n]\/12\)\) In[25]:=
FROMwordTOscale[word_]:=
  Table[FROMletterTOnote[Part[word,i]],{i,1,Length[word]}]
In[26]:= (*** a musical piece is a list each element of which is
itself a list of two
    elemente: as note and its duration ***)
In[27]:=
FROMscaleTOpiece[scale_]:=Table[{Part[scale,i],crotchet},{i,1,Length[scale]}]
In[28]:=
FROMpieceTOscale[piece_]:=Table[Part[piece,i,1],{i,1,Length[piece]}]
In[29]:= (*** the following instructions implement the two
operation of translation
    and inversion at the basis of atonal music as well as the relative
    tests ***)
In[30]:=
translation[word_,n_]:=Table[Mod[Part[word,i]+n,12],{i,1,Length[word]}]
In[31]:= translationequivalenceofwordsQ[w1_,w2_]:=
  Not[Equal[Table[Equal[w2,translation[w1,n]],{n,0,11}],Table[False,{12}]]]
In[32]:=
inversion[word_]:=Table[Mod[12-Part[word,i],12],{i,1,Length[word]}]
In[33]:=
inversionequivalenceofwordsQ[w1_,w2_]:=Equal[w1,inversion[w2]]
In[34]:= inversioninvarianceofawordQ[w_]:=Equal[w,inversion[w]]
In[35]:= \!\(\( (*\(**\ mode[word, i]\)\ gives\ the\ i -
      th\ mode\ of\ a\ words\ on\ the\ musical\ alphabet\ \ \ Z\_12\ ***) \
  \)\)
In[36]:= mode[word_,i_]:=RotateLeft[word,i-1]; In[37]:= (*** the
following instructions introduces the most used words in Music
***) In[38]:=
majorword[0]:={letter[0],letter[2],letter[4],letter[5],letter[7],letter[9],
    letter[11]}
In[39]:=
majorword[n_]:=If[n==0,majorword[0],translation[majorword[0],n]]
In[40]:=
minorword[0]:={letter[0],letter[2],letter[3],letter[5],letter[7],letter[8],
    letter[10]}
In[41]:=
minorword[n_]:=If[n==0,minorword[0],translation[minorword[0],n]]
In[42]:=
harmonicminorword[0]:={letter[0],letter[2],letter[3],letter[5],letter[7],
    letter[8],letter[11]}
In[43]:= harmonicminorword[n_]:=
  If[n==0,harmonicminorword[0],translation[harmonicminorword[0],n]]
In[44]:=
dorianword[0]:=Mod[translation[mode[majorword[0],2],-2],12]
In[45]:=
dorianword[n_]:=If[n==0,dorianword[0],translation[dorianword[0],n]]
In[46]:=
phrigianword[0]:=Mod[translation[mode[majorword[0],3],-4],12]
In[47]:=
phrigianword[n_]:=If[n==0,phrigianword[0],translation[phrigianword[0],n]]
In[48]:=
lydianword[0]:=Mod[translation[mode[majorword[0],4],-5],12]
In[49]:=
lydianword[n_]:=If[n==0,lydianword[0],translation[lydianword[0],n]]
In[50]:=
mixolydianword[0]:=Mod[translation[mode[majorword[0],5],-7],12]
In[51]:= mixolydianword[n_]:=
  If[n==0,mixolydianword[0],translation[mixolydianword[0],n]]
In[52]:=
locrianword[0]:=Mod[translation[mode[majorword[0],7],-11],12]
In[53]:=
locrianword[n_]:=If[n==0,locrianword[0],translation[locrianword[0],n]]
In[54]:=
tziganminorword[0]:={letter[0],letter[2],letter[3],letter[6],letter[7],
    letter[8],letter[11]}
In[55]:= tziganminorword[n_]:=
  If[n==0,tziganminorword[0],translation[tziganminorword[0],n]]
In[56]:=
jewishword[0]:={letter[0],letter[1],letter[4],letter[5],letter[7],letter[8],
    letter[10]}
In[57]:=
jewishword[n_]:=If[n==0,jewishword[0],translation[jewishword[0],n]]
In[58]:=
indianword[0]:={letter[0],letter[1],letter[4],letter[5],letter[7],letter[8],
    letter[11]}
In[59]:=
indianword[n_]:=If[n==0,indianword[0],translation[indianword[0],n]]
In[60]:=
majorpentatonicword[0]:={letter[0],letter[2],letter[4],letter[7],letter[9]}
In[61]:= majorpentatonicword[n_]:=
  If[n==0,majorpentatonicword[0],translation[majorpentatonicword[0],n]]
In[62]:= minorpentatonicword[0]:=
  Mod[translation[mode[majorpentatonicword[0],6],-9] ,12]
In[63]:= minorpentatonicword[n_]:=
  If[n==0,minorpentatonicword[0],translation[minorpentatonicword[0],n]]
In[64]:=
bluesword[0]:={letter[0],letter[3],letter[5],letter[6],letter[7],letter[10]}
In[65]:=
bluesword[n_]:=If[n==0,bluesword[0],translation[bluesword[0],n]]
In[66]:=
esatonalword[0]:={letter[0],letter[2],letter[4],letter[6],letter[8],
    letter[10]}
In[67]:=
esatonalword[n_]:=If[n==0,esatonalword[0],translation[esatonalword[0],n]]
In[68]:=
augmentedword[0]:={letter[0],letter[3],letter[4],letter[7],letter[8],
    letter[11]}
In[69]:=
augmentedword[n_]:=If[n==0,augmentedword[0],translation[augmentedword[0],n]]
In[70]:=
halfwholediminishedword[0]:={letter[0],letter[1],letter[3],letter[4],
    letter[6],letter[7],letter[9],letter[10]}
In[71]:= halfwholediminishedword[n_]:=
  If[n==0,halfwholediminishedword[0],
    translation[halfwholediminishedword[0],n]]
In[72]:=
wholehalfdiminishedword[0]:={letter[0],letter[2],letter[3],letter[5],
    letter[6],letter[8],letter[9],letter[11]}
In[73]:= wholehalfdiminishedword[n_]:=
  If[n==0,wholehalfdiminishedword[0],
    translation[wholehalfdiminishedword[0],n]]
In[74]:=
wholetonediminishedword[0]:={letter[0],letter[1],letter[3],letter[4],
    letter[6],letter[8],letter[10]}
In[75]:= wholetonediminishedword[n_]:=
  If[n==0,wholetonediminishedword[0],
    translation[wholetonediminishedword[0],n]]
In[76]:=
bebopmajorword[0]:={letter[0],letter[2],letter[4],letter[5],letter[7],
    letter[8],letter[9],letter[11]}
In[77]:= bebopmajorword[n_]:=
  If[n==0,bebopmajorword[0],translation[bebopmajorword[0],n]]
In[78]:=
bebopdominant[0]:={letter[0],letter[2],letter[4],letter[5],letter[7],
    letter[9],letter[10],letter[11]}
In[79]:= bebopdominantword[n_]:=
  If[n==0,bebopdominantword[0],translation[bebopdominantword[0],n]]
In[80]:= chromaticword[0]:=alphabet In[81]:=
chromaticword[n_]:=If[n==0,chromaticword[0],translation[chromaticword[0],n]]
In[82]:= (*** randomword[n] gives a pseudo-
    random word on the musical alphabet of length n ***)
In[83]:= randomword[n_]:=Table[letter[ Random[Integer,{0,11}]
],{i,1,n}] In[84]:= (*** words[n] is the list of all the words of
length n in lexicographic
    ordering ***)
In[85]:= words[n_]:=Strings[chromaticword[0],n] In[86]:=
nonrepetitiveQ[l_]:=
  Equal[Table[
      If [Equal[Part[l,i],Part[l,j]],1,0],{i,1,Length[l]},{j,1,Length[l]}],
    IdentityMatrix[Length[l]]]
In[87]:= nonrepetitivewords[n_]:=Select[words[n],nonrepetitiveQ]
In[88]:= (*** "wordsupto[n]" is the list of all the words of
length less or equal to n
    in lexicographic ordering ***)
In[89]:=
wordsupto[n_]:=If[n==1,words[1],Join[wordsupto[n-1],words[n]]]
In[90]:= (*** "subwords[x]" gives the list of all the sub-
    words of the word x in lexicographic ordering ***)
In[91]:= subwords[x_]:=LexicographicSubsets[x] In[92]:= (***
"word[n]" gives the n^th word in  global lexicographic ordering
***) (*** N.B.: up to now it is defined only for n < 16 ***)
In[93]:= word[n_]:=Part[wordsupto[15],n] In[94]:= (***
"locallexicographicnumber[x]" gives the  local lexicographic
number of
    the word x,
  i.e. the lexicographic number of x w.r.t. to the set of the words having
    the same length of x ***)
In[95]:=
locallexicographicnumber[x_]:=Part[Flatten[Position[words[Length[x]],x]],1]
In[96]:= (*** "globallexicographicnumber[x]" gives the global
lexicographic number of
    the word x,
  i.e. the lexicographic number w.r.t. to the set of the words of any
    length   ***)
In[97]:= globallexicographicnumber[x_]:=
  If[Length[x]==1,locallexicographicnumber[x],
    Length[wordsupto[Length[x]-1]]+locallexicographicnumber[x]]
In[98]:= (*** "quasilexicographicnumber[x]" gives the
quasilexicographic number of
    the word x,
  i.e. the lexicographic number w.r.t. to the set of the words of any length
    considering also the null string  ***)
In[99]:=
quasilexicographicnumber[x_]:=If[x=={},0,globallexicographicnumber[x]+1]
In[100]:= (*** prefix[x,n] is the n-long prefix of the word x ***)
In[101]:= prefix[x_,n_]:=Part[x,Table[i,{i,1,n}]] In[102]:= (***
prefixes[x] is the list of the prefixes of the word x ***)
In[103]:= prefixes[x_]:=Table[prefix[x,n],{n,1,Length[x]-1}]
In[104]:= (*** prefixQ[x,y] is the predicate stating that x is a
prefix of y ***) In[105]:= prefixQ[x_,y_]:=MemberQ[prefixes[y],x]
In[106]:= (*** the following instructions introduced the notes
without the Mod-12
      constraint ***)
In[107]:= c[1]=FROMletterTOnote[0]; In[108]:= c[n_]:=2*c[n-1]
In[109]:= c\[Sharp][1]=FROMletterTOnote[1]; In[110]:=
c\[Sharp][n_]:=2*c\[Sharp][n-1] In[111]:=
d[1]=FROMletterTOnote[2]; In[112]:= d[n_]:=2*d[n-1] In[113]:=
d\[Sharp][1]=FROMletterTOnote[3]; In[114]:=
d\[Sharp][n_]:=2*d\[Sharp][n-1] In[115]:=
e[1]=FROMletterTOnote[4]; In[116]:= e[n_]:=2*e[n-1] In[117]:=
f[1]=FROMletterTOFROMletterTOnote[5]; In[118]:= f[n_]:=2*f[n-1]
In[119]:= f\[Sharp][1]=FROMletterTOnote[6]; In[120]:=
f\[Sharp][n_]:=2*f\[Sharp][n-1] In[121]:=
g[1]=FROMletterTOnote[7]; In[122]:= g[n_]:=2*g[n-1] In[123]:=
g\[Sharp][1]=FROMletterTOnote[8]; In[124]:=
g\[Sharp][n_]:=2*g\[Sharp][n-1] In[125]:=
a[1]=FROMletterTOnote[9]; In[126]:= a[n_]:=2*a[n-1] In[127]:=
a\[Sharp][1]=FROMletterTOnote[10]; In[128]:=
a\[Sharp][n_]:=2*a\[Sharp][n-1] In[129]:=
b[1]=FROMletterTOnote[11]; In[130]:= b[n_]:=2*b[n-1] In[131]:=
(*** insert the "referencetime" ***) In[132]:= referencetime=4;
In[133]:= semibreve=1*referencetime; In[134]:= minim=
(1/2)*referencetime; In[135]:= crotchet= (1/4)*referencetime;
In[136]:= quaver=(1/8)*referencetime; In[137]:=
semiquaver=(1/16)*referencetime; In[138]:=
demisemiquaver=(1/32)*referencetime; In[139]:=
hemidemisemiquaver=(1/64)*referencetime; In[140]:= (*** the
following instructions introduce the formalism of Euler's musical
    space ***)
In[141]:=
eulercoordination[x_List]:=
  Product[Power[Prime[i],Part[x,i]], {i,1,Length[x]}]
In[142]:= FROMeulerpointTOnote[eulpoint_List] :=
  referencenote*2^Part[eulpoint,1]*3^Part[eulpoint,2]*5^Part[eulpoint,3]
In[143]:= \!\(FROMnoteTOpitch[nu_] := 1200\/Log[2]*Log[nu\/\(\
referencenote\)\ ]\) In[144]:= \!\(FROMpitchTOnote[pitch_] := \ \
referencenote*Exp[Log[2]\/1200*pitch]\) In[145]:=
Hprime=Table[Log[Prime[i]],{i,1,3}]; In[146]:=
\!\(FROMeulerpointTOpitch[eulpoint_List] :=
    1200\/Log[2]*\ Dot[eulpoint, Hprime]\)
In[147]:=
lastpartisequaltosomethingQ[x_,something_]:=Equal[Last[x],something]
In[148]:= rationalsupto[n_]:=Flatten[
Table[Rational[i,j],{i,0,n},{j,1,n}] ] In[149]:=
FROMindexTOrational[index_,n_]:=Part[rationalsupto[n],index]
In[150]:= FROMlistofindicesTOlistofrationals[list_,n_]:=
  Table[FROMindexTOrational[Part[list,i],n],{i,1,Length[list]}]
In[151]:= FROMnoteTOeulerpoint[nu_,n_]:=
  FROMlistofindicesTOlistofrationals[
    First[generalizedselect[
        Flatten[Table[ {{i,j,k},
              Dot[{Part[ rationalsupto[2],i ],Part[ rationalsupto[2],j],
                  Part[ rationalsupto[2],k ]},Hprime]},{i,1,
              Length[ rationalsupto[2] ]},{j,1,Length[ rationalsupto[2] 
]},{k,
              1,Length[ rationalsupto[2] ]}],2]  ,
        lastpartisequaltosomethingQ,FROMnoteTOpitch[nu]]],n]
In[152]:= FROMwordTOlistofeulerpoints[w_]:=
  Table[FROMnoteTOeulerpoint[Part[FROMwordTOscale[w] ,i]],{i,1,Length[w]}]
In[153]:= octaveepoint={1,0,0}; In[154]:= fifthpoint={0,1,0};
In[155]:= thirdpoint={0,0,1}; In[156]:=
canonicalnotesbasis:={octaveepoint,fifthpoint,thirdpoint}
In[157]:=
canonicalintervalsbasis:={octaveepoint,fifthpoint-octaveepoint,
    thirdpoint-2octaveepoint}
In[158]:= fifthcomma =-7*octaveepoint+ 12
*(fifthpoint-octaveepoint); In[159]:= thirdcomma=
  2*octaveepoint-4*(fifthpoint-octaveepoint)+thirdpoint-2octaveepoint;
In[160]:= gradussuavitatis[n_Integer]:=
  1+Sum[  ( FactorInteger[n][[i]][[1]]-1)* FactorInteger[n][[i]][[2]]  , 
{i,1,
        Length[FactorInteger[n]]}]
In[161]:= gradussuavitatis[r_Rational]:=
  gradussuavitatis[
    Times[Numerator[r]*Denominator[r],
      Power[GCD[Numerator[r],Power[Denominator[r],2]] , -1 ]]]
In[162]:= (***  chord[word,i,level] gives the chord of the i-
    th degree of a word at a chosen level ***)
In[163]:= chord[word_,i_,level_]:=
  Table[Part[mode[word,i],
      
If[Mod[2n+1,Length[word]]==0,Length[word],Mod[2n+1,Length[word]]]],{n,0,
      1+level}]
In[164]:= (*** numberoflevels of a word gives the possible chords
level of that
    word ***)
In[165]:= \!\(numbersoflevels[word_] :=
    1 + IntegerPart[\(Length[word] - 5\)\/2] -
      \(\((\(-1\))\)\^Length[word] - 1\)\/2*IntegerPart[Length[word]\/2]\)
In[166]:= (*** the following instructions introduce the basic
concept of Tonal
    Harmony ***)
In[167]:=
tonality[word_,level_]:=Table[chord[word,i,level],{i,1,Length[word]}]
In[168]:= pivotaldegrees[t1_,t2_]:=Intersection[t1,t2] In[169]:=
harmonicwords[t_,n_]:=Strings[t,n] In[170]:=
harmonicwordsupto[t_,n_]:=
  If[n==1,harmonicwords[t,1],
    Join[harmonicwordsupto[t,n-1],harmonicwords[t,n]]]
In[171]:= harmonicwordintermofdegrees[t_,listofdegrees_]:=
  Table[Part[t,Part[listofdegrees,i]],{i,1,Length[listofdegrees]}]
In[172]:=
degreeofachordinatonality[w_,t_]:=Part[Flatten[Position[t,w]],1]
In[173]:= (*** "pivotaldegreesintermofdegrees[t1,t2]" give a list
of two elements,
  the first being the pivotal degrees of "t1" and "t2" expressed as degrees
    of "t1", the second being the pivotal degrees of "t1" and "t2" expressed
    as degrees of t2" ***)
In[174]:= pivotaldegreesintermofdegrees[t1_,t2_]:={
    Table[degreeofachordinatonality[Part[pivotaldegrees[t1,t2],i],t1],{i,1,
        Length[pivotaldegrees[t1,t2]]}],
    Table[degreeofachordinatonality[Part[pivotaldegrees[t1,t2],i],t2],{i,1,
        Length[pivotaldegrees[t1,t2]]}]}
In[175]:= FROMharmonicwordTOphysicalchord[hw_]:=
  Table[FROMwordTOscale[Part[hw,i]],{i,1,Length[hw]}]
In[176]:=
tonalitymembershipQ[hw_,t_]:=MemberQ[harmonicwords[t,Length[hw]],hw]
In[177]:=
setofthemajortonalities[n_]:=Table[tonality[majorword[i],n],{i,0,11}];
In[178]:=
setoftheminortonalities[n_]:=Table[tonality[minorword[i],n],{i,0,11}];
In[179]:= setoftheclassicaltonalities[n_]:=
  Union[setofthemajortonalities[n],setoftheminortonalities[n]];
In[180]:= setofthegregoriantonalities[n_]:=
  Flatten[Table[tonality[mode[majorword[i],j],n],{j,1,7},{i,0,11}],1];
In[181]:= mazzolawords:=nonrepetitivewords[7] In[182]:=
mazzolaword[n_]:=Part[mazzolawords,n] In[183]:=
setofthemazzolatonalities[n_]:=
  
Table[tonality[Part[mazzolawords,i],n],{i,1,Length[nonrepetitivewords[7]]}];
In[184]:=
setofthejewishtonalities[n_]:=Table[tonality[jewishword[i],n],{i,0,11}];
In[185]:= truthQ[x_]:=Equal[x,True] In[186]:=
cadenceQ[hw_,setoftonalities_]:=
  If[Length[
        Select[Table[
            tonalitymembershipQ[hw,Part[setoftonalities,i]],{i,1,
              Length[setoftonalities]}],truthQ]]==1,True,False]
In[187]:= subharmonicwords[hw_]:=LexicographicSubsets[hw]
In[188]:= cadences[t_,n_,setoftonalities_]:=
  generalizedselect[harmonicwordsupto[t,n],cadenceQ,setoftonalities]
In[189]:= cadencesintermofdegrees[t_,n_,setoftonalities_]:=
  Table[degreeofachordinatonality[
      Part[Part[cadences[t,n,setoftonalities],i],j],t],{i,1,
      Length[cadences[t,n,setoftonalities] ]},{j,1,
      Length[ Part[cadences[t,n,setoftonalities],i] ]}]
In[190]:=
minimalcadenceQ[hw_,setoftonalities_]:=And[cadenceQ[hw,setoftonalities],
        
Equal[Table[cadenceQ[prefix[hw,i],setoftonalities],{i,1,Length[hw]-1}],
      Table[False,{Length[hw]-1}]]]
In[191]:= minimalcadences[t_,n_,setoftonalities_]:=
  generalizedselect[harmonicwords[t,n],minimalcadenceQ,setoftonalities]
In[192]:= minimalcadencesintermofdegrees[t_,n_,setoftonalities_]:=
  Table[degreeofachordinatonality[
      Part[Part[minimalcadences[t,n,setoftonalities],i],j],t],{i,1,
      Length[minimalcadences[t,n,setoftonalities] ]},{j,1,
      Length[ Part[minimalcadences[t,n,setoftonalities],i] ]}]
In[193]:= (*** a monodic piece is a list each element of which is
itself a list of the
    form {note ,duration}   ***)
In[194]:= playmonodic[piece_]:=
  Do[Play[Sin[piece[[i]][[1]] *2*\[Pi]*t],{t,0,piece[[i]][[2]]}] , {i,1,
      Length[piece]}]
In[195]:=
playwwordasarpeggio[w_]:=playmonodic[FROMscaleTOpiece[FROMwordTOscale[w]]]
In[196]:= playwordaschord[w_,time_]:=
  Play[Sum[Sin[Part[FROMletterTOnote[w],i]*2*\[Pi]*t],{i,1,Length[w]}],{t,0,
      time}]
In[197]:= playharmonicword[hw_,time_]:=
  Do[playwordaschord[Part[hw,i],time],{i,1,Length[hw]}]
In[198]:= (*** a polyphonic piece is a list each element of which
is itself a list of
    the form {word ,duration} ***)
In[199]:= playpolyphonic[piece_]:=
  Do[playwordaschord[piece[[i]][[1]],piece[[i]][[2]]] , {i,1,Length[piece]}]
In[200]:= translationofharmonicword[hw_,n_]:=
  Table[translation[Part[hw,i],n],{i,1,Length[hw]}]
In[201]:=
inversionofharmonicword[hw_]:=Table[inversion[Part[hw,i]],{i,1,Length[hw]}]
In[202]:= translationequivalenceofharmonicwordsQ[hw1_,hw2_]:=
  And[Length[hw1]==Length[hw2],
    Equal[Table[
        translationequivalenceofwordsQ[Part[hw1,i],Part[hw2,i]],{i,1,
          Length[hw1]}],Table[True,{Length[hw1]}]]]
In[203]:= inversionequivalenceofharmonicwordsQ[hw1_,hw2_]:=
  And[Length[hw1]==Length[hw2],
    Equal[Table[
        inversionequivalenceofwordsQ[Part[hw1,i],Part[hw2,i]],{i,1,
          Length[hw1]}],Table[True,{Length[hw1]}]]]
In[204]:= inversioninvarianceofanharmonicwordQ[hw_]:=
  Equal[Table[inversioninvarianceofawordQ[Part[hw,i]],{i,1,Length[hw]}],
    Table[True,{Length[w]}]]
In[205]:= stepalongfifthscycle[n_]:=Mod[n+7,12] In[206]:=
fifthscycle=NestList[stepalongfifthscycle,0,12]; In[207]:=
modulation[hw1_,hw2_,t1_,t2_,np_,cad2_,setoftonalities_]:=
  If[ Not[tonalitymembershipQ[hw1,t1]],False,
    If[ Not[tonalitymembershipQ[hw2,t2]],False,
      If[Not[MemberQ[pivotaldegrees[t1,t2],np ]],False,
        If[Not[MemberQ[cadences[t2,Length[np],setofonalities],cad2]],False,
          Join[hw1,np,cad2,hw2]]]]]
In[208]:= Mazzolamodulator[t1_,t2_]:=
  If[translationequivalenceofharmonicwordsQ[t1,t2],
    If[inversioninvarianceofanharmonicwordQ[t1],
      letter[n_]
        \[Rule]inversion[translation[ letter[n_] ,Part[t2,1]-Part[t1,1]]],
      letter[n_] \[Rule]translation[ letter[n_] ,Part[t2,1]-Part[t1,1]] ,
      "undefined"]]
In[209]:= MazzolamodulationQ[l_,t1_,t2_,n_,setoftonalities_]:=
  And[Equal[Part[l,1],Mazzolamodulator[t1,t2]],
    MemberQ[Part[l,2],mimimalcadentialsets[t,n,setoftonalities]]]
In[210]:= symmetrytransformation[g_]:=Rule[x_,g[x]] In[211]:=
symmetricQ[w_,g_]:=Equal[w,ReplaceAll[w,symmetrytransformation[g]]]
In[212]:= (*** an arbitarry law of harmony may be always codified
in terms of a
    predicate concerning tonalities and harmonic words ***)
In[213]:= lawofresolutiononthetonicQ[hw_,t_]:=
  And[tonalitymembershipQ[hw,t], Equal[Last[hw],Part[t,1]]]
In[214]:= specialtransformationrule[n_]:=
Rule[majorword[x_],majorword[x+n]] In[215]:=
specialinvarianceQ[n_,lawofharmomyQ_,argumentsofthelawofharmony__]:=
  Equal[lawofharmonyQ,
    ReplaceAll[lawofharmonyQ[argumentsofthelawofharmony],
      specialtransformationrule[n]]]
In[216]:= morethanspecialtransformationrule[n_]:=
Rule[mazzolaword[x_],mazzolaword[x+n]] In[217]:=
morethanspecialinvarianceQ[n_,lawofharmomyQ_,
argumentsofthelawofharmony__]:=
  Equal[lawofharmonyQ,
    ReplaceAll[lawofharmonyQ[argumentsofthelawofharmony],
      specialtransformationrule[n]]]
In[218]:= generaltransformationrule[y_scale]:=Rule[x_scale , y]
In[219]:=
generalinvariance[lawofharmonyQ,argumentsofthelawofharmony__,y_scale]:=
  Equal[x,ReplaceAll[lawofharmonyQ[argumentsofthelawofharmony],
      generaltransformationrule[y]]]
In[220]:= FROMnoteTOpulsation[nu_]:=Times[nu,Power[2*Pi,-1]]
In[221]:= FROMpulsationTOnote[omega_]:=2*Pi*omega In[222]:=
musicalinstrumentQ[a_]:=Sum[a[n]*a[-n],{n,-Infty,+Infty}]< + Infty
In[223]:= (*** a sound may be codified as a list of the form
{omega,a} where "omega" is
    a  pulsation,while "a" is a musical instrument ****)
In[224]:= FROMpulsationTOsound[omega_,instrument_]:=
  If[musicalinstrumentQ[instrument],{omega,instrument} , "undefined"]
In[225]:= FROMsoundTOpulsation[sound_]:=Part[sound,1] In[226]:=
physicalindexofconsonanceofsounds[sound1_,sound2_]:=
  Sum[KroneckerDelta[n*Part[sound1,1]+m*Part[sound2,1]]*(
        Part[sound1,2][n]+Part[sound2,2][m]),{n,-Infty,+Infty},{
      m,-Infty,+Infty}]
In[227]:=
physicalindexofconsonanceofpulsation[omega1_,omega2_,instrument_]:=
  If[musicalinstrumentQ[instrument],
    
physicalindexofconsonanceofsounds[FROMpulsationTOsound[omega1,instrument],
      FROMpulsationTOsound[omega2,instrument]]  , "undefined"]
In[228]:= commensurabilityQ[omega1_,omega2_]:=
  Equal[ Head[Times[omega1,Power[omega2,-1]]],Rational]
In[229]:= (*** the following instructions implement the pytagoric
ansatz ***) In[230]:= pytagoricletter[i_,n_]:={letter[i],n}
In[231]:= pytagoricalphabetuptofifthcycles[m_]:=
  Flatten[Table[pytagoricletter[n,j],{j,0,m},{n,0,11}],1]
In[232]:= FROMpytagoricletterTOnote[x_]:=
  Part[putinorder[
      Take[scaleatfixedinterval[referencenote,3,12*Part[x,2]+11],-12]],
    Part[x,1]+1]
In[233]:= pytagoricscaleuptofifthcycles[m_]:=
  Map[FROMpytagoricletterTOnote,pytagoricalphabetuptofifthcycles[m]]
In[234]:=
pytagoricword[w_,n_]:=Table[pytagoricletter[Part[w,i],n],{i,1,Length[w]}]
In[235]:= \!\(empiricalsimplicitymeasure[r_Rational] := \
    1\/GCD[Numerator[r], Denominator[r]]*
      \(Numerator[r] + Denominator[r]\)\/\(Numerator[r]*Denominator[r]\)\)
In[236]:=
vogelchromaticjustlistofeulerpoints:={{0,0,0},{4,-1,-1},{-3,2,0},{1,1,-1},{-2,
      0,1},{2,-1,0},{-5,2,1},{-1,1,0},{3,0,-1},{0,-1,1},{4,-2,0}{-3,1,1}}
In[237]:= cycleraising[pytagoricword_,i_]:=
  Table[If[n==i,{Part[pytagoricword,n,1],Part[pytagoricword,n,2]+1},
      Part[pytagoricword,n]],{n,1,Length[pytagoricword]}]
In[238]:= cyclelowering[pytagoricword_,i_]:=
  Table[If[n==i,{Part[pytagoricword,n,1],Part[pytagoricword,n,2]-1},
      Part[pytagoricword,n]],{n,1,Length[pytagoricword]}]


\end{verbatim}
\newpage
\subsection{Some useful palette}
\begin{description}
    \item[Basic words' palette]

\bigskip

\begin{verbatim}

(***********************************************************************

                    Mathematica-Compatible Notebook

This notebook can be used on any computer system with Mathematica
3.0, MathReader 3.0, or any compatible application. The data for
the notebook starts with the line of stars above.

To get the notebook into a Mathematica-compatible application, do
one of the following:

* Save the data starting with the line of stars above into a file
  with a name ending in .nb, then open the file inside the application;

* Copy the data starting with the line of stars above to the
  clipboard, then use the Paste menu command inside the application.

Data for notebooks contains only printable 7-bit ASCII and can be
sent directly in email or through ftp in text mode.  Newlines can
be CR, LF or CRLF (Unix, Macintosh or MS-DOS style).

NOTE: If you modify the data for this notebook not in a
Mathematica- compatible application, you must delete the line
below containing the word CacheID, otherwise
Mathematica-compatible applications may try to use invalid cache
data.

For more information on notebooks and Mathematica-compatible
applications, contact Wolfram Research:
  web: http://www.wolfram.com
  email: info@wolfram.com
  phone: +1-217-398-0700 (U.S.)

Notebook reader applications are available free of charge from
Wolfram Research.
***********************************************************************)

(*CacheID: 232*)


(*NotebookFileLineBreakTest NotebookFileLineBreakTest*)
(*NotebookOptionsPosition[      4116,        117]*)
(*NotebookOutlinePosition[      5172,        156]*) (*
CellTagsIndexPosition[      5128,        152]*)
(*WindowFrame->Palette*)



Notebook[{ Cell[BoxData[GridBox[{
        {
          ButtonBox[\(majorword[\[SelectionPlaceholder]]\)]},
        {
          ButtonBox[\(minorword[\[SelectionPlaceholder]]\)]},
        {
          ButtonBox[\(harmonicminorword[\[SelectionPlaceholder]]\)]},
        {
          ButtonBox[\(dorianword[\[SelectionPlaceholder]]\)]},
        {
          ButtonBox[\(phrigianword[\[SelectionPlaceholder]]\)]},
        {
          ButtonBox[\(lydianword[\[SelectionPlaceholder]]\)]},
        {
          ButtonBox[\(mixolydianword[\[SelectionPlaceholder]]\)]},
        {
          ButtonBox[\(locrianword[\[SelectionPlaceholder]]\)]},
        {
          ButtonBox[\(tziganminorword[\[SelectionPlaceholder]]\)]},
        {
          ButtonBox[\(jewishword[\[SelectionPlaceholder]]\)]},
        {
          ButtonBox[\(majorpentatonicword[\[SelectionPlaceholder]]\)]},
        {
          ButtonBox[\(minorpentatonicword[\[SelectionPlaceholder]]\)]},
        {
          ButtonBox[\(bluesword[\[SelectionPlaceholder]]\)]},
        {
          ButtonBox[\(esatonalword[\[SelectionPlaceholder]]\)]},
        {
          ButtonBox[\(augmentedword[\[SelectionPlaceholder]]\)]},
        {
          ButtonBox[\(halfwholediminishedword[\[SelectionPlaceholder]]\)]},
        {
          ButtonBox[\(wholehalfdiminishedword[\[SelectionPlaceholder]]\)]},
        {
          ButtonBox[\(wholetonediminishedword[\[SelectionPlaceholder]]\)]},
        {
          ButtonBox[\(bebopmajorword[\[SelectionPlaceholder]]\)]},
        {
          ButtonBox[\(bebopdominant[\[SelectionPlaceholder]]\)]},
        {
          ButtonBox[\(chromaticword[\[SelectionPlaceholder]]\)]},
        {
          ButtonBox[\(randomword[\[SelectionPlaceholder]]\)]}
        },
      RowSpacings->0,
      ColumnSpacings->0,
      GridDefaultElement:>ButtonBox[ "\\[Placeholder]"]]], NotebookDefault,
  CellMargins->{{Inherited, Inherited}, {5, Inherited}},
  Evaluatable->True,
  CellGroupingRules->"InputGrouping",
  PageBreakAbove->True,
  PageBreakWithin->False,
  GroupPageBreakWithin->False,
  CellLabelMargins->{{11, Inherited}, {Inherited, Inherited}},
  DefaultFormatType->DefaultInputFormatType,
  LineSpacing->{1.25, 0},
  AutoItalicWords->{},
  FormatType->InputForm,
  ScriptMinSize->9,
  ShowStringCharacters->True,
  NumberMarks->True,
  CounterIncrements->"Input",
  StyleMenuListing->None,
  FontFamily->"Courier",
  FontWeight->"Bold"]
}, FrontEndVersion->"Microsoft Windows 3.0", ScreenRectangle->{{0,
1024}, {0, 712}}, Editable->False, WindowToolbars->{},
PageWidth->441, WindowSize->{Fit, Fit}, WindowMargins->{{60,
Automatic}, {Automatic, 33}}, WindowFrame->"Palette",
WindowElements->{}, WindowFrameElements->"CloseBox",
WindowClickSelect->False,
ScrollingOptions->{"PagewiseScrolling"->True},
ShowCellBracket->False, CellMargins->{{0, 0}, {Inherited, 0}},
Active->True, CellOpen->True, ShowCellLabel->False,
ShowCellTags->False, ImageMargins->{{0, Inherited}, {Inherited,
0}}, Magnification->1 ]


(***********************************************************************
Cached data follows.  If you edit this Notebook file directly, not
using Mathematica, you must remove the line containing CacheID at
the top of the file.  The cache data will then be recreated when
you save this file from within Mathematica.
***********************************************************************)

(*CellTagsOutline CellTagsIndex->{} *)

(*CellTagsIndex CellTagsIndex->{} *)

(*NotebookFileOutline Notebook[{ Cell[1710, 49, 2402, 66, 362,
NotebookDefault,
  Evaluatable->True,
  CellGroupingRules->"InputGrouping",
  PageBreakAbove->True,
  PageBreakWithin->False,
  CounterIncrements->"Input"]
} ] *)




(***********************************************************************
End of Mathematica Notebook file.
***********************************************************************)

\end{verbatim}

\bigskip

\item[basic harmonic construction's palette]

\bigskip

\begin{verbatim}

(***********************************************************************

                    Mathematica-Compatible Notebook

This notebook can be used on any computer system with Mathematica
3.0, MathReader 3.0, or any compatible application. The data for
the notebook starts with the line of stars above.

To get the notebook into a Mathematica-compatible application, do
one of the following:

* Save the data starting with the line of stars above into a file
  with a name ending in .nb, then open the file inside the application;

* Copy the data starting with the line of stars above to the
  clipboard, then use the Paste menu command inside the application.

Data for notebooks contains only printable 7-bit ASCII and can be
sent directly in email or through ftp in text mode.  Newlines can
be CR, LF or CRLF (Unix, Macintosh or MS-DOS style).

NOTE: If you modify the data for this notebook not in a
Mathematica- compatible application, you must delete the line
below containing the word CacheID, otherwise
Mathematica-compatible applications may try to use invalid cache
data.

For more information on notebooks and Mathematica-compatible
applications, contact Wolfram Research:
  web: http://www.wolfram.com
  email: info@wolfram.com
  phone: +1-217-398-0700 (U.S.)

Notebook reader applications are available free of charge from
Wolfram Research.
***********************************************************************)

(*CacheID: 232*)


(*NotebookFileLineBreakTest NotebookFileLineBreakTest*)
(*NotebookOptionsPosition[      8129,        231]*)
(*NotebookOutlinePosition[      9185,        270]*) (*
CellTagsIndexPosition[      9141,        266]*)
(*WindowFrame->Palette*)



Notebook[{ Cell[BoxData[GridBox[{
        {
          ButtonBox[\(letter[\[SelectionPlaceholder]]\)],

          ButtonBox[
            \(cadences[\[SelectionPlaceholder], \[SelectionPlaceholder],
              \[SelectionPlaceholder]]\)]},
        {
          ButtonBox[\(words[\[SelectionPlaceholder]]\)],

          ButtonBox[
            \(cadencesintermofdegrees[\[SelectionPlaceholder],
              \[SelectionPlaceholder], \[SelectionPlaceholder]]\)]},
        {
          ButtonBox[\(wordsupto[\[SelectionPlaceholder]]\)],

          ButtonBox[
            \(minimalcadences[\[SelectionPlaceholder],
              \[SelectionPlaceholder], \[SelectionPlaceholder]]\)]},
        {
          ButtonBox[\(nonrepetitivewords[\[SelectionPlaceholder]]\)],

          ButtonBox[
            \(minimalcadencesintermofdegrees[\[SelectionPlaceholder],
              \[SelectionPlaceholder], \[SelectionPlaceholder]]\)]},
        {

          ButtonBox[
            \(translation[\[SelectionPlaceholder], \[SelectionPlaceholder]]
              \)],

          ButtonBox[
            \(translationofharmonicword[\[SelectionPlaceholder],
              \[SelectionPlaceholder]]\)]},
        {
          ButtonBox[\(inversion[\[SelectionPlaceholder]]\)],
          ButtonBox[\(inversionofharmonicword[\[SelectionPlaceholder]]\)]},
        {

          ButtonBox[
            \(translationequivalenceofwordsQ[\[SelectionPlaceholder],
              \[SelectionPlaceholder]]\)],

          ButtonBox[
            
\(translationequivalenceofharmonicwordsQ[\[SelectionPlaceholder],
              \[SelectionPlaceholder]]\)]},
        {

          ButtonBox[
            \(inversionequivalenceofwordsQ[\[SelectionPlaceholder],
              \[SelectionPlaceholder]]\)],

          ButtonBox[
            \(inversionequivalenceofharmonicwordsQ[\[SelectionPlaceholder],
              \[SelectionPlaceholder]]\)]},
        {

          ButtonBox[
            \(inversioninvarianceofawordQ[\[SelectionPlaceholder]]\)],

          ButtonBox[
            \(inversioninvarianceofanharmonicwordQ[\[SelectionPlaceholder]]
              \)]},
        {

          ButtonBox[
            \(mode[\[SelectionPlaceholder], \[SelectionPlaceholder]]\)],
          ButtonBox[
            \(Mazzolamodulator[\[SelectionPlaceholder],
              \[SelectionPlaceholder]]\)]},
        {

          ButtonBox[
            \(chord[\[SelectionPlaceholder], \[SelectionPlaceholder],
              \[SelectionPlaceholder]]\)],

          ButtonBox[
            \(MazzolamodulationQ[\[SelectionPlaceholder],
              \[SelectionPlaceholder], \[SelectionPlaceholder],
              \[SelectionPlaceholder], \[SelectionPlaceholder]]\)]},
        {

          ButtonBox[
            \(tonality[\[SelectionPlaceholder], \[SelectionPlaceholder]]\)],

          ButtonBox[
            \(modulation[\[SelectionPlaceholder], \[SelectionPlaceholder],
              \[SelectionPlaceholder], \[SelectionPlaceholder],
              \[SelectionPlaceholder], \[SelectionPlaceholder],
              \[SelectionPlaceholder]]\)]},
        {

          ButtonBox[
            \(pivotaldegrees[\[SelectionPlaceholder], 
\[SelectionPlaceholder]]
              \)],

          ButtonBox[
            \(minimalcadenceQ[\[SelectionPlaceholder],
              \[SelectionPlaceholder], \[SelectionPlaceholder]]\)]},
        {

          ButtonBox[
            \(harmonicwords[\[SelectionPlaceholder], 
\[SelectionPlaceholder]]
              \)],
          ButtonBox[\({\[SelectionPlaceholder]}\)]},
        {

          ButtonBox[
            \(harmonicwordsupto[\[SelectionPlaceholder],
              \[SelectionPlaceholder]]\)],
          ButtonBox[\({\[SelectionPlaceholder], 
\[SelectionPlaceholder]}\)]},
        {

          ButtonBox[
            \(harmonicwordintermofdegrees[\[SelectionPlaceholder],
              \[SelectionPlaceholder]]\)],

          ButtonBox[
            \({\[SelectionPlaceholder], \[SelectionPlaceholder],
              \[SelectionPlaceholder]}\)]},
        {

          ButtonBox[
            \(tonalitymembershipQ[\[SelectionPlaceholder],
              \[SelectionPlaceholder]]\)],

          ButtonBox[
            \({\[SelectionPlaceholder], \[SelectionPlaceholder],
              \[SelectionPlaceholder], \[SelectionPlaceholder]}\)]},
        {

          ButtonBox[
            \(degreeofachordinatonality[\[SelectionPlaceholder],
              \[SelectionPlaceholder]]\)],

          ButtonBox[
            \({\[SelectionPlaceholder], \[SelectionPlaceholder],
              \[SelectionPlaceholder], \[SelectionPlaceholder],
              \[SelectionPlaceholder]}\)]},
        {

          ButtonBox[
            \(pivotaldegreesintermofdegrees[\[SelectionPlaceholder],
              \[SelectionPlaceholder]]\)],

          ButtonBox[
            \({\[SelectionPlaceholder], \[SelectionPlaceholder],
              \[SelectionPlaceholder], \[SelectionPlaceholder],
              \[SelectionPlaceholder], \[SelectionPlaceholder]}\)]},
        {

          ButtonBox[
            \(cadenceQ[\[SelectionPlaceholder], \[SelectionPlaceholder]]\)],

          ButtonBox[
            \({\[SelectionPlaceholder], \[SelectionPlaceholder],
              \[SelectionPlaceholder], \[SelectionPlaceholder],
              \[SelectionPlaceholder], \[SelectionPlaceholder],
              \[SelectionPlaceholder]}\)]}
        },
      RowSpacings->0,
      ColumnSpacings->0,
      GridDefaultElement:>ButtonBox[ "\\[Placeholder]"]]], NotebookDefault,
  CellMargins->{{Inherited, Inherited}, {5, Inherited}},
  Evaluatable->True,
  CellGroupingRules->"InputGrouping",
  PageBreakAbove->True,
  PageBreakWithin->False,
  GroupPageBreakWithin->False,
  CellLabelMargins->{{11, Inherited}, {Inherited, Inherited}},
  DefaultFormatType->DefaultInputFormatType,
  LineSpacing->{1.25, 0},
  AutoItalicWords->{},
  FormatType->InputForm,
  ScriptMinSize->9,
  ShowStringCharacters->True,
  NumberMarks->True,
  CounterIncrements->"Input",
  StyleMenuListing->None,
  FontFamily->"Courier",
  FontWeight->"Bold"]
}, FrontEndVersion->"Microsoft Windows 3.0", ScreenRectangle->{{0,
1024}, {0, 712}}, Editable->False, WindowToolbars->{},
PageWidth->820, WindowSize->{Fit, Fit}, WindowMargins->{{30,
Automatic}, {Automatic, 48}}, WindowFrame->"Palette",
WindowElements->{}, WindowFrameElements->"CloseBox",
WindowClickSelect->False,
ScrollingOptions->{"PagewiseScrolling"->True},
ShowCellBracket->False, CellMargins->{{0, 0}, {Inherited, 0}},
Active->True, CellOpen->True, ShowCellLabel->False,
ShowCellTags->False, ImageMargins->{{0, Inherited}, {Inherited,
0}}, Magnification->1 ]


(***********************************************************************
Cached data follows.  If you edit this Notebook file directly, not
using Mathematica, you must remove the line containing CacheID at
the top of the file.  The cache data will then be recreated when
you save this file from within Mathematica.
***********************************************************************)

(*CellTagsOutline CellTagsIndex->{} *)

(*CellTagsIndex CellTagsIndex->{} *)

(*NotebookFileOutline Notebook[{ Cell[1710, 49, 6415, 180, 330,
NotebookDefault,
  Evaluatable->True,
  CellGroupingRules->"InputGrouping",
  PageBreakAbove->True,
  PageBreakWithin->False,
  CounterIncrements->"Input"]
} ] *)




(***********************************************************************
End of Mathematica Notebook file.
***********************************************************************)

\end{verbatim}

\bigskip

\item[basic set of tonalities' palette]

\bigskip

\begin{verbatim}

(***********************************************************************

                    Mathematica-Compatible Notebook

This notebook can be used on any computer system with Mathematica
3.0, MathReader 3.0, or any compatible application. The data for
the notebook starts with the line of stars above.

To get the notebook into a Mathematica-compatible application, do
one of the following:

* Save the data starting with the line of stars above into a file
  with a name ending in .nb, then open the file inside the application;

* Copy the data starting with the line of stars above to the
  clipboard, then use the Paste menu command inside the application.

Data for notebooks contains only printable 7-bit ASCII and can be
sent directly in email or through ftp in text mode.  Newlines can
be CR, LF or CRLF (Unix, Macintosh or MS-DOS style).

NOTE: If you modify the data for this notebook not in a
Mathematica- compatible application, you must delete the line
below containing the word CacheID, otherwise
Mathematica-compatible applications may try to use invalid cache
data.

For more information on notebooks and Mathematica-compatible
applications, contact Wolfram Research:
  web: http://www.wolfram.com
  email: info@wolfram.com
  phone: +1-217-398-0700 (U.S.)

Notebook reader applications are available free of charge from
Wolfram Research.
***********************************************************************)

(*CacheID: 232*)


(*NotebookFileLineBreakTest NotebookFileLineBreakTest*)
(*NotebookOptionsPosition[      2997,         89]*)
(*NotebookOutlinePosition[      4053,        128]*) (*
CellTagsIndexPosition[      4009,        124]*)
(*WindowFrame->Palette*)



Notebook[{ Cell[BoxData[GridBox[{
        {
          ButtonBox[\(setofthemajortonalities[\[SelectionPlaceholder]]\)]},
        {
          ButtonBox[\(setoftheminortonalities[\[SelectionPlaceholder]]\)]},
        {

          ButtonBox[
            \(setoftheclassicaltonalities[\[SelectionPlaceholder]]\)]},
        {

          ButtonBox[
            \(setofthegregoriantonalities[\[SelectionPlaceholder]]\)]},
        {
          
ButtonBox[\(setofthemazzolatonalities[\[SelectionPlaceholder]]\)]},
        {
          ButtonBox[\(setofthejewishtonalities[\[SelectionPlaceholder]]\)]}
        },
      RowSpacings->0,
      ColumnSpacings->0,
      GridDefaultElement:>ButtonBox[ "\\[Placeholder]"]]], NotebookDefault,
  CellMargins->{{Inherited, Inherited}, {5, Inherited}},
  Evaluatable->True,
  CellGroupingRules->"InputGrouping",
  PageBreakAbove->True,
  PageBreakWithin->False,
  GroupPageBreakWithin->False,
  CellLabelMargins->{{11, Inherited}, {Inherited, Inherited}},
  DefaultFormatType->DefaultInputFormatType,
  LineSpacing->{1.25, 0},
  AutoItalicWords->{},
  FormatType->InputForm,
  ScriptMinSize->9,
  ShowStringCharacters->True,
  NumberMarks->True,
  CounterIncrements->"Input",
  StyleMenuListing->None,
  FontFamily->"Courier",
  FontWeight->"Bold"]
}, FrontEndVersion->"Microsoft Windows 3.0", ScreenRectangle->{{0,
1024}, {0, 712}}, Editable->False, WindowToolbars->{},
PageWidth->468, WindowSize->{Fit, Fit}, WindowMargins->{{30,
Automatic}, {Automatic, 48}}, WindowFrame->"Palette",
WindowElements->{}, WindowFrameElements->"CloseBox",
WindowClickSelect->False,
ScrollingOptions->{"PagewiseScrolling"->True},
ShowCellBracket->False, CellMargins->{{0, 0}, {Inherited, 0}},
Active->True, CellOpen->True, ShowCellLabel->False,
ShowCellTags->False, ImageMargins->{{0, Inherited}, {Inherited,
0}}, Magnification->1 ]


(***********************************************************************
Cached data follows.  If you edit this Notebook file directly, not
using Mathematica, you must remove the line containing CacheID at
the top of the file.  The cache data will then be recreated when
you save this file from within Mathematica.
***********************************************************************)

(*CellTagsOutline CellTagsIndex->{} *)

(*CellTagsIndex CellTagsIndex->{} *)

(*NotebookFileOutline Notebook[{ Cell[1710, 49, 1283, 38, 106,
NotebookDefault,
  Evaluatable->True,
  CellGroupingRules->"InputGrouping",
  PageBreakAbove->True,
  PageBreakWithin->False,
  CounterIncrements->"Input"]
} ] *)




(***********************************************************************
End of Mathematica Notebook file.
***********************************************************************)

\end{verbatim}

\bigskip

\item[play palette]

\bigskip

\begin{verbatim}

(***********************************************************************

                    Mathematica-Compatible Notebook

This notebook can be used on any computer system with Mathematica
3.0, MathReader 3.0, or any compatible application. The data for
the notebook starts with the line of stars above.

To get the notebook into a Mathematica-compatible application, do
one of the following:

* Save the data starting with the line of stars above into a file
  with a name ending in .nb, then open the file inside the application;

* Copy the data starting with the line of stars above to the
  clipboard, then use the Paste menu command inside the application.

Data for notebooks contains only printable 7-bit ASCII and can be
sent directly in email or through ftp in text mode.  Newlines can
be CR, LF or CRLF (Unix, Macintosh or MS-DOS style).

NOTE: If you modify the data for this notebook not in a
Mathematica- compatible application, you must delete the line
below containing the word CacheID, otherwise
Mathematica-compatible applications may try to use invalid cache
data.

For more information on notebooks and Mathematica-compatible
applications, contact Wolfram Research:
  web: http://www.wolfram.com
  email: info@wolfram.com
  phone: +1-217-398-0700 (U.S.)

Notebook reader applications are available free of charge from
Wolfram Research.
***********************************************************************)

(*CacheID: 232*)


(*NotebookFileLineBreakTest NotebookFileLineBreakTest*)
(*NotebookOptionsPosition[      2940,         89]*)
(*NotebookOutlinePosition[      3996,        128]*) (*
CellTagsIndexPosition[      3952,        124]*)
(*WindowFrame->Palette*)



Notebook[{ Cell[BoxData[GridBox[{
        {
          ButtonBox[\(playmonodic[\[SelectionPlaceholder]]\)]},
        {
          ButtonBox[\(playpolyphonic[\[SelectionPlaceholder]]\)]},
        {
          ButtonBox[\(playwwordasarpeggio[\[SelectionPlaceholder]]\)]},
        {

          ButtonBox[
            \(playwordaschord[\[SelectionPlaceholder],
              \[SelectionPlaceholder]]\)]},
        {

          ButtonBox[
            \(playharmonicword[\[SelectionPlaceholder],
              \[SelectionPlaceholder]]\)]}
        },
      RowSpacings->0,
      ColumnSpacings->0,
      GridDefaultElement:>ButtonBox[ "\\[Placeholder]"]]], NotebookDefault,
  CellMargins->{{Inherited, Inherited}, {5, Inherited}},
  Evaluatable->True,
  CellGroupingRules->"InputGrouping",
  PageBreakAbove->True,
  PageBreakWithin->False,
  GroupPageBreakWithin->False,
  CellLabelMargins->{{11, Inherited}, {Inherited, Inherited}},
  DefaultFormatType->DefaultInputFormatType,
  LineSpacing->{1.25, 0},
  AutoItalicWords->{},
  FormatType->InputForm,
  ScriptMinSize->9,
  ShowStringCharacters->True,
  NumberMarks->True,
  CounterIncrements->"Input",
  StyleMenuListing->None,
  FontFamily->"Courier",
  FontWeight->"Bold"]
}, FrontEndVersion->"Microsoft Windows 3.0", ScreenRectangle->{{0,
1024}, {0, 712}}, Editable->False, WindowToolbars->{},
PageWidth->413, WindowSize->{Fit, Fit}, WindowMargins->{{90,
Automatic}, {Automatic, 18}}, WindowFrame->"Palette",
WindowElements->{}, WindowFrameElements->"CloseBox",
WindowClickSelect->False,
ScrollingOptions->{"PagewiseScrolling"->True},
ShowCellBracket->False, CellMargins->{{0, 0}, {Inherited, 0}},
Active->True, CellOpen->True, ShowCellLabel->False,
ShowCellTags->False, ImageMargins->{{0, Inherited}, {Inherited,
0}}, Magnification->1 ]


(***********************************************************************
Cached data follows.  If you edit this Notebook file directly, not
using Mathematica, you must remove the line containing CacheID at
the top of the file.  The cache data will then be recreated when
you save this file from within Mathematica.
***********************************************************************)

(*CellTagsOutline CellTagsIndex->{} *)

(*CellTagsIndex CellTagsIndex->{} *)

(*NotebookFileOutline Notebook[{ Cell[1710, 49, 1226, 38, 90,
NotebookDefault,
  Evaluatable->True,
  CellGroupingRules->"InputGrouping",
  PageBreakAbove->True,
  PageBreakWithin->False,
  CounterIncrements->"Input"]
} ] *)




(***********************************************************************
End of Mathematica Notebook file.
***********************************************************************)

\end{verbatim}

\bigskip

\item[times'palette]

\bigskip

\begin{verbatim}

(***********************************************************************

                    Mathematica-Compatible Notebook

This notebook can be used on any computer system with Mathematica
3.0, MathReader 3.0, or any compatible application. The data for
the notebook starts with the line of stars above.

To get the notebook into a Mathematica-compatible application, do
one of the following:

* Save the data starting with the line of stars above into a file
  with a name ending in .nb, then open the file inside the application;

* Copy the data starting with the line of stars above to the
  clipboard, then use the Paste menu command inside the application.

Data for notebooks contains only printable 7-bit ASCII and can be
sent directly in email or through ftp in text mode.  Newlines can
be CR, LF or CRLF (Unix, Macintosh or MS-DOS style).

NOTE: If you modify the data for this notebook not in a
Mathematica- compatible application, you must delete the line
below containing the word CacheID, otherwise
Mathematica-compatible applications may try to use invalid cache
data.

For more information on notebooks and Mathematica-compatible
applications, contact Wolfram Research:
  web: http://www.wolfram.com
  email: info@wolfram.com
  phone: +1-217-398-0700 (U.S.)

Notebook reader applications are available free of charge from
Wolfram Research.
***********************************************************************)

(*CacheID: 232*)


(*NotebookFileLineBreakTest NotebookFileLineBreakTest*)
(*NotebookOptionsPosition[      2744,         87]*)
(*NotebookOutlinePosition[      3801,        126]*) (*
CellTagsIndexPosition[      3757,        122]*)
(*WindowFrame->Palette*)



Notebook[{ Cell[BoxData[GridBox[{
        {
          ButtonBox["semibreve"]},
        {
          ButtonBox["minim"]},
        {
          ButtonBox["crotchet"]},
        {
          ButtonBox["quaver"]},
        {
          ButtonBox["semiquaver"]},
        {
          ButtonBox["demisemiquaver"]},
        {
          ButtonBox["hemidemisemiquaver"]}
        },
      RowSpacings->0,
      ColumnSpacings->0,
      GridDefaultElement:>ButtonBox[ "\\[Placeholder]"]]], NotebookDefault,
  CellMargins->{{Inherited, Inherited}, {5, Inherited}},
  Evaluatable->True,
  CellGroupingRules->"InputGrouping",
  PageBreakAbove->True,
  PageBreakWithin->False,
  GroupPageBreakWithin->False,
  CellLabelMargins->{{11, Inherited}, {Inherited, Inherited}},
  DefaultFormatType->DefaultInputFormatType,
  LineSpacing->{1.25, 0},
  AutoItalicWords->{},
  FormatType->InputForm,
  ScriptMinSize->9,
  ShowStringCharacters->True,
  NumberMarks->True,
  CounterIncrements->"Input",
  StyleMenuListing->None,
  FontFamily->"Courier",
  FontWeight->"Bold"]
}, FrontEndVersion->"Microsoft Windows 3.0", ScreenRectangle->{{0,
1024}, {0, 712}}, Editable->False, WindowToolbars->{},
PageWidth->384, WindowSize->{Fit, Fit}, WindowMargins->{{100,
Automatic}, {Automatic, 33}}, WindowFrame->"Palette",
WindowElements->{}, WindowFrameElements->"CloseBox",
WindowClickSelect->False,
ScrollingOptions->{"PagewiseScrolling"->True},
ShowCellBracket->False, CellMargins->{{0, 0}, {Inherited, 0}},
Active->True, CellOpen->True, ShowCellLabel->False,
ShowCellTags->False, ImageMargins->{{0, Inherited}, {Inherited,
0}}, Magnification->1 ]


(***********************************************************************
Cached data follows.  If you edit this Notebook file directly, not
using Mathematica, you must remove the line containing CacheID at
the top of the file.  The cache data will then be recreated when
you save this file from within Mathematica.
***********************************************************************)

(*CellTagsOutline CellTagsIndex->{} *)

(*CellTagsIndex CellTagsIndex->{} *)

(*NotebookFileOutline Notebook[{ Cell[1710, 49, 1030, 36, 122,
NotebookDefault,
  Evaluatable->True,
  CellGroupingRules->"InputGrouping",
  PageBreakAbove->True,
  PageBreakWithin->False,
  CounterIncrements->"Input"]
} ] *)




(***********************************************************************
End of Mathematica Notebook file.
***********************************************************************)

\end{verbatim}

\bigskip

\item[Euler space's palette]

\bigskip

\begin{verbatim}

(***********************************************************************

                    Mathematica-Compatible Notebook

This notebook can be used on any computer system with Mathematica
3.0, MathReader 3.0, or any compatible application. The data for
the notebook starts with the line of stars above.

To get the notebook into a Mathematica-compatible application, do
one of the following:

* Save the data starting with the line of stars above into a file
  with a name ending in .nb, then open the file inside the application;

* Copy the data starting with the line of stars above to the
  clipboard, then use the Paste menu command inside the application.

Data for notebooks contains only printable 7-bit ASCII and can be
sent directly in email or through ftp in text mode.  Newlines can
be CR, LF or CRLF (Unix, Macintosh or MS-DOS style).

NOTE: If you modify the data for this notebook not in a
Mathematica- compatible application, you must delete the line
below containing the word CacheID, otherwise
Mathematica-compatible applications may try to use invalid cache
data.

For more information on notebooks and Mathematica-compatible
applications, contact Wolfram Research:
  web: http://www.wolfram.com
  email: info@wolfram.com
  phone: +1-217-398-0700 (U.S.)

Notebook reader applications are available free of charge from
Wolfram Research.
***********************************************************************)

(*CacheID: 232*)


(*NotebookFileLineBreakTest NotebookFileLineBreakTest*)
(*NotebookOptionsPosition[      3619,        107]*)
(*NotebookOutlinePosition[      4675,        146]*) (*
CellTagsIndexPosition[      4631,        142]*)
(*WindowFrame->Palette*)



Notebook[{ Cell[BoxData[GridBox[{
        {
          ButtonBox[\(eulercoordination[\[SelectionPlaceholder]]\)],
          ButtonBox["fifthpoint"]},
        {
          ButtonBox[\(FROMeulerpointTOnote[\[SelectionPlaceholder]]\)],
          ButtonBox["thirdpoint"]},
        {
          ButtonBox[\(FROMnoteTOpitch[\[SelectionPlaceholder]]\)],
          ButtonBox["canonicalnotesbasis"]},
        {
          ButtonBox[\(FROMpitchTOnote[\[SelectionPlaceholder]]\)],
          ButtonBox["canonicalintervalsbasis"]},
        {
          ButtonBox[\(FROMeulerpointTOpitch[\[SelectionPlaceholder]]\)],
          ButtonBox["fifthcomma"]},
        {

          ButtonBox[
            \(FROMnoteTOeulerpoint[\[SelectionPlaceholder],
              \[SelectionPlaceholder]]\)],
          ButtonBox["thirdcomma"]},
        {

          ButtonBox[
            \(FROMwordTOlistofeulerpoints[\[SelectionPlaceholder],
              \[SelectionPlaceholder]]\)],

          ButtonBox[
            \(gradussuavitatis[\[SelectionPlaceholder]  _Integer]\)]},
        {
          ButtonBox["octaveepoint"],

          ButtonBox[
            \(gradussuavitatis[\[SelectionPlaceholder]  _Rational]\)]}
        },
      RowSpacings->0,
      ColumnSpacings->0,
      GridDefaultElement:>ButtonBox[ "\\[Placeholder]"]]], NotebookDefault,
  CellMargins->{{Inherited, Inherited}, {5, Inherited}},
  Evaluatable->True,
  CellGroupingRules->"InputGrouping",
  PageBreakAbove->True,
  PageBreakWithin->False,
  GroupPageBreakWithin->False,
  CellLabelMargins->{{11, Inherited}, {Inherited, Inherited}},
  DefaultFormatType->DefaultInputFormatType,
  LineSpacing->{1.25, 0},
  AutoItalicWords->{},
  FormatType->InputForm,
  ScriptMinSize->9,
  ShowStringCharacters->True,
  NumberMarks->True,
  CounterIncrements->"Input",
  StyleMenuListing->None,
  FontFamily->"Courier",
  FontWeight->"Bold"]
}, FrontEndVersion->"Microsoft Windows 3.0", ScreenRectangle->{{0,
1024}, {0, 712}}, Editable->False, WindowToolbars->{},
PageWidth->694, WindowSize->{Fit, Fit}, WindowMargins->{{10,
Automatic}, {Automatic, 18}}, WindowFrame->"Palette",
WindowElements->{}, WindowFrameElements->"CloseBox",
WindowClickSelect->False,
ScrollingOptions->{"PagewiseScrolling"->True},
ShowCellBracket->False, CellMargins->{{0, 0}, {Inherited, 0}},
Active->True, CellOpen->True, ShowCellLabel->False,
ShowCellTags->False, ImageMargins->{{0, Inherited}, {Inherited,
0}}, Magnification->1 ]


(***********************************************************************
Cached data follows.  If you edit this Notebook file directly, not
using Mathematica, you must remove the line containing CacheID at
the top of the file.  The cache data will then be recreated when
you save this file from within Mathematica.
***********************************************************************)

(*CellTagsOutline CellTagsIndex->{} *)

(*CellTagsIndex CellTagsIndex->{} *)

(*NotebookFileOutline Notebook[{ Cell[1710, 49, 1905, 56, 142,
NotebookDefault,
  Evaluatable->True,
  CellGroupingRules->"InputGrouping",
  PageBreakAbove->True,
  PageBreakWithin->False,
  CounterIncrements->"Input"]
} ] *)




(***********************************************************************
End of Mathematica Notebook file.
***********************************************************************)

\end{verbatim}

\bigskip

\item[musical relativity's palette]

\bigskip

\begin{verbatim}

(***********************************************************************

                    Mathematica-Compatible Notebook

This notebook can be used on any computer system with Mathematica
3.0, MathReader 3.0, or any compatible application. The data for
the notebook starts with the line of stars above.

To get the notebook into a Mathematica-compatible application, do
one of the following:

* Save the data starting with the line of stars above into a file
  with a name ending in .nb, then open the file inside the application;

* Copy the data starting with the line of stars above to the
  clipboard, then use the Paste menu command inside the application.

Data for notebooks contains only printable 7-bit ASCII and can be
sent directly in email or through ftp in text mode.  Newlines can
be CR, LF or CRLF (Unix, Macintosh or MS-DOS style).

NOTE: If you modify the data for this notebook not in a
Mathematica- compatible application, you must delete the line
below containing the word CacheID, otherwise
Mathematica-compatible applications may try to use invalid cache
data.

For more information on notebooks and Mathematica-compatible
applications, contact Wolfram Research:
  web: http://www.wolfram.com
  email: info@wolfram.com
  phone: +1-217-398-0700 (U.S.)

Notebook reader applications are available free of charge from
Wolfram Research.
***********************************************************************)

(*CacheID: 232*)


(*NotebookFileLineBreakTest NotebookFileLineBreakTest*)
(*NotebookOptionsPosition[      3513,        107]*)
(*NotebookOutlinePosition[      4569,        146]*) (*
CellTagsIndexPosition[      4525,        142]*)
(*WindowFrame->Palette*)



Notebook[{ Cell[BoxData[GridBox[{
        {
          ButtonBox[\(symmetrytransformation[\[SelectionPlaceholder]]\)]},
        {

          ButtonBox[
            \(symmetricQ[\[SelectionPlaceholder], \[SelectionPlaceholder]]
              \)]},
        {

          ButtonBox[
            \(lawofresolutiononthetonicQ[\[SelectionPlaceholder],
              \[SelectionPlaceholder]]\)]},
        {
          
ButtonBox[\(specialtransformationrule[\[SelectionPlaceholder]]\)]},
        {

          ButtonBox[
            \(specialinvarianceQ[\[SelectionPlaceholder],
              \[SelectionPlaceholder]]\)]},
        {

          ButtonBox[
            
\(morethanspecialtransformationrule[\[SelectionPlaceholder]]\)]},
        {

          ButtonBox[
            \(morethanspecialinvarianceQ[\[SelectionPlaceholder],
              \[SelectionPlaceholder]]\)]},
        {
          
ButtonBox[\(generaltransformationrule[\[SelectionPlaceholder]]\)]},
        {
          ButtonBox["generalinvariance"]},
        {
          ButtonBox[\(\[SelectionPlaceholder]  _scale\)]}
        },
      RowSpacings->0,
      ColumnSpacings->0,
      GridDefaultElement:>ButtonBox[ "\\[Placeholder]"]]], NotebookDefault,
  CellMargins->{{Inherited, Inherited}, {5, Inherited}},
  Evaluatable->True,
  CellGroupingRules->"InputGrouping",
  PageBreakAbove->True,
  PageBreakWithin->False,
  GroupPageBreakWithin->False,
  CellLabelMargins->{{11, Inherited}, {Inherited, Inherited}},
  DefaultFormatType->DefaultInputFormatType,
  LineSpacing->{1.25, 0},
  AutoItalicWords->{},
  FormatType->InputForm,
  ScriptMinSize->9,
  ShowStringCharacters->True,
  NumberMarks->True,
  CounterIncrements->"Input",
  StyleMenuListing->None,
  FontFamily->"Courier",
  FontWeight->"Bold"]
}, FrontEndVersion->"Microsoft Windows 3.0", ScreenRectangle->{{0,
1024}, {0, 712}}, Editable->False, WindowToolbars->{},
PageWidth->511, WindowSize->{Fit, Fit}, WindowMargins->{{30,
Automatic}, {Automatic, 48}}, WindowFrame->"Palette",
WindowElements->{}, WindowFrameElements->"CloseBox",
WindowClickSelect->False,
ScrollingOptions->{"PagewiseScrolling"->True},
ShowCellBracket->False, CellMargins->{{0, 0}, {Inherited, 0}},
Active->True, CellOpen->True, ShowCellLabel->False,
ShowCellTags->False, ImageMargins->{{0, Inherited}, {Inherited,
0}}, Magnification->1 ]


(***********************************************************************
Cached data follows.  If you edit this Notebook file directly, not
using Mathematica, you must remove the line containing CacheID at
the top of the file.  The cache data will then be recreated when
you save this file from within Mathematica.
***********************************************************************)

(*CellTagsOutline CellTagsIndex->{} *)

(*CellTagsIndex CellTagsIndex->{} *)

(*NotebookFileOutline Notebook[{ Cell[1710, 49, 1799, 56, 172,
NotebookDefault,
  Evaluatable->True,
  CellGroupingRules->"InputGrouping",
  PageBreakAbove->True,
  PageBreakWithin->False,
  CounterIncrements->"Input"]
} ] *)




(***********************************************************************
End of Mathematica Notebook file.
***********************************************************************)

\end{verbatim}

\bigskip

\item[monodic palette]

\bigskip

\begin{verbatim}

(***********************************************************************

                    Mathematica-Compatible Notebook

This notebook can be used on any computer system with Mathematica
3.0, MathReader 3.0, or any compatible application. The data for
the notebook starts with the line of stars above.

To get the notebook into a Mathematica-compatible application, do
one of the following:

* Save the data starting with the line of stars above into a file
  with a name ending in .nb, then open the file inside the application;

* Copy the data starting with the line of stars above to the
  clipboard, then use the Paste menu command inside the application.

Data for notebooks contains only printable 7-bit ASCII and can be
sent directly in email or through ftp in text mode.  Newlines can
be CR, LF or CRLF (Unix, Macintosh or MS-DOS style).

NOTE: If you modify the data for this notebook not in a
Mathematica- compatible application, you must delete the line
below containing the word CacheID, otherwise
Mathematica-compatible applications may try to use invalid cache
data.

For more information on notebooks and Mathematica-compatible
applications, contact Wolfram Research:
  web: http://www.wolfram.com
  email: info@wolfram.com
  phone: +1-217-398-0700 (U.S.)

Notebook reader applications are available free of charge from
Wolfram Research.
***********************************************************************)

(*CacheID: 232*)


(*NotebookFileLineBreakTest NotebookFileLineBreakTest*)
(*NotebookOptionsPosition[      8418,        189]*)
(*NotebookOutlinePosition[      9474,        228]*) (*
CellTagsIndexPosition[      9430,        224]*)
(*WindowFrame->Palette*)



Notebook[{ Cell[BoxData[GridBox[{
        {
          ButtonBox[\(playmonodic[\[SelectionPlaceholder]]\)]},
        {

          ButtonBox[
            \({{\[SelectionPlaceholder], \[SelectionPlaceholder]}}\)]},
        {

          ButtonBox[
            \({{\[SelectionPlaceholder], \[SelectionPlaceholder]}, {
                \[SelectionPlaceholder], \[SelectionPlaceholder]}}\)]},
        {

          ButtonBox[
            \({{\[SelectionPlaceholder], \[SelectionPlaceholder]}, {
                \[SelectionPlaceholder], \[SelectionPlaceholder]}, {
                \[SelectionPlaceholder], \[SelectionPlaceholder]}}\)]},
        {

          ButtonBox[
            \({{\[SelectionPlaceholder], \[SelectionPlaceholder]}, {
                \[SelectionPlaceholder], \[SelectionPlaceholder]}, {
                \[SelectionPlaceholder], \[SelectionPlaceholder]}, {
                \[SelectionPlaceholder], \[SelectionPlaceholder]}}\)]},
        {

          ButtonBox[
            \({{\[SelectionPlaceholder], \[SelectionPlaceholder]}, {
                \[SelectionPlaceholder], \[SelectionPlaceholder]}, {
                \[SelectionPlaceholder], \[SelectionPlaceholder]}, {
                \[SelectionPlaceholder], \[SelectionPlaceholder]}, {
                \[SelectionPlaceholder], \[SelectionPlaceholder]}}\)]},
        {

          ButtonBox[
            \({{\[SelectionPlaceholder], \[SelectionPlaceholder]}, {
                \[SelectionPlaceholder], \[SelectionPlaceholder]}, {
                \[SelectionPlaceholder], \[SelectionPlaceholder]}, {
                \[SelectionPlaceholder], \[SelectionPlaceholder]}, {
                \[SelectionPlaceholder], \[SelectionPlaceholder]}, {
                \[SelectionPlaceholder], \[SelectionPlaceholder]}}\)]},
        {

          ButtonBox[
            \({{\[SelectionPlaceholder], \[SelectionPlaceholder]}, {
                \[SelectionPlaceholder], \[SelectionPlaceholder]}, {
                \[SelectionPlaceholder], \[SelectionPlaceholder]}, {
                \[SelectionPlaceholder], \[SelectionPlaceholder]}, {
                \[SelectionPlaceholder], \[SelectionPlaceholder]}, {
                \[SelectionPlaceholder], \[SelectionPlaceholder]}, {
                \[SelectionPlaceholder], \[SelectionPlaceholder]}}\)]},
        {

          ButtonBox[
            \({{\[SelectionPlaceholder], \[SelectionPlaceholder]}, {
                \[SelectionPlaceholder], \[SelectionPlaceholder]}, {
                \[SelectionPlaceholder], \[SelectionPlaceholder]}, {
                \[SelectionPlaceholder], \[SelectionPlaceholder]}, {
                \[SelectionPlaceholder], \[SelectionPlaceholder]}, {
                \[SelectionPlaceholder], \[SelectionPlaceholder]}, {
                \[SelectionPlaceholder], \[SelectionPlaceholder]}, {
                \[SelectionPlaceholder], \[SelectionPlaceholder]}}\)]},
        {

          ButtonBox[
            \({{\[SelectionPlaceholder], \[SelectionPlaceholder]}, {
                \[SelectionPlaceholder], \[SelectionPlaceholder]}, {
                \[SelectionPlaceholder], \[SelectionPlaceholder]}, {
                \[SelectionPlaceholder], \[SelectionPlaceholder]}, {
                \[SelectionPlaceholder], \[SelectionPlaceholder]}, {
                \[SelectionPlaceholder], \[SelectionPlaceholder]}, {
                \[SelectionPlaceholder], \[SelectionPlaceholder]}, {
                \[SelectionPlaceholder], \[SelectionPlaceholder]}, {
                \[SelectionPlaceholder], \[SelectionPlaceholder]}}\)]},
        {

          ButtonBox[
            \({{\[SelectionPlaceholder], \[SelectionPlaceholder]}, {
                \[SelectionPlaceholder], \[SelectionPlaceholder]}, {
                \[SelectionPlaceholder], \[SelectionPlaceholder]}, {
                \[SelectionPlaceholder], \[SelectionPlaceholder]}, {
                \[SelectionPlaceholder], \[SelectionPlaceholder]}, {
                \[SelectionPlaceholder], \[SelectionPlaceholder]}, {
                \[SelectionPlaceholder], \[SelectionPlaceholder]}, {
                \[SelectionPlaceholder], \[SelectionPlaceholder]}, {
                \[SelectionPlaceholder], \[SelectionPlaceholder]}, {
                \[SelectionPlaceholder], \[SelectionPlaceholder]}}\)]},
        {

          ButtonBox[
            \({{\[SelectionPlaceholder], \[SelectionPlaceholder]}, {
                \[SelectionPlaceholder], \[SelectionPlaceholder]}, {
                \[SelectionPlaceholder], \[SelectionPlaceholder]}, {
                \[SelectionPlaceholder], \[SelectionPlaceholder]}, {
                \[SelectionPlaceholder], \[SelectionPlaceholder]}, {
                \[SelectionPlaceholder], \[SelectionPlaceholder]}, {
                \[SelectionPlaceholder], \[SelectionPlaceholder]}, {
                \[SelectionPlaceholder], \[SelectionPlaceholder]}, {
                \[SelectionPlaceholder], \[SelectionPlaceholder]}, {
                \[SelectionPlaceholder], \[SelectionPlaceholder]}, {
                \[SelectionPlaceholder], \[SelectionPlaceholder]}}\)]},
        {

          ButtonBox[
            \({{\[SelectionPlaceholder], \[SelectionPlaceholder]}, {
                \[SelectionPlaceholder], \[SelectionPlaceholder]}, {
                \[SelectionPlaceholder], \[SelectionPlaceholder]}, {
                \[SelectionPlaceholder], \[SelectionPlaceholder]}, {
                \[SelectionPlaceholder], \[SelectionPlaceholder]}, {
                \[SelectionPlaceholder], \[SelectionPlaceholder]}, {
                \[SelectionPlaceholder], \[SelectionPlaceholder]}, {
                \[SelectionPlaceholder], \[SelectionPlaceholder]}, {
                \[SelectionPlaceholder], \[SelectionPlaceholder]}, {
                \[SelectionPlaceholder], \[SelectionPlaceholder]}, {
                \[SelectionPlaceholder], \[SelectionPlaceholder]}, {
                \[SelectionPlaceholder], \[SelectionPlaceholder]}}\)]}
        },
      RowSpacings->0,
      ColumnSpacings->0,
      GridDefaultElement:>ButtonBox[ "\\[Placeholder]"]]], NotebookDefault,
  CellMargins->{{Inherited, Inherited}, {5, Inherited}},
  Evaluatable->True,
  CellGroupingRules->"InputGrouping",
  PageBreakAbove->True,
  PageBreakWithin->False,
  GroupPageBreakWithin->False,
  CellLabelMargins->{{11, Inherited}, {Inherited, Inherited}},
  DefaultFormatType->DefaultInputFormatType,
  LineSpacing->{1.25, 0},
  AutoItalicWords->{},
  FormatType->InputForm,
  ScriptMinSize->9,
  ShowStringCharacters->True,
  NumberMarks->True,
  CounterIncrements->"Input",
  StyleMenuListing->None,
  FontFamily->"Courier",
  FontWeight->"Bold"]
}, FrontEndVersion->"Microsoft Windows 3.0", ScreenRectangle->{{0,
1024}, {0, 712}}, Editable->False, WindowToolbars->{},
PageWidth->866, WindowSize->{Fit, Fit}, WindowMargins->{{10,
Automatic}, {Automatic, 18}}, WindowFrame->"Palette",
WindowElements->{}, WindowFrameElements->"CloseBox",
WindowClickSelect->False,
ScrollingOptions->{"PagewiseScrolling"->True},
ShowCellBracket->False, CellMargins->{{0, 0}, {Inherited, 0}},
Active->True, CellOpen->True, ShowCellLabel->False,
ShowCellTags->False, ImageMargins->{{0, Inherited}, {Inherited,
0}}, Magnification->1 ]


(***********************************************************************
Cached data follows.  If you edit this Notebook file directly, not
using Mathematica, you must remove the line containing CacheID at
the top of the file.  The cache data will then be recreated when
you save this file from within Mathematica.
***********************************************************************)

(*CellTagsOutline CellTagsIndex->{} *)

(*CellTagsIndex CellTagsIndex->{} *)

(*NotebookFileOutline Notebook[{ Cell[1710, 49, 6704, 138, 218,
NotebookDefault,
  Evaluatable->True,
  CellGroupingRules->"InputGrouping",
  PageBreakAbove->True,
  PageBreakWithin->False,
  CounterIncrements->"Input"]
} ] *)




(***********************************************************************
End of Mathematica Notebook file.
***********************************************************************)

\end{verbatim}

\bigskip

\item[notes' palette]

\bigskip

\begin{verbatim}

(***********************************************************************

                    Mathematica-Compatible Notebook

This notebook can be used on any computer system with Mathematica
3.0, MathReader 3.0, or any compatible application. The data for
the notebook starts with the line of stars above.

To get the notebook into a Mathematica-compatible application, do
one of the following:

* Save the data starting with the line of stars above into a file
  with a name ending in .nb, then open the file inside the application;

* Copy the data starting with the line of stars above to the
  clipboard, then use the Paste menu command inside the application.

Data for notebooks contains only printable 7-bit ASCII and can be
sent directly in email or through ftp in text mode.  Newlines can
be CR, LF or CRLF (Unix, Macintosh or MS-DOS style).

NOTE: If you modify the data for this notebook not in a
Mathematica- compatible application, you must delete the line
below containing the word CacheID, otherwise
Mathematica-compatible applications may try to use invalid cache
data.

For more information on notebooks and Mathematica-compatible
applications, contact Wolfram Research:
  web: http://www.wolfram.com
  email: info@wolfram.com
  phone: +1-217-398-0700 (U.S.)

Notebook reader applications are available free of charge from
Wolfram Research.
***********************************************************************)

(*CacheID: 232*)


(*NotebookFileLineBreakTest NotebookFileLineBreakTest*)
(*NotebookOptionsPosition[      3654,        109]*)
(*NotebookOutlinePosition[      4709,        148]*) (*
CellTagsIndexPosition[      4665,        144]*)
(*WindowFrame->Palette*)



Notebook[{ Cell[BoxData[GridBox[{
        {
          ButtonBox[\(c[1]\)],
          ButtonBox[\(c[\[SelectionPlaceholder]]\)]},
        {
          ButtonBox[\(c\[Sharp][1]\)],
          ButtonBox[\(c\[Sharp][\[SelectionPlaceholder]]\)]},
        {
          ButtonBox[\(d[1]\)],
          ButtonBox[\(d[\[SelectionPlaceholder]]\)]},
        {
          ButtonBox[\(d\[Sharp][1]\)],
          ButtonBox[\(d\[Sharp][\[SelectionPlaceholder]]\)]},
        {
          ButtonBox[\(e[1]\)],
          ButtonBox[\(e[\[SelectionPlaceholder]]\)]},
        {
          ButtonBox[\(f[1]\)],
          ButtonBox[\(f[\[SelectionPlaceholder]]\)]},
        {
          ButtonBox[\(f\[Sharp][1]\)],
          ButtonBox[\(f\[Sharp][\[SelectionPlaceholder]]\)]},
        {
          ButtonBox[\(g[1]\)],
          ButtonBox[\(g[\[SelectionPlaceholder]]\)]},
        {
          ButtonBox[\(g\[Sharp][1]\)],
          ButtonBox[\(g\[Sharp][\[SelectionPlaceholder]]\)]},
        {
          ButtonBox[\(a[1]\)],
          ButtonBox[\(a[\[SelectionPlaceholder]]\)]},
        {
          ButtonBox[\(a\[Sharp][1]\)],
          ButtonBox[\(a\[Sharp][\[SelectionPlaceholder]]\)]},
        {
          ButtonBox[\(b[1]\)],
          ButtonBox[\(b[\[SelectionPlaceholder]]\)]}
        },
      RowSpacings->0,
      ColumnSpacings->0,
      GridDefaultElement:>ButtonBox[ "\\[Placeholder]"]]], NotebookDefault,
  CellMargins->{{Inherited, Inherited}, {5, Inherited}},
  Evaluatable->True,
  CellGroupingRules->"InputGrouping",
  PageBreakAbove->True,
  PageBreakWithin->False,
  GroupPageBreakWithin->False,
  CellLabelMargins->{{11, Inherited}, {Inherited, Inherited}},
  DefaultFormatType->DefaultInputFormatType,
  LineSpacing->{1.25, 0},
  AutoItalicWords->{},
  FormatType->InputForm,
  ScriptMinSize->9,
  ShowStringCharacters->True,
  NumberMarks->True,
  CounterIncrements->"Input",
  StyleMenuListing->None,
  FontFamily->"Courier",
  FontWeight->"Bold"]
}, FrontEndVersion->"Microsoft Windows 3.0", ScreenRectangle->{{0,
1024}, {0, 712}}, Editable->False, WindowToolbars->{},
PageWidth->334, WindowSize->{Fit, Fit}, WindowMargins->{{80,
Automatic}, {Automatic, 3}}, WindowFrame->"Palette",
WindowElements->{}, WindowFrameElements->"CloseBox",
WindowClickSelect->False,
ScrollingOptions->{"PagewiseScrolling"->True},
ShowCellBracket->False, CellMargins->{{0, 0}, {Inherited, 0}},
Active->True, CellOpen->True, ShowCellLabel->False,
ShowCellTags->False, ImageMargins->{{0, Inherited}, {Inherited,
0}}, Magnification->1 ]


(***********************************************************************
Cached data follows.  If you edit this Notebook file directly, not
using Mathematica, you must remove the line containing CacheID at
the top of the file.  The cache data will then be recreated when
you save this file from within Mathematica.
***********************************************************************)

(*CellTagsOutline CellTagsIndex->{} *)

(*CellTagsIndex CellTagsIndex->{} *)

(*NotebookFileOutline Notebook[{ Cell[1710, 49, 1940, 58, 202,
NotebookDefault,
  Evaluatable->True,
  CellGroupingRules->"InputGrouping",
  PageBreakAbove->True,
  PageBreakWithin->False,
  CounterIncrements->"Input"]
} ] *)




(***********************************************************************
End of Mathematica Notebook file.
***********************************************************************)

\end{verbatim}

    \item[pytagoric palette]

\begin{verbatim}

(***********************************************************************

                    Mathematica-Compatible Notebook

This notebook can be used on any computer system with Mathematica
3.0, MathReader 3.0, or any compatible application. The data for
the notebook starts with the line of stars above.

To get the notebook into a Mathematica-compatible application, do
one of the following:

* Save the data starting with the line of stars above into a file
  with a name ending in .nb, then open the file inside the application;

* Copy the data starting with the line of stars above to the
  clipboard, then use the Paste menu command inside the application.

Data for notebooks contains only printable 7-bit ASCII and can be
sent directly in email or through ftp in text mode.  Newlines can
be CR, LF or CRLF (Unix, Macintosh or MS-DOS style).

NOTE: If you modify the data for this notebook not in a
Mathematica- compatible application, you must delete the line
below containing the word CacheID, otherwise
Mathematica-compatible applications may try to use invalid cache
data.

For more information on notebooks and Mathematica-compatible
applications, contact Wolfram Research:
  web: http://www.wolfram.com
  email: info@wolfram.com
  phone: +1-217-398-0700 (U.S.)

Notebook reader applications are available free of charge from
Wolfram Research.
***********************************************************************)

(*CacheID: 232*)


(*NotebookFileLineBreakTest NotebookFileLineBreakTest*)
(*NotebookOptionsPosition[      3377,        105]*)
(*NotebookOutlinePosition[      4433,        144]*) (*
CellTagsIndexPosition[      4389,        140]*)
(*WindowFrame->Palette*)



Notebook[{ Cell[BoxData[GridBox[{
        {

          ButtonBox[
            \(pytagoricletter[\[SelectionPlaceholder],
              \[SelectionPlaceholder]]\)]},
        {

          ButtonBox[
            \(pytagoricalphabetuptofifthcycles[\[SelectionPlaceholder]]\)]},
        {

          ButtonBox[
            \(pytagoricscaleuptofifthcycles[\[SelectionPlaceholder]]\)]},
        {
          
ButtonBox[\(FROMpytagoricletterTOnote[\[SelectionPlaceholder]]\)]},
        {
          ButtonBox[\(FROMeulerpointTOpitch[fifthcomma]\)]},
        {

          ButtonBox[
            \(pytagoricword[\[SelectionPlaceholder], 
\[SelectionPlaceholder]]
              \)]},
        {

          ButtonBox[
            \(cycleraising[\[SelectionPlaceholder], \[SelectionPlaceholder]]
              \)]},
        {

          ButtonBox[
            \(cyclelowering[\[SelectionPlaceholder], 
\[SelectionPlaceholder]]
              \)]}
        },
      RowSpacings->0,
      ColumnSpacings->0,
      GridDefaultElement:>ButtonBox[ "\\[Placeholder]"]]], NotebookDefault,
  CellMargins->{{Inherited, Inherited}, {5, Inherited}},
  Evaluatable->True,
  CellGroupingRules->"InputGrouping",
  PageBreakAbove->True,
  PageBreakWithin->False,
  GroupPageBreakWithin->False,
  CellLabelMargins->{{11, Inherited}, {Inherited, Inherited}},
  DefaultFormatType->DefaultInputFormatType,
  LineSpacing->{1.25, 0},
  AutoItalicWords->{},
  FormatType->InputForm,
  ScriptMinSize->9,
  ShowStringCharacters->True,
  NumberMarks->True,
  CounterIncrements->"Input",
  StyleMenuListing->None,
  FontFamily->"Courier",
  FontWeight->"Bold"]
}, FrontEndVersion->"Microsoft Windows 3.0", ScreenRectangle->{{0,
1024}, {0, 712}}, Editable->False, WindowToolbars->{},
PageWidth->504, WindowSize->{Fit, Fit}, WindowMargins->{{30,
Automatic}, {Automatic, 48}}, WindowFrame->"Palette",
WindowElements->{}, WindowFrameElements->"CloseBox",
WindowClickSelect->False,
ScrollingOptions->{"PagewiseScrolling"->True},
ShowCellBracket->False, CellMargins->{{0, 0}, {Inherited, 0}},
Active->True, CellOpen->True, ShowCellLabel->False,
ShowCellTags->False, ImageMargins->{{0, Inherited}, {Inherited,
0}}, Magnification->1 ]


(***********************************************************************
Cached data follows.  If you edit this Notebook file directly, not
using Mathematica, you must remove the line containing CacheID at
the top of the file.  The cache data will then be recreated when
you save this file from within Mathematica.
***********************************************************************)

(*CellTagsOutline CellTagsIndex->{} *)

(*CellTagsIndex CellTagsIndex->{} *)

(*NotebookFileOutline Notebook[{ Cell[1710, 49, 1663, 54, 138,
NotebookDefault,
  Evaluatable->True,
  CellGroupingRules->"InputGrouping",
  PageBreakAbove->True,
  PageBreakWithin->False,
  CounterIncrements->"Input"]
} ] *)




(***********************************************************************
End of Mathematica Notebook file.
***********************************************************************)

\end{verbatim}

    \item[physical consonance's palette]

\begin{verbatim}

    (***********************************************************************

                    Mathematica-Compatible Notebook

This notebook can be used on any computer system with Mathematica
3.0, MathReader 3.0, or any compatible application. The data for
the notebook starts with the line of stars above.

To get the notebook into a Mathematica-compatible application, do
one of the following:

* Save the data starting with the line of stars above into a file
  with a name ending in .nb, then open the file inside the application;

* Copy the data starting with the line of stars above to the
  clipboard, then use the Paste menu command inside the application.

Data for notebooks contains only printable 7-bit ASCII and can be
sent directly in email or through ftp in text mode.  Newlines can
be CR, LF or CRLF (Unix, Macintosh or MS-DOS style).

NOTE: If you modify the data for this notebook not in a
Mathematica- compatible application, you must delete the line
below containing the word CacheID, otherwise
Mathematica-compatible applications may try to use invalid cache
data.

For more information on notebooks and Mathematica-compatible
applications, contact Wolfram Research:
  web: http://www.wolfram.com
  email: info@wolfram.com
  phone: +1-217-398-0700 (U.S.)

Notebook reader applications are available free of charge from
Wolfram Research.
***********************************************************************)

(*CacheID: 232*)


(*NotebookFileLineBreakTest NotebookFileLineBreakTest*)
(*NotebookOptionsPosition[      3080,         93]*)
(*NotebookOutlinePosition[      4136,        132]*) (*
CellTagsIndexPosition[      4092,        128]*)
(*WindowFrame->Palette*)



Notebook[{ Cell[BoxData[GridBox[{
        {

          ButtonBox[
            \(physicalindexofconsonanceofsounds[\[SelectionPlaceholder],
              \[SelectionPlaceholder]]\)]},
        {

          ButtonBox[
            \(physicalindexofconsonanceofpulsation[\[SelectionPlaceholder],
              \[SelectionPlaceholder], \[SelectionPlaceholder]]\)]},
        {

          ButtonBox[
            \(commensurabilityQ[\[SelectionPlaceholder],
              \[SelectionPlaceholder], \[SelectionPlaceholder]\ ]\)]},
        {

          ButtonBox[
            \(empiricalsimplicitymeasure[\[SelectionPlaceholder]  _Rational]
              \)]}
        },
      RowSpacings->0,
      ColumnSpacings->0,
      GridDefaultElement:>ButtonBox[ "\\[Placeholder]"]]], NotebookDefault,
  CellMargins->{{Inherited, Inherited}, {5, Inherited}},
  Evaluatable->True,
  CellGroupingRules->"InputGrouping",
  PageBreakAbove->True,
  PageBreakWithin->False,
  GroupPageBreakWithin->False,
  CellLabelMargins->{{11, Inherited}, {Inherited, Inherited}},
  DefaultFormatType->DefaultInputFormatType,
  LineSpacing->{1.25, 0},
  AutoItalicWords->{},
  FormatType->InputForm,
  ScriptMinSize->9,
  ShowStringCharacters->True,
  NumberMarks->True,
  CounterIncrements->"Input",
  StyleMenuListing->None,
  FontFamily->"Courier",
  FontWeight->"Bold"]
}, FrontEndVersion->"Microsoft Windows 3.0", ScreenRectangle->{{0,
1024}, {0, 712}}, Editable->False, WindowToolbars->{},
PageWidth->569, WindowSize->{324, 71}, WindowMargins->{{Automatic,
119}, {Automatic, 13}}, WindowFrame->"Palette",
WindowElements->{}, WindowFrameElements->"CloseBox",
WindowClickSelect->False,
ScrollingOptions->{"PagewiseScrolling"->True},
ShowCellBracket->False, CellMargins->{{0, 0}, {Inherited, 0}},
Active->True, CellOpen->True, ShowCellLabel->False,
ShowCellTags->False, ImageMargins->{{0, Inherited}, {Inherited,
0}}, Magnification->1 ]


(***********************************************************************
Cached data follows.  If you edit this Notebook file directly, not
using Mathematica, you must remove the line containing CacheID at
the top of the file.  The cache data will then be recreated when
you save this file from within Mathematica.
***********************************************************************)

(*CellTagsOutline CellTagsIndex->{} *)

(*CellTagsIndex CellTagsIndex->{} *)

(*NotebookFileOutline Notebook[{ Cell[1710, 49, 1366, 42, 76,
NotebookDefault,
  Evaluatable->True,
  CellGroupingRules->"InputGrouping",
  PageBreakAbove->True,
  PageBreakWithin->False,
  CounterIncrements->"Input"]
} ] *)




(***********************************************************************
End of Mathematica Notebook file.
***********************************************************************)

\end{verbatim}

    \item[form's conversion's palette]

\begin{verbatim}

(***********************************************************************

                    Mathematica-Compatible Notebook

This notebook can be used on any computer system with Mathematica
3.0, MathReader 3.0, or any compatible application. The data for
the notebook starts with the line of stars above.

To get the notebook into a Mathematica-compatible application, do
one of the following:

* Save the data starting with the line of stars above into a file
  with a name ending in .nb, then open the file inside the application;

* Copy the data starting with the line of stars above to the
  clipboard, then use the Paste menu command inside the application.

Data for notebooks contains only printable 7-bit ASCII and can be
sent directly in email or through ftp in text mode.  Newlines can
be CR, LF or CRLF (Unix, Macintosh or MS-DOS style).

NOTE: If you modify the data for this notebook not in a
Mathematica- compatible application, you must delete the line
below containing the word CacheID, otherwise
Mathematica-compatible applications may try to use invalid cache
data.

For more information on notebooks and Mathematica-compatible
applications, contact Wolfram Research:
  web: http://www.wolfram.com
  email: info@wolfram.com
  phone: +1-217-398-0700 (U.S.)

Notebook reader applications are available free of charge from
Wolfram Research.
***********************************************************************)

(*CacheID: 232*)


(*NotebookFileLineBreakTest NotebookFileLineBreakTest*)
(*NotebookOptionsPosition[      3856,        112]*)
(*NotebookOutlinePosition[      4914,        151]*) (*
CellTagsIndexPosition[      4870,        147]*)
(*WindowFrame->Palette*)



Notebook[{ Cell[BoxData[GridBox[{
        {
          ButtonBox[\(FROMletterTOnote[\[SelectionPlaceholder]]\)]},
        {
          ButtonBox[\(FROMwordTOscale[\[SelectionPlaceholder]]\)]},
        {
          ButtonBox[\(FROMscaleTOpiece[\[SelectionPlaceholder]]\)]},
        {
          ButtonBox[\(FROMpieceTOscale[\[SelectionPlaceholder]]\)]},
        {

          ButtonBox[
            \(FROMharmonicwordTOphysicalchord[\[SelectionPlaceholder]]\)]},
        {
          ButtonBox[\(FROMlistTOcolumnvector[\[SelectionPlaceholder]]\)]},
        {
          ButtonBox[\(FROMcolumnvectorTOlist[\[SelectionPlaceholder]]\)]},
        {
          ButtonBox[\(FROMnoteTOpulsation[\[SelectionPlaceholder]]\)]},
        {
          ButtonBox[\(FROMpulsationTOnote[\[SelectionPlaceholder]]\)]},
        {

          ButtonBox[
            \(FROMpulsationTOsound[\[SelectionPlaceholder],
              \[SelectionPlaceholder]]\)]},
        {
          ButtonBox[\(FROMsoundTOpulsation[\[SelectionPlaceholder]]\)]},
        {
          ButtonBox[\(FROMnoteTOpitch[\[SelectionPlaceholder]]\)]},
        {
          ButtonBox[\(FROMpitchTOnote[\[SelectionPlaceholder]]\)]},
        {
          ButtonBox[\(FROMeulerpointTOpitch[\[SelectionPlaceholder]]\)]},
        {
          ButtonBox[\(FROMnoteTOeulerpoint[\[SelectionPlaceholder]]\)]},
        {

          ButtonBox[
            \(FROMwordTOlistofeulerpoints[\[SelectionPlaceholder]]\)]}
        },
      RowSpacings->0,
      ColumnSpacings->0,
      GridDefaultElement:>ButtonBox[ "\\[Placeholder]"]]], NotebookDefault,
  CellMargins->{{Inherited, Inherited}, {5, Inherited}},
  Evaluatable->True,
  CellGroupingRules->"InputGrouping",
  PageBreakAbove->True,
  PageBreakWithin->False,
  GroupPageBreakWithin->False,
  CellLabelMargins->{{11, Inherited}, {Inherited, Inherited}},
  DefaultFormatType->DefaultInputFormatType,
  LineSpacing->{1.25, 0},
  AutoItalicWords->{},
  FormatType->InputForm,
  ScriptMinSize->9,
  ShowStringCharacters->True,
  NumberMarks->True,
  CounterIncrements->"Input",
  StyleMenuListing->None,
  FontFamily->"Courier",
  FontWeight->"Bold"]
}, FrontEndVersion->"Microsoft Windows 3.0", ScreenRectangle->{{0,
1024}, {0, 712}}, Editable->False, WindowToolbars->{},
PageWidth->497, WindowSize->{252, 261},
WindowMargins->{{Automatic, 247}, {Automatic, 117}},
WindowFrame->"Palette", WindowElements->{},
WindowFrameElements->"CloseBox", WindowClickSelect->False,
ScrollingOptions->{"PagewiseScrolling"->True},
ShowCellBracket->False, CellMargins->{{0, 0}, {Inherited, 0}},
Active->True, CellOpen->True, ShowCellLabel->False,
ShowCellTags->False, ImageMargins->{{0, Inherited}, {Inherited,
0}}, Magnification->1 ]


(***********************************************************************
Cached data follows.  If you edit this Notebook file directly, not
using Mathematica, you must remove the line containing CacheID at
the top of the file.  The cache data will then be recreated when
you save this file from within Mathematica.
***********************************************************************)

(*CellTagsOutline CellTagsIndex->{} *)

(*CellTagsIndex CellTagsIndex->{} *)

(*NotebookFileOutline Notebook[{ Cell[1710, 49, 2142, 61, 266,
NotebookDefault,
  Evaluatable->True,
  CellGroupingRules->"InputGrouping",
  PageBreakAbove->True,
  PageBreakWithin->False,
  CounterIncrements->"Input"]
} ] *)




(***********************************************************************
End of Mathematica Notebook file.
***********************************************************************)

\end{verbatim}

\end{description}

\end{document}